\documentclass[a4paper,11pt,reqno]{amsart}
\usepackage[centering,includeheadfoot,margin=2.5 cm]{geometry}
\usepackage{graphics,enumitem,epsfig,textcomp}
\usepackage{amsfonts,amsmath,amssymb,euscript,color,mathrsfs}
\usepackage[utf8]{inputenc}	

\usepackage[symbol]{footmisc}

\usepackage
[
pdfauthor   = {Lisa Maria Kreusser},
pdftitle    = {Title},
pdfsubject  = {},
pdfkeywords = {},
bookmarks=true,
bookmarksopen=true,
pagebackref=false,
hyperindex=true,
colorlinks=true,
linkcolor=blue,
citecolor=blue,
filecolor=blue,
urlcolor=blue,
]
{hyperref}					
\usepackage{textcomp}	
\usepackage{color}			
\usepackage{graphicx}		

\usepackage[caption=false,subrefformat=parens,labelformat=parens]{subfig}

\newcommand\N{\mathbb{N}}
\newcommand\Z{\mathbb{Z}}
\newcommand\R{\mathbb{R}}

\newcommand\C{\mathbb{C}}

\newcommand\bl{\left(}
\newcommand\br{\right)}
\newcommand{\length}{\ell}

\newcommand*\di{\mathop{}\!\mathrm{d}}

\renewcommand\epsilon{\varepsilon}
\renewcommand\theta{\vartheta}

\newtheoremstyle{mytheoremstyle} 
{6pt}                    
{6pt}                    
{\itshape}                   
{}						
{\bf}                   
{}                          
{.5em}                       
{}  
\theoremstyle{mytheoremstyle}
\newtheorem{satz}{Satz}[section]
\newtheorem{lemma}[satz]{Lemma}
\newtheorem{proposition}[satz]{Proposition}
\newtheorem{corollary}[satz]{Corollary}
\newtheorem{theorem}[satz]{Theorem}
\newtheorem{assumption}[satz]{Assumption}
\newtheorem{remark}[satz]{Remark}

\newtheoremstyle{mytdefintionstyle} 
{6pt}                    
{6pt}                    
{\rm}                   
{}						
{\bf}                   
{}                          
{.5em}                       
{}  
\theoremstyle{mytdefintionstyle}

\graphicspath{{./Numerics/}}

\numberwithin{equation}{section}
\usepackage{cleveref}

\title{Stability analysis of line patterns of an anisotropic interaction model}

\begin{document}
	
	\maketitle

	\centerline{
		{\large Jos\'{e} A. Carrillo}\footnote{Department of Mathematics, Imperial College London, London SW7 2AZ, United Kingdom; 
			{\it carrillo@imperial.ac.uk}}\quad
		{\large Bertram D\"uring}\footnote{Department of Mathematics, University of Sussex, Pevensey II, Brighton BN1 9QH, United Kingdom;
			{\it b.during@sussex.ac.uk}}\quad
		{\large Lisa Maria Kreusser}\footnote{Department of Applied Mathematics and Theoretical Physics (DAMTP), University of Cambridge, Wilberforce Road, Cambridge CB3 0WA, United Kingdom;
			{\it L.M.Kreusser@damtp.cam.ac.uk}}\quad
		{\large Carola-Bibiane Sch\"{o}nlieb}\footnote{Department of Applied Mathematics and Theoretical Physics (DAMTP), University of Cambridge, Wilberforce Road, Cambridge CB3 0WA, United Kingdom; {\it C.B.Schoenlieb@damtp.cam.ac.uk}}	
	}
	\vskip 10mm

	\noindent{\bf Abstract.}
Motivated by the formation of fingerprint patterns we consider a class of interacting particle models with anisotropic, repulsive-attractive interaction forces whose orientations depend on an underlying tensor field. This class of models can be regarded as a generalization of a gradient flow of a nonlocal interaction potential which has a local repulsion and a long-range attraction structure. In addition, the underlying tensor field introduces an anisotropy leading to complex patterns which do not occur in isotropic models. Central to this pattern formation are straight line patterns. For a given spatially homogeneous tensor field, we show that there exists a preferred direction of straight lines, i.e.\ straight vertical lines can be stable for sufficiently many particles, while many other rotations of the straight lines are unstable steady states, both for a sufficiently large number of particles and in the continuum limit. For straight vertical lines we consider specific force coefficients for the stability analysis of steady states, show that stability can be achieved for exponentially decaying force coefficients for a sufficiently large number of particles and relate these results to the K\"ucken-Champod model for simulating fingerprint patterns. The mathematical analysis of the steady states is completed with numerical results.

	\vskip 3mm
	
	
	
	\vskip 7mm

\noindent{\bf Keywords.}
	Aggregation, swarming, pattern formation, dynamical systems.

\noindent{\bf AMS.}
	35B36, 35Q92, 70F10, 70F45, 82C22

\section{Introduction}

Mathematical models for biological aggregation describing the collective behaviour of large numbers of individuals have given us many tools to understand pattern formation in nature. Typical examples include models for explaining the complex phenomena observed in swarms of insects, flocks of birds, schools of fish or colonies of bacteria see for instance
\cite{Ballerini,swarmequilibria, Birnir2007,selforganization, cuckersmale,Cavagna,fishbehavior,bacterialcolonies,selfpropelled_particles,migratorylocusts, nonlocal_swarm, scienceanialaggregation,order_of_chaos}. Some continuum models have been derived from individual based descriptions \cite{Blanchet2008,Blanchet2006,Boi2000163,burger2007aggregation,swarmdynamics,SwarmingPatterns,TBL,PredictingPatternFormation},  see also the reviews \cite{CFTV,review}, leading to an understanding of the stability of patterns at different levels \cite{secondordermodel,Bertozzi2015,CHM2,Carrillo2014NonlinFlock,stabilityringpatterns}.

A key feature of many of these models is that the  communication between individuals takes place at different scales, i.e.\ each individual can interact not only with its neighbours but also with individuals further away. This can be described by short- and long-range interactions  \cite{swarmequilibria,migratorylocusts,nonlocal_swarm}. In most models the interactions are assumed to be \textit{isotropic} for simplicity. However, pattern formation in nature is usually \textit{anisotropic} \cite{patternsinnature}. Motivated by the simulation of fingerprint patterns we consider a class of interacting particle models with \textit{anisotropic} interaction forces in this paper. In particular, these anisotropic interaction models capture  important swarming behaviours, neglected in the simplified isotropic interaction model, such as anisotropic steady states.

The simplest form of isotropic interaction models is based on  radial interaction potentials  \cite{nonlocalinteraction}. In this case one can consider the stationary points of the $N$-particle interaction energy
\begin{align*}
E(x_1,\ldots,x_N)=\frac{1}{2N^2}\sum_{\substack{j,k=1\\k\neq j}}^N W\bl x_j-x_k\br.
\end{align*}
Here, $W(d)=\overline{W}(|d|)$ denotes the radially symmetric interaction potential and $x_j=x_j(t)\in \R^n$ for $j=1,\ldots,N$ are the positions of the particles at time $t\geq 0$  \cite{Bertozzi2015, stabilityringpatterns}. One can easily show that the associated gradient flow reads:
\begin{align}\label{eq:standardmodel}
\frac{\di x_j}{\di t}
=\frac{1}{N}\sum_{\substack{k=1\\k\neq j}}^N F(x_j-x_k),
\end{align}
where  $F(x_j-x_k)$ is a conservative force, aligned along the distance vector $x_j-x_k$ with ${F(d)=-\nabla W(d)}$. In many biological applications the number of interacting particles is large and one may consider the underlying continuum formulation of \eqref{eq:standardmodel}, which is known as the aggregation equation  \cite{bertozzi2009,  Bertozzi2015,stabilityringpatterns} and of the form
\begin{align}\label{eq:standardmodelmacroscopic}
\rho_t+\nabla\cdot\bl \rho u\br=0,\qquad u=-\nabla W\ast \rho,
\end{align} 
where  $u=u(t,x)$ is the macroscopic velocity field and $\rho=\rho(t,x)$ denotes the density of particles at location $x\in \R^n$ at time $t>0$.
The aggregation equation \eqref{eq:standardmodelmacroscopic} has  been studied extensively recently, mainly in terms of its gradient flow structure \cite{gradientflows,Carrillo:2003,  Carrillo2006,Li2004, opac-b1122739}, the blow-up dynamics for fully 	attractive potentials \cite{bertozzi2009, bertozzi2012,carrillo2011,  Carrillo2016304}, and the rich variety of steady states for repulsive-attractive potentials \cite{balague_preprint, nonlocalinteraction, Confinement, swarmequilibria, bertozzi2012, Canizo2015, Carrillo2016,Carrillo2012550, Carrillo2012306,CH,CHM2,Carrillo2014NonlinFlock,Fellner2010, Fellner20111436,nonlocalInteractionEquation,vonBrecht2012, PredictingPatternFormation}. 

In biological applications, the interactions determined by the force $F$, or equivalently the interaction potential $W$, are usually described by  short-range repulsion,  preventing collisions between the individuals, as well as  long-range attraction,  keeping the swarm cohesive \cite{Mogilner2003, Okubo}. In this case, the associated radially symmetric potentials $\overline{W}$  first decrease  and then increase as a function of the radius. Due to the repulsive forces these potentials lead to possibly more  steady states than the purely attractive potentials. In particular, these repulsive-attractive potentials can be considered as a minimal model for pattern formation in large systems of individuals \cite{nonlocalinteraction,review} and the references therein.

Pattern formation in multiple dimensions is studied in \cite{Bertozzi2015, stabilityringpatterns, vonBrecht2012, PredictingPatternFormation,CHM2} for repulsive-attractive potentials. The instabilities of the sphere and ring solutions are studied in \cite{Bertozzi2015, vonBrecht2012, PredictingPatternFormation}. The linear stability of ring equilibria is analysed and  conditions on the potential are derived to classify the different instabilities. A numerical study of the $N$-particle interaction model for specific repulsion-attraction potentials is also performed in \cite{Bertozzi2015,stabilityringpatterns} leading to a wide range of radially symmetric patterns such as rings, annuli and uniform circular patches, as well as more complex patterns. Based on this analysis the  stability of flock solutions and mill rings in the associated second-order model can be studied, see  \cite{secondordermodel} and \cite{Carrillo2014NonlinFlock} for the linear and  nonlinear stability of flocks, respectively. 

In this work, we consider a generalization of the particle model \eqref{eq:standardmodel} by introducing an anisotropy given by a tensor field $T$. This leads to an extended particle model of the form
\begin{align}\label{eq:particlemodel}
\frac{\di x_j}{\di t}=\frac{1}{N}\sum_{\substack{k=1\\k\neq j}}^N F(x_j-x_k,T(x_j))
\end{align}
where we prescribe initial data $x_j(0)=x_j^{in},~j=1,\ldots,N$, for given scalars $x_j^{in},~j=1,\ldots,N$.  A special instance of this model  has been  introduced in \cite{Merkel} for simulating fingerprint patterns. The particle model in its general form \eqref{eq:particlemodel} has been studied in \cite{patternformationanisotropicmodel,anisotropicfingerprint}.
Here, the positions of each of the $N$ particles at time $t$ are denoted by  
$x_j=x_j(t)\in\mathbb{R}^2$, $j=1,\ldots,N,$
and   $F(x_j-x_k,T(x_j))$  denotes the total force that particle $k$ exerts on particle $j$ subject to an underlying stress tensor field $T(x_j)$ at $x_j$, given by 
\begin{align}\label{eq:tensorfield}
T(x):=\chi s(x)\otimes s(x) +l(x)\otimes l(x)\in\R^{2,2}
\end{align}
for orthonormal  vector fields $s=s(x)$ and $l=l(x)\in\R^2$ and $\chi\in[0,1]$. Here, the outer product $v\otimes w$ for two vectors $v,w\in\R^2$ equals the matrix multiplication $vw^T$ and results in a matrix of size $\R^{2,2}$. The parameter $\chi$ introduces an anisotropy in the direction $s$  in the definition of the tensor field.

For  repulsive forces along $s$ and short-range repulsive, long-range attractive forces along $l$ the numerical simulations in \cite{patternformationanisotropicmodel} suggest that straight vertical line patterns formed by the interacting particles at positions $x_j$ are stable for a certain spatially homogeneous tensor field, specified later. In this paper, we want to rigorously study this empirical observation by providing a linear stability analysis of such patterns where particles distribute equidistantly along straight lines. 

The stability analysis of steady states of the particle model \eqref{eq:particlemodel} is important for understanding the robustness of the patterns that arise from applying \eqref{eq:particlemodel} for numerical simulation, for instance as for its originally intended application to fingerprint simulation in \cite{Merkel}. Indeed, in what follows,  we will show  that for spatially homogeneous tensor fields $T$ the solution formed by a number of vertical straight lines (referred to as ridges) is a stationary solution, whereas ridge bifurcations, i.e.\ a single ridge dividing into two ridges as typically appearing in fingerprint patterns, is not.

The aim of this paper is to prove that  sufficiently large numbers of 
particles distributed equidistantly along straight vertical lines  are stable steady states to the particle model \eqref{eq:particlemodel} for short-range repulsive, long-range attractive forces along $l$ and  repulsive forces along $s$. All other rotations of straight lines are unstable steady states for this choice of force coefficients for a sufficiently large number of particles and for the continuum limit. We focus on this very simple class of steady states as a first step towards understanding stable formations that can be achieved by model \eqref{eq:particlemodel}. Note that  the continuum straight line is a steady state of the associated continuum model
\begin{align*}
\begin{split}
\partial_t \rho(t,x)+\nabla_x\cdot \left[ \rho(t,x)\bl F\bl\cdot,T(x)\br \ast \rho(t,\cdot)\br\bl x\br\right]=0\qquad \text{in }\R_+\times \R^2,
\end{split}
\end{align*} 
see \cite{patternformationanisotropicmodel}, but its asymptotic stability cannot be concluded from the linear stability analysis for finitely many particles.

The paper is organized as follows. In Section \ref{sec:modeldescription} we describe a  general formulation of an anisotropic interaction model, based on the model proposed  by K\"ucken and Champod \cite{Merkel}. Section \ref{sec:stationarystates} is devoted to a high-wave number stability analysis of line patterns for the continuum limit $N\to\infty$, including vertical, horizontal and rotated straight lines for spatially homogeneous tensor fields. Due to the instability of arbitrary rotations except for vertical straight lines for the considered tensor field we focus on the stability analysis of straight vertical lines for particular forces for any $N\in\N$ in Section \ref{sec:verticallines}. Section \ref{sec:numerics} illustrates the form of the steady states in case the derived stability conditions are not satisfied.


\section{Description of the model}\label{sec:modeldescription}
In this section, we describe a general formulation of the anisotropic microscopic model \eqref{eq:particlemodel} and relate it to the K\"{u}cken-Champod particle model \cite{Merkel}. 
K\"{u}cken and Champod consider the particle model \eqref{eq:particlemodel} where the total force $F$  is given by
\begin{align}\label{eq:totalforce}
F(d(x_j,x_k),T(x_j))=F_A(d(x_j,x_k),T(x_j))+F_R(d(x_j,x_k))
\end{align} 
for the distance vector $d(x_j,x_k)=x_j-x_k\in\R^2$. Here, $F_R$ denotes the repulsion force that particle $k$ exerts on particle $j$ and $F_A$ is the attraction force particle $k$ exerts on particle $j$. The repulsion and attraction forces are of the form
\begin{align}\label{eq:repulsionforce}
F_R(d=d(x_j,x_k))=f_R(|d|)d
\end{align}
and
\begin{align}\label{eq:attractionforce}
F_A(d=d(x_j,x_k),T(x_j))=f_A(|d|)T(x_j)d,
\end{align}
respectively, with coefficient functions $f_R$ and $f_A$, where, again, $d=d(x_j,x_k)=x_j-x_k\in\R^2$. Note that the repulsion and attraction force coefficients $f_R,f_A$ are radially symmetric.
The direction of the interaction forces is determined by the parameter $\chi\in[0,1]$ in the definition of  $T$ in \eqref{eq:tensorfield}. 
Motivated by plugging \eqref{eq:tensorfield} into the definition of the total force \eqref{eq:totalforce},  we consider a more general form of the total force, given by
\begin{align}\label{eq:totalforcenew}
F(d=d(x_j,x_k),T(x_j))=f_s(|d|)(s(x_j)\cdot d)s(x_j)+f_l(|d|)(l(x_j)\cdot d)l(x_j),
\end{align}
where the total force is decomposed into  forces along the direction $s$ and along the direction $l$.
In particular, the force coefficients in the K\"ucken-Champod model \eqref{eq:particlemodel} with repulsive and attractive forces $F_R$ and $F_A$ in \eqref{eq:repulsionforce} and \eqref{eq:attractionforce}, respectively, can be recovered for 
\begin{align*}
f_l(|d|)=f_A(|d|)+f_R(|d|) \quad\text{and}\quad f_s(|d|)=\chi f_A(|d|)+f_R(|d|).
\end{align*}

Since a steady state of the particle model \eqref{eq:particlemodel}  for any spatially homogeneous tensor field $\tilde{T}$ can be regarded as a coordinate transform of the steady state of the particle model \eqref{eq:particlemodel} for the  tensor field $T$, see \cite{patternformationanisotropicmodel} for details, we restrict ourselves to the study of  steady states for the spatially homogeneous tensor field $T$ given by the orthonormal vectors $s=(0,1)$ and $l=(1,0)$, i.e.
\begin{align}\label{eq:tensorfieldhom}
T=\begin{pmatrix}
1& 0\\ 0 &\chi
\end{pmatrix}.
\end{align}
The total force in the K\"ucken-Champod model \eqref{eq:totalforce} and the generalised total force \eqref{eq:totalforcenew}   reduce to 
\begin{align}\label{eq:totalforcenewhom}
F(d)=\begin{pmatrix}
\bl f_A(|d|)+f_R(|d|)\br d_1\\\ \bl\chi f_A(|d|)+f_R(|d|)\br d_2
\end{pmatrix}
\end{align} 
and 
\begin{align}\label{eq:totalforcehom}
F(d)=\begin{pmatrix}
f_l(|d|)d_1\\\ f_s(|d|) d_2
\end{pmatrix}\quad \text{for~} d=(d_1,d_2)\in \R^2,
\end{align} 
respectively, for the spatially homogeneous tensor field $T$ in \eqref{eq:tensorfieldhom}.

In the sequel, we consider the particle model \eqref{eq:particlemodel} on the torus $\mathbb{T}^2$, or equivalently, on the unit square $[0,1]^2$ with periodic boundary conditions. This can be achieved by considering the full force \eqref{eq:totalforcehom} on $[-0.5,0.5]^2$, extending it periodically on $\R^2$, and requiring that the force coefficients are differentiable and vanish on $\partial [-0.5,0.5]^2$ for physically realistic dynamics.
That is, we use \eqref{eq:totalforcehom} to define its periodic extension $\bar{F}\colon \R^2\to \R^2$ by
\begin{align}\label{eq:totalforcehomperiod}
\begin{split}
\bar{F}(d)&:=F(d) \quad \text{for }d\in [-0.5,0.5]^2,\\
\bar{F}(d+k)&:=\bar{F}(d) \quad \text{for }d\in [-0.5,0.5]^2,k\in\Z^2.
\end{split}
\end{align}	
Then, the particle model \eqref{eq:particlemodel} can be rewritten as
\begin{align}\label{eq:particlemodelperiodic}
\frac{\di x_j}{\di t}=\frac{1}{N}\sum_{\substack{k=1\\k\neq j}}^N \bar{F}(x_j-x_k)
\end{align}
for  $x_j\in\R^2$ where the right-hand side can be regarded as the force acting on  particle $j$.
For physically realistic forces, the force $\bar{F}$ has to vanish for any $d\in\partial [-0.5,0.5]^2$, implying that $f_l(0.5)=f_s(0.5)=0$ for $f_l,f_s$ in \eqref{eq:totalforcehom} and hence $f_l(|d|)=f_s(|d|)=0$ for  $d\in\R^2$ with $|d|=0.5$. Thus, we require that $\bar{F}(d)=0$ for all $d\in \partial [-0.5,0.5]^2$ for physically relevant forces.
To guarantee that the resulting force coefficient is differentiable which is required for the stability analysis we construct a differentiable approximation of the given force coefficient $f$ by considering $f(|d|)$ for $|d|\leq 0.5-\epsilon$ for some $\epsilon>0$, a cubic polynomial on $(0.5-\epsilon,0.5)$ and the constant zero function for $|d|\geq 0.5$  such that the resulting function is continuously differentiable on $(0,\infty)$. 
Motivated by this, we also consider smaller values of the cutoff radius $R_c\in(0,0.5]$ and adapt  the  force coefficients as
\begin{align}\label{eq:forcecutoff}
f^\epsilon(|d|)=\begin{cases} f(|d|), & |d|\in[0,R_c-\epsilon],\\
f'(R_c-\epsilon)\bl \frac{(|d|-R_c)^3}{\epsilon^2}+  \frac{\left(|d|-R_c \right)^2}{\epsilon}\br &\\\quad+ f(R_c-\epsilon)\bl 2\frac{(|d|-R_c)^3}{\epsilon^3}+  3\frac{\left(|d|-R_c \right)^2}{\epsilon^2} \br,  & |d| \in (R_c-\epsilon,R_c), 
\\ 0, & |d| \geq R_c.
\end{cases}
\end{align}
Note that this definition results in a  differentiable function whose absolute value and its derivative vanish for $|d|=R_c$.
This is in analogy to the notion of cutoff and is only a small modification compared to the original definition provided $f(s)$ for $s \in (R_c- \epsilon, R_c)$ is of order $\mathcal{O}(\epsilon)$ and $f'(R_c-\epsilon)$ is of order $\mathcal{O}(1)$. In this case, both the original force coefficients and its adaptation $f^\epsilon$ are of order $\mathcal{O}(\epsilon)$ on $(R_c-\epsilon,R_c)$. Further note that the interaction forces on distances $|d|\ll R_c-\epsilon$ are significantly larger than on $(R_c-\epsilon,R_c)$ and hence, the dynamics are mainly determined by interactions of range $|d|\ll R_c$. In particular, this allows us to replace $f_l$ and $f_s$ in \eqref{eq:totalforcehom} by differentiable approximations $f_l^\epsilon$ and $f^\epsilon_s$, defined as in \eqref{eq:forcecutoff}, if necessary.  

Note that the  assumption to consider the unit square $[0,1]^2$ with periodic boundary conditions is not restrictive and by rescaling in time our analysis extends to any domain $[0,\delta]^2$ with a cutoff radius $R_c\in(0,\frac{\delta}{2}]$ for $\delta\in \R_+$ where the cutoff of any force coefficient $f$ is defined in \eqref{eq:forcecutoff}.

The coefficient function $f_R$ of the repulsion force $F_R$ in \eqref{eq:repulsionforce} in  the K\"{u}cken-Champod model is originally of the form 
\begin{align}\label{eq:repulsionforcekc}
f_R(|d|)=(\alpha |d|^2+\beta)\exp(-e_R |d|)
\end{align}
for  $d\in\R^2$ and nonnegative parameters $\alpha$, $\beta$ and $e_R$. 
The coefficient function $f_A$  of the attraction force $F_A$ in \eqref{eq:attractionforce} 
is of the form 
\begin{align}\label{eq:attractionforcekc}
f_A(|d|)=-\gamma|d|\exp(-e_A|d|)
\end{align}
for $d\in\R^2$ and nonnegative constants $\gamma$ and $e_A$. To be as close as possible to the work by K\"ucken and Champod \cite{Merkel} we assume that the total force \eqref{eq:totalforce} exhibits short-range repulsion and long-range attraction along $l$ and one can choose the parameters as:
\begin{align}\label{eq:parametervaluesRepulsionAttraction}
\begin{split}
\alpha&=270, \quad \beta=0.1, \quad \gamma=35, \quad
e_A=95, \quad e_R=100, \quad \chi\in[0,1]
\end{split}
\end{align}  
as proposed in \cite{patternformationanisotropicmodel}.
Based on the adaptations of the force coefficients in \eqref{eq:forcecutoff}, we consider the modified K\"ucken-Champod force coefficients  in the sequel, given by
\begin{align}\label{eq:repulsionforcemodel}
f_R^\epsilon(|d|)=\begin{cases} f_R(|d|), & |d|\in[0,R_c-\epsilon],\\
f_R'(R_c-\epsilon)\bl \frac{(|d|-R_c)^3}{\epsilon^2}+  \frac{\left(|d|-R_c \right)^2}{\epsilon}\br &\\\quad+ f_R(R_c-\epsilon)\bl 2\frac{(|d|-R_c)^3}{\epsilon^3}+  3\frac{\left(|d|-R_c \right)^2}{\epsilon^2} \br,  & |d| \in (R_c-\epsilon,R_c), 
\\ 0, & |d| \geq R_c.
\end{cases}
\end{align}
and 
\begin{align}\label{eq:attractionforcemodel}
f_A^\epsilon(|d|)=\begin{cases} f_A(|d|), & |d|\in[0,R_c-\epsilon],\\
f_A'(R_c-\epsilon)\bl \frac{(|d|-R_c)^3}{\epsilon^2}+  \frac{\left(|d|-R_c \right)^2}{\epsilon}\br &\\\quad+ f_A(R_c-\epsilon)\bl 2\frac{(|d|-R_c)^3}{\epsilon^3}+  3\frac{\left(|d|-R_c \right)^2}{\epsilon^2} \br,  & |d| \in (R_c-\epsilon,R_c), 
\\ 0, & |d| \geq R_c.
\end{cases}
\end{align}
Here,  $f_R,f_A$ are very small in a neighborhood of the cutoff $R_c=0.5$ for the parameters in \eqref{eq:parametervaluesRepulsionAttraction}, or more generally, for $e_R$ and $e_A$ sufficiently large. Since the derivatives $f_R'$ and $f_A'$ also contain the exponential decaying terms $\exp(-e_R |d|)$ and $\exp(-e_A|d|)$, respectively, and are scaled by a factor $\mathcal{O}(\epsilon)$ in \eqref{eq:repulsionforcemodel} and \eqref{eq:attractionforcemodel}, respectively, the differences between $f_R^\epsilon$ and $f_R$, and $f_A^\epsilon$ and $f_A$, respectively, are very small compared to the size of the interaction forces at distances $|d|\ll R_c-\epsilon$ and the total force exerted on particle $x_j$, given by the right-hand side of \eqref{eq:particlemodelperiodic}. In particular, $f_R^\epsilon,f_A^\epsilon$ can be regarded as differentiable approximations of $f_R,f_A$.



\begin{figure}[htbp]
	\centering
	\includegraphics[width=0.45\textwidth]{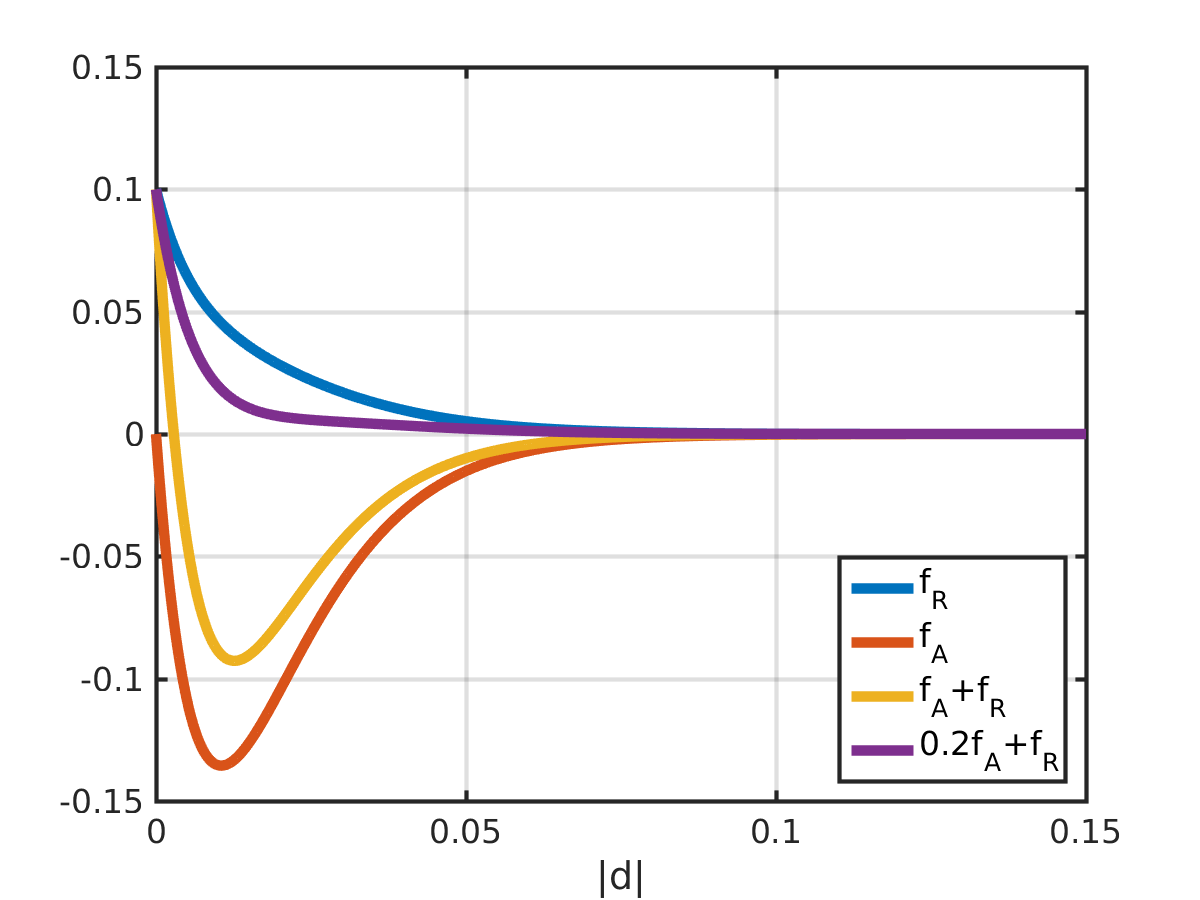}
	\includegraphics[width=0.45\textwidth]{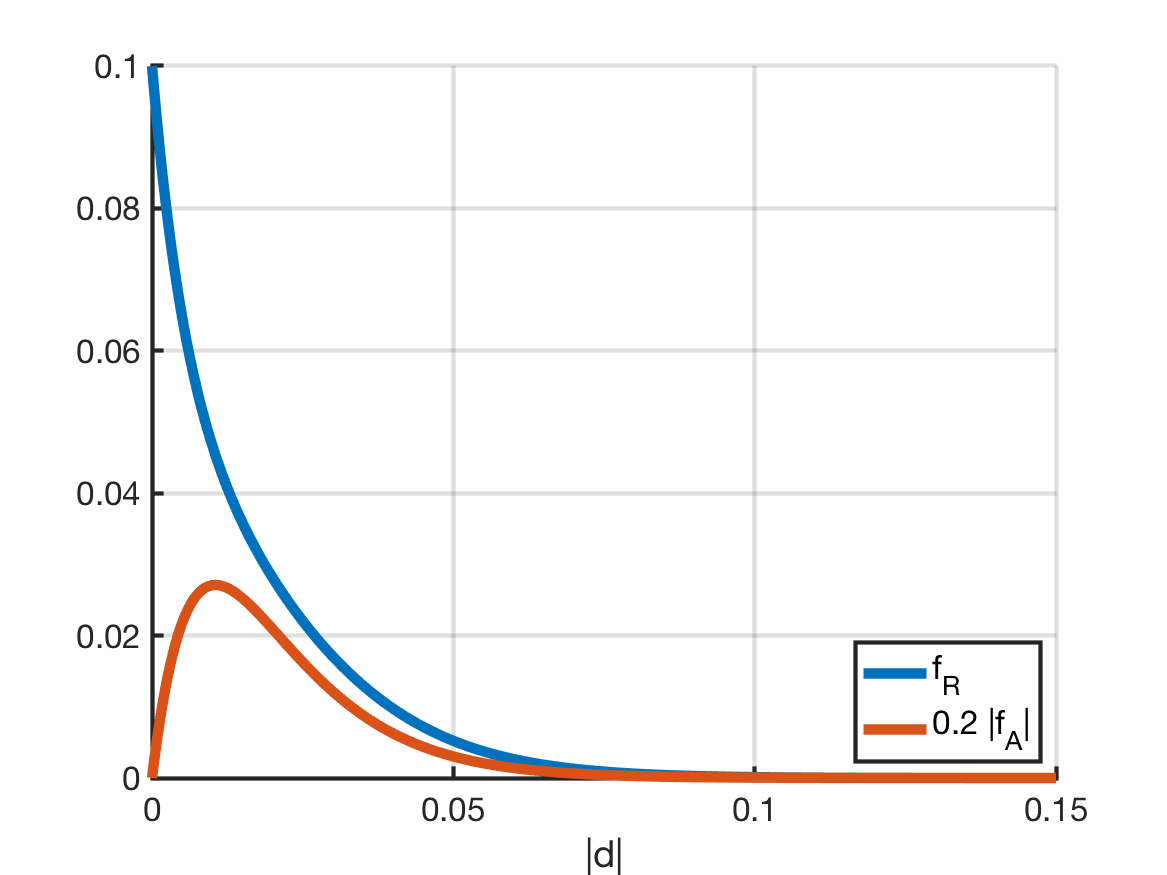}
	\caption{Coefficients $f_R$ in \eqref{eq:repulsionforcekc} and $f_A$ in \eqref{eq:attractionforcekc} of repulsion force \eqref{eq:repulsionforce} and attraction force \eqref{eq:attractionforce}, respectively, as well as  the force coefficients along $s=(0,1)$ and $l=(1,0)$  (i.e.\ $f_A+f_R$ and $0.2f_A+f_R$) for parameter values in  \eqref{eq:parametervaluesRepulsionAttraction}.\label{fig:forces}}
\end{figure}

For the particle model \eqref{eq:particlemodelperiodic} with differentiable coefficient functions $f_R^\epsilon,f_A^\epsilon$ and parameters  \eqref{eq:parametervaluesRepulsionAttraction}, we plot the original coefficient functions $f_R,f_A$ of the total force \eqref{eq:totalforcenewhom} for a spatially homogeneous underlying tensor field $T$ with $s=(0,1)$ and $l=(1,0)$ in Figure \ref{fig:forces}. However, note that $f_R\approx \lim_{\epsilon\to 0} f_R^\epsilon$ and $f_A\approx \lim_{\epsilon\to 0} f_A^\epsilon$. Moreover, we show the resulting coefficient functions $\chi f_A+f_R$ with $\chi=0.2$ and $f_A+f_R$ along $s=(0,1)$ and $l=(1,0)$, respectively,  in Figure~\ref{fig:forces}. Note that the repulsive force  coefficient $f_R$ is positive and the attractive force coefficient $f_A$ is negative. Repulsion dominates for short distances along $l$ to prevent  collisions of the particles. Besides, the total force exhibits   long-range attraction along $l$ whose absolute value decreases with the distance between particles. Along $s$,  the particles are purely repulsive for $\chi=0.2$ and the  repulsion force gets weaker for longer distances. 


\section{Stability/Instability of straight lines}\label{sec:stationarystates}
In this section, we consider the total force $\bar{F}$ in \eqref{eq:totalforcehomperiod}, defined on $\R^2$  by periodic extension of $F$  on $[-0.5,0.5]^2$ in \eqref{eq:totalforcehom}. This total force $\bar{F}$ can be described by (periodically extending) a short-range repulsive, long-range attractive force coefficient $f_l$ along $l$ and a purely repulsive force coefficient $f_s$ along $s$. Without loss of generality we may assume that the force coefficients $f_l,f_s$ are differentiable since otherwise they may be replaced by $f_l^\epsilon,f_s^\epsilon$, defined as in \eqref{eq:forcecutoff} for given functions $f_l,f_s$. Motivated by this we require in the sequel:

\begin{assumption}\label{ass:newforcegeneral}
	Let $f_l,f_s$ be continuously differentiable  functions  on $[0,\infty)$. Let  $f_s$ be purely repulsive, i.e.\ $f_s\geq 0$ with $f_s(0)>0$ for $s\in[0,R_c)$ and $f_s(s)=0$ for $s \geq R_c$, implying $\int_0^{R_c} f_s\di s>0$. Further let $f_l$ be short-range repulsive, long-range attractive with $f_l(R_c)= 0$.
\end{assumption}

As shown in \cite{patternformationanisotropicmodel} for the analysis of steady states with general spatially homogeneous tensor fields, it is sufficient to restrict  ourselves to the spatially homogeneous tensor field $T$ with  $s=(0,1)$ and $l=(1,0)$  in the sequel.

\subsection{Straight line}
In this section, we consider line patterns as steady states which were observed in the numerical simulations in \cite{patternformationanisotropicmodel}. 
For $x_j\in \R^2,j=1,\ldots,N$, evolving according to the particle model \eqref{eq:particlemodelperiodic}, we have 
\begin{align*}
\frac{\di}{\di t}\sum_{j=1}^N x_j=0,
\end{align*}
implying that  the centre of mass is conserved. Hence, we can assume without loss of generality that the centre of mass is in $\Z^2$. By identifying $\R^2$ with $\C$, we  make the ansatz
\begin{align}\label{eq:straightlineansatzgenerallength}
\bar{x}_k=\frac{k}{N}\exp(i\theta)\length(\theta),\quad k=1,\ldots,N.
\end{align}
Here, $\theta$ denotes the angle of rotation. The length of the line pattern is denoted by $\length=\length(\theta)>0$ and can be regarded as a multiplicative factor with $\length(0)=\length\bl\frac{\pi}{2}\br=1$ and $\length\bl \frac{\pi}{4}\br=\length\bl\frac{3\pi}{4}\br=\sqrt{2}$. Note that it is sufficient to restrict ourselves to $\theta\in[0,\pi)$ since ansatz \eqref{eq:straightlineansatzgenerallength} for $\theta$ and $\theta+k\pi$ with $k\in\Z$ leads to the same straight line after periodic extension on $\R^2$ and hence also on the torus $\mathbb{T}^2$. Depending on the choice of $\theta$, ansatz \eqref{eq:straightlineansatzgenerallength} might lead to multiple windings on the torus $\mathbb{T}^2$. To guarantee that ansatz \eqref{eq:straightlineansatzgenerallength} satisfies the periodic boundary conditions, we require that the winding number of the straight lines in \eqref{eq:straightlineansatzgenerallength} is a natural number and hence we can restrict ourselves to ansatz \eqref{eq:straightlineansatzgenerallength} on the torus $\mathbb{T}^2$ for $\theta\in\mathcal{A}$ where 
\begin{align}\label{eq:straightlineanglerotation}
\begin{split}
\mathcal{A}&:=\left\{0,\frac{\pi}{4},\frac{\pi}{2},\frac{3\pi}{4}\right\}\cup\left\{\psi\in\bl 0,\frac{\pi}{4}\br\cup\bl \frac{3\pi}{4},\pi\br\colon \cot(\psi)\in\Z\right\}\\&\qquad\cup\left\{\psi\in\bl \frac{\pi}{4},\frac{3\pi}{4}\br\colon \tan(\psi)\in\Z\right\}.
\end{split}
\end{align}
Note that considering the torus  $\mathbb{T}^2$ as the domain, i.e.\ the unit square with periodic boundary conditions, or equivalently $\R^2$ by periodic extension, is not restrictive due to the discussion in Section \ref{sec:modeldescription}.

For a single vertical straight line we have $\theta=\frac{\pi}{2}$ and ansatz \eqref{eq:straightlineansatzgenerallength} reduces to 
\begin{align}\label{eq:straightlineansatz}
\bar{x}_k=\frac{k}{N}i,\quad k=1,\ldots,N,
\end{align}
and for a horizontal line with $\theta=0$ we have
\begin{align}\label{eq:straightlinehorizontalansatz}
\bar{x}_k=\frac{k}{N},\quad k=1,\ldots,N.
\end{align}
Note that the winding number is one  for  \eqref{eq:straightlineansatzgenerallength} with $\theta\in \left\{0,\frac{\pi}{4},\frac{\pi}{2},\frac{3\pi}{4}\right\}$, while  
the winding number is larger than one for $\theta\in \mathcal{A}\backslash \left\{0,\frac{\pi}{4},\frac{\pi}{2},\frac{3\pi}{4}\right\}$. Translations of the ansatz \eqref{eq:straightlineansatzgenerallength}  result in steady states with a shifted centre of mass. Besides, parallel equidistant straight line patterns, obtained from considering \eqref{eq:straightlineansatzgenerallength} for a fixed rotation angle  \eqref{eq:straightlineanglerotation} and certain translations, may also
lead to steady states.

For equilibria $\bar{x}_j\in \R^2, j=1,\ldots,N,$ to the particle model \eqref{eq:particlemodelperiodic} we require that
\begin{align*}
\frac{1}{N}\sum_{\substack{k=1\\k\neq j}}^{N} \bar{F}(\bar{x}_j-\bar{x}_k,T)=0\quad \text{for all~}j=1,\ldots,N.
\end{align*}
Setting $\bar{x}_k$ for $k\in \Z$ as in \eqref{eq:straightlineansatzgenerallength}, we have 
$\bar{F}(\bar{x}_j-\bar{x}_k)=\bar{F}(\bar{x}_j-\bar{x}_{k+nN})$ for $j,k=1,\ldots,N$ and any $n\in\Z$ by the periodicity of $\bar{F}$.
Since the particles are uniformly distributed along straight lines by ansatz \eqref{eq:straightlineansatzgenerallength}, it is sufficient to require  
\begin{align}\label{eq:straightlinecondition}
\sum_{k=1}^{N-1} \bar{F}(\bar{x}_N-\bar{x}_k,T)=0
\end{align}
for  steady states.
Note that $\bar{F}(\bar{x}_N-\bar{x}_k,T)=-\bar{F}(\bar{x}_N-\bar{x}_{N-k},T)$ for $k=1,\ldots,\lceil N/2\rceil-1$ and for $N$ even we have $\bar{F}(\bar{x}_N-\bar{x}_{N/2},T)=0$ by the definition of the cutoff $R_c$.
Hence, \eqref{eq:straightlinecondition} is satisfied for the ansatz \eqref{eq:straightlineansatzgenerallength} for $\theta\in\mathcal{A}$, provided the length  $\length(\theta)$ of the lines is set such that the particles are distributed uniformly along the entire axis of angle $\theta$. 

\subsection{Stability conditions}
In this section we derive stability conditions for equilibria of the particle model \eqref{eq:particlemodelperiodic}, based on a linearised stability analysis. The  real parts of the eigenvalues of a stability matrix play a crucial role and we denote the real part of eigenvalue $\lambda \in \C$ by $\Re(\lambda)$ in the following.

\begin{proposition}
	For finite $N\in\N$, the steady state $\bar{x}_j,~j=1,\ldots,N$,  of the particle model \eqref{eq:particlemodelperiodic} is asymptotically  stable if the eigenvalues $\lambda$ of the stability matrix 
	\begin{align}\label{eq:stabilitymatrix}
	M=M(j,m)=\begin{pmatrix}
	I_1(j,m)  &I_2(j,m)
	\end{pmatrix}\in\C^{2,2}
	\end{align}
	satisfy $\Re(\lambda)< 0$ for all $j=1,\ldots,N$ and $m=1,\ldots,N-1$
	where
	\begin{align}\label{eq:stabilitymatrixentries}
	\begin{split}
	I_1(j,m)&=\frac{1}{N}\sum_{k\neq j}\bl 1-\exp(im(\phi_k-\phi_j))\br\frac{\partial \bar{F}}{\partial d_1}(\bar{x}_j-\bar{x}_k)\\&=\frac{1}{N}\sum_{k\neq j}\bl 1-\exp\bl \frac{2\pi im(k-j)}{N}\br\br\frac{\partial \bar{F}}{\partial d_1}(\bar{x}_j-\bar{x}_k),\\
	I_2(j,m)&=\frac{1}{N}\sum_{k\neq j}\bl 1-\exp(im(\phi_k-\phi_j))\br\frac{\partial \bar{F}}{\partial d_2}(\bar{x}_j-\bar{x}_k)\\&=\frac{1}{N}\sum_{k\neq j}\bl 1-\exp\bl \frac{2\pi im(k-j)}{N}\br\br\frac{\partial \bar{F}}{\partial d_2}(\bar{x}_j-\bar{x}_k)
	\end{split}
	\end{align}
	for $j=1,\ldots,N$ and $m=1,\ldots,N$. 
\end{proposition}

\begin{proof}
	Let $\bar{x}_j,~j=1,\ldots,N,$ denote a steady state  of \eqref{eq:particlemodelperiodic}. We define the perturbation ${g_j=g_j(t)}$, $h_j=h_j(t)\in\R$ of $\bar{x}_j$ by
	\begin{align*}
	x_j=\bar{x}_j+\begin{pmatrix}
	g_j\\h_j
	\end{pmatrix},\quad j=1,\ldots,N.
	\end{align*}
	Linearising \eqref{eq:particlemodelperiodic} around the steady state $\bar{x}_j$ gives
	\begin{align}\label{eq:lineraisedstabilityeq}
	\frac{\di}{\di t}\begin{pmatrix}
	g_j\\h_j
	\end{pmatrix}=\frac{1}{N}\sum_{k\neq j}(g_j-g_k)\frac{\partial \bar{F}}{\partial d_1}(\bar{x}_j-\bar{x}_k)+\frac{1}{N}\sum_{k\neq j}(h_j-h_k)\frac{\partial \bar{F}}{\partial d_2}(\bar{x}_j-\bar{x}_k).
	\end{align}
	We choose the ansatz functions
	\begin{align*}
	g_j=\zeta_g\bl \exp(im\phi_j)+\exp(-im\phi_j)\br,\quad h_j=\zeta_h\bl \exp(im\phi_j)+\exp(-im\phi_j)\br,\\ j=1,\ldots,N,\quad m=1,\ldots,N,
	\end{align*}
	where $\zeta_g=\zeta_g(t),\zeta_h=\zeta_h(t)$ and $\phi_j=\frac{2\pi j}{N}$.  Note that $g_j,h_j\in\R$ for all $j=1,\ldots,N$ and
	\begin{align*}
	\sum_{j=1}^N \exp(im\phi_j)=\sum_{j=1}^N \bl\exp\bl \frac{2\pi im }{N}\br\br^j=\begin{cases} 0, & m=1,\ldots,N-1,\\N, &m=N,\end{cases}
	\end{align*}
	since $\phi_j$ are the roots of $r^N=1$ and 
	\begin{align*}
	\sum_{j=0}^{N-1} r^j=\frac{1-r^N}{1-r}.
	\end{align*}
	This implies 
	\begin{align*}
	\sum_{j=1}^N g_j(t)=\sum_{j=1}^N h_j(t)=\begin{cases}0, & m=1,\ldots,N-1,\\N, &m=N,\end{cases}
	\end{align*}
	for all times $t\geq 0$, i.e. the centre of mass of the perturbations $g_j,h_j$ is preserved.
	We have
	\begin{align*}
	g_j-g_k&=\zeta_g\bl \exp(im\phi_j)+\exp(-im\phi_j)\br\bl 1-\exp(im(\phi_k-\phi_j))\br,\\ h_j-h_k&=\zeta_h\bl \exp(im\phi_j)+\exp(-im\phi_j)\br\bl 1-\exp(im(\phi_k-\phi_j))\br.
	\end{align*}
	Plugging this into \eqref{eq:lineraisedstabilityeq} and collecting like terms in $\exp(im\phi_j),\exp(-im\phi_j)$ results in
	\begin{align*}
	\frac{\di}{\di t}\begin{pmatrix}
	\zeta_g\\\zeta_h
	\end{pmatrix}&=\frac{\zeta_g}{N}\sum_{k\neq j}\bl 1-\exp(im(\phi_k-\phi_j))\br\frac{\partial \bar{F}}{\partial d_1}(\bar{x}_j-\bar{x}_k)\\\nonumber&\quad +\frac{\zeta_h}{N}\sum_{k\neq j}\bl 1-\exp(im(\phi_k-\phi_j))\br\frac{\partial \bar{F}}{\partial d_2}(\bar{x}_j-\bar{x}_k),
	\end{align*}
	i.e.\
	\begin{align}\label{eq:stabilitysystem}
	\frac{\di}{\di t}\begin{pmatrix}
	\zeta_g\\\zeta_h
	\end{pmatrix}=M\begin{pmatrix}
	\zeta_g\\\zeta_h
	\end{pmatrix},
	\end{align}
	where the stability matrix $M\in\C^{2,2}$ is defined in \eqref{eq:stabilitymatrix}. The ansatz $\zeta_g=\xi_g\exp(\lambda t)$, ${\zeta_h=\xi_h\exp(\lambda t)}$ solves the system \eqref{eq:stabilitysystem} for any eigenvalue $\lambda\in\C$ of the stability matrix $M=M(j,m)$. Note that the stability matrix $M$ is the zero matrix for $m=N$ and any $j=1,\ldots,N$. Hence,  we have $\lambda=0$ for $m=N$ and all $j=1,\ldots,N$, corresponding to translations along the vertical and horizontal axis. Thus, the straight  line $\bar{x}_j,j=1,\ldots,N,$ is stable if $\Re(\lambda)<0$ for any $j=1,\ldots,N$ and $m=1,\ldots,N-1$.
\end{proof}

\subsection{Stability  of a single  vertical straight line}\label{sec:verticalstraightline}

To study the stability of a single vertical straight line of the form \eqref{eq:straightlineansatz} we determine the eigenvalues of the stability matrix \eqref{eq:stabilitymatrix} and derive  stability conditions for steady states $\bar{x}_j,j=1,\ldots,N,$ satisfying \eqref{eq:straightlinecondition}.
In the continuum limit $N\to\infty$ the steady state condition \eqref{eq:straightlinecondition} becomes
$$\int_{-0.5}^{0.5} F((0,s),T)\di s=\int_{-0.5}^{0.5} \bar{F}((0,s),T)\di s=0.$$
Due to the cutoff radius $R_c\in(0,0.5]$ it is sufficient to require
\begin{align}\label{eq:straightlineconditioncont}
\int_{-R_c}^{R_c} F((0,s),T)\di s=0
\end{align}
for equilibria. This condition is clearly satisfied for forces of the form \eqref{eq:totalforcehom} and in particular for forces of the form \eqref{eq:totalforcenewhom}.

\begin{theorem}\label{th:stabilitystraightline}
	For finite $N\in\N$, the single vertical straight line \eqref{eq:straightlineansatz} is an asymptotically stable steady state of the particle model \eqref{eq:particlemodelperiodic} with total force \eqref{eq:totalforcehom}  if $\Re(\lambda_{i,N}(m))< 0$ for $i=1,2$ and all $m=1,\ldots,N-1$ where the eigenvalues $\lambda_{i,N}=\lambda_{i,N}(m)$ of the stability matrix \eqref{eq:stabilitymatrix} are given by
	\begin{align}\label{eq:eigenvalueslinegeneraldiscrete}
	\begin{split}
	\lambda_{1,N}(m)&=\frac{1}{N}\sum_{k=\lceil \frac{N}{2}\rceil}^{N-1+\lceil \frac{N}{2}\rceil} f_l(|d_{Nk}|)\bl 1-\exp\bl \frac{2\pi imk}{N}\br\br,\\
	\lambda_{2,N}(m)&=\frac{1}{N}\sum_{k=\lceil \frac{N}{2}\rceil}^{N-1+\lceil \frac{N}{2}\rceil} \bl f_s(|d_{Nk}|)+ f_s'(|d_{Nk}|)|d_{Nk}|\br\bl 1-\exp\bl \frac{2\pi imk}{N}\br\br
	\end{split}
	\end{align}
	with $$d_{Nk}=\begin{pmatrix}
	0\\\frac{N-k}{N}
	\end{pmatrix}$$ 
	for $k\in\N$.
	Denoting the cutoff radius by $R_c\in(0,0.5]$, steady states satisfying the steady state condition \eqref{eq:straightlineconditioncont} in the continuum limit $N\to\infty$ are  unstable if $\Re(\lambda_i(m))> 0$ for some $m\in\N$ and some $i\in\{1,2\}$ where the eigenvalues $\lambda_i=\lambda_i(m),i=1,2$, of the stability matrix \eqref{eq:stabilitymatrix}  are given by
	\begin{align}\label{eq:eigenvalueslinegeneral}
	\begin{split}
	\lambda_1(m)&=\int_{-R_c}^{R_c}  f_l(|s|)\bl 1-\exp\bl - 2\pi ims\br\br\di s,\\
	\lambda_2(m)&=\int_{-R_c}^{R_c}
	\bl  f_s(|s|)+ f_s'(|s|)|s|
	\br\bl 1-\exp\bl -2\pi ims\br\br\di s.
	\end{split}
	\end{align}
	In particular, 
	\begin{align}\label{eq:eigenvaluesreallinegeneral}
	\begin{split}
	\Re(\lambda_1)(m)&=2\int_{0}^{R_c}  f_l(s)\bl 1-\cos\bl - 2\pi ms\br\br\di s,\\
	\Re(\lambda_2)(m)&=2\int_{0}^{R_c}
	\bl  f_s(s)+ f_s'(s)s
	\br\bl 1-\cos\bl -2\pi ms\br\br\di s.
	\end{split}
	\end{align}
\end{theorem}

\begin{proof}
	For the spatially homogeneous tensor field $T$, defined by $s=(0,1)$ and $l=(1,0)$, the derivatives of the total force \eqref{eq:totalforcehom} are given by 
	\begin{align}\label{eq:straightlinederivativeforce}
	\frac{\partial \bar{F}}{\partial d_1} (d)=\begin{pmatrix}
	f_l(|d|)+f_l'(|d|)\frac{d_1^2}{|d|}\\ f_s'(|d|)\frac{d_1d_2}{|d|}
	\end{pmatrix}, \qquad
	\frac{\partial \bar{F}}{\partial d_2} (d)=\begin{pmatrix}
	f_l'(|d|)\frac{d_1d_2}{|d|}\\ f_s(|d|)+ f_s'(|d|)\frac{d_2^2}{|d|}
	\end{pmatrix}
	\end{align}
	for $d=(d_1,d_2)\in [-0.5,0.5]^2$ and its periodic extension $\frac{\partial \bar{F}}{\partial d_i}(d+k)=\frac{\partial \bar{F}}{\partial d_i}(d)$ for $i=1,2$, ${d\in [-0.5,0.5]^2}$ and $k\in\Z^2$. Note that $f_l, f_s$ are differentiable due to the smoothing assumptions at the cutoff $R_c$ in \eqref{eq:forcecutoff} and their derivatives vanish for $d\in [-0.5,0.5]^2$ with $|d|\geq R_c$. 
	Using  ansatz \eqref{eq:straightlineansatz} for a single  vertical straight line, we obtain: 
	\begin{align}\label{eq:derivativeforcesrepattr}
	\begin{split}
	\frac{\partial \bar{F}}{\partial d_1}(d_{jk})&=\begin{pmatrix}f_l(|d_{jk}|)\\0\end{pmatrix},\\
	\frac{\partial \bar{F}}{\partial d_2}(d_{jk})&=\begin{pmatrix}0 \\f_s(|d_{jk}|)+f_s'(|d_{jk}|)|d_{jk}|\end{pmatrix}
	\end{split}
	\end{align}
	for $d_{jk}\in [-0.5,0.5]^2$ and note that $\frac{\partial \bar{F}}{\partial d_i}(d_{jk})=\frac{\partial \bar{F}}{\partial d_i}(d_{j,k+nN})$ for $i=1,2$, $j,k=1,\ldots,N$, and $n\in\N$.
	This implies that the particles along the straight vertical line are indistinguishable and it suffices to consider $j=N$.  The entries \eqref{eq:stabilitymatrixentries} of the stability matrix \eqref{eq:stabilitymatrix} are given by 
	\begin{align*}
	\begin{split}
	I_1(m)&=\frac{1}{N}\sum_{k=1}^{N}\bl 1-\exp\bl \frac{2\pi imk}{N}\br\br\frac{\partial \bar{F}}{\partial d_1}\bl d_{Nk}\br,\\
	I_2(m)&=\frac{1}{N}\sum_{k=1}^{N}\bl 1-\exp\bl \frac{2\pi imk}{N}\br\br\frac{\partial \bar{F}}{\partial d_2}\bl d_{Nk}\br. 
	\end{split}
	\end{align*}
	Note that for $k=\lceil \frac{N}{2}\rceil,\ldots,N$, we have $d_{Nk}\in\{0\}\times [0,0.5]\subset [-0.5,0.5]^2$, implying that the derivatives of $\bar{F}$ are given by \eqref{eq:derivativeforcesrepattr} where $\bar{F}(d_{N,\lceil \frac{N}{2}\rceil})=0$ by the definition of the cutoff $R_c$ for $N$ even. Since $\frac{\partial \bar{F}}{\partial d_{i}}(d_{Nk})=\frac{\partial \bar{F}}{\partial d_{i}}(d_{N,N+k})$ for $i=1,2$, $k=1,\ldots,\lceil N/2\rceil-1$, and $d_{N,N+k}\in\{0\}\times (-0.5,0)\subset [-0.5,0.5]^2$, we can replace the sum over $k\in\{1,\ldots,N\}$ by the sum over $k\in\{\lceil \frac{N}{2}\rceil,\ldots, N-1+\lceil \frac{N}{2}\rceil\}$, resulting in
	\begin{align}\label{eq:stabilitymatrixentriesline}
	\begin{split}
	I_1(m)&=\frac{1}{N}\sum_{k=\lceil \frac{N}{2}\rceil}^{N-1+\lceil \frac{N}{2}\rceil}\bl 1-\exp\bl \frac{2\pi imk}{N}\br\br\frac{\partial \bar{F}}{\partial d_1}\bl d_{Nk}\br,\\
	I_2(m)&=\frac{1}{N}\sum_{k=\lceil \frac{N}{2}\rceil}^{N-1+\lceil \frac{N}{2}\rceil}\bl 1-\exp\bl \frac{2\pi imk}{N}\br\br\frac{\partial \bar{F}}{\partial d_2}\bl d_{Nk}\br. 
	\end{split}
	\end{align}
	
	Note that the stability matrix \eqref{eq:stabilitymatrix} is a diagonal matrix whose eigenvalues are   the non-trivial entries in \eqref{eq:stabilitymatrixentriesline} and are given by \eqref{eq:eigenvalueslinegeneraldiscrete}. Since the sums in \eqref{eq:stabilitymatrixentriesline} are Riemannian sums   we can pass to the continuum limit $N\to \infty$.  
	Note that $\frac{k}{n}\in [0.5,1.5]$
	for $k\in \{\lceil \frac{N}{2}\rceil,\ldots,N-1+\lceil \frac{N}{2}\rceil\}$ appears in the entries of the stability matrix \eqref{eq:stabilitymatrixentriesline}. 
	For passing to the limit $N\to\infty$ in \eqref{eq:stabilitymatrixentriesline}, we consider the domain of integration $[0.5,1.5]$ and do a change of variables resulting in
	\begin{align*}
	I_i(m)&=\int_{0.5}^{1.5}\frac{\partial \bar{F} }{\partial d_i}\bl (0,1-s)\br\bl 1-\exp\bl 2\pi ims\br\br\di s=\int_{-0.5}^{0.5}\frac{\partial \bar{F}}{\partial d_i} \bl (0,s)\br\bl 1-\exp\bl - 2\pi ims\br\br\di s \\&=\int_{-0.5}^{0.5}\frac{\partial F}{\partial d_i} \bl (0,s)\br\bl 1-\exp\bl - 2\pi ims\br\br\di s
	\end{align*}
	for $i=1,2$ and all $m\in\N$.
	Clearly the stability matrix \eqref{eq:stabilitymatrix} with entries $I_i,i=1,2,$ is again a diagonal matrix and the eigenvalues $\lambda_i=\lambda_i(m), i=1,2$, in \eqref{eq:eigenvalueslinegeneral} are given by the diagonal entries of the stability matrix \eqref{eq:stabilitymatrix}.
\end{proof}

\begin{remark}
	In Theorem \ref{th:stabilitystraightline} we study the stability of the straight vertical line for the dynamical system \eqref{eq:particlemodelperiodic} for a finite number of particles $N$ where the differentiability of $\bar{F}$ at the cutoff $R_c$ is necessary for the definition of the eigenvalues in the discrete setting in \eqref{eq:eigenvalueslinegeneraldiscrete}. Note that we cannot conclude stability/instability if $\Re(\lambda_i(m))\leq 0$ for $i=1,2$ and all $m=1,\ldots,N-1$. By the assumptions on the force coefficients $f_s,f_l$ in Assumption \ref{ass:newforcegeneral} we can pass to the continuum limit $N\to\infty$ in the definition of the eigenvalues of the stability matrix and study the stability of the steady states of the particle model \eqref{eq:particlemodelperiodic} in the continuum limit $N\to\infty$. If there exists $m\in \N$ for some $i\in\{1,2\}$ such that $\Re(\lambda_i(m ))>0$, then the steady state is unstable in the continuum limit.  However, if $\Re(\lambda_i(m))\leq 0$ for $i\in\{1,2\}$ and all $m\in\N$ stability/instability of the steady state cannot be concluded since it is difficult to give general conditions for $\Re(\lambda_i(m))\to \sigma$ as $m\to \infty$ with $\sigma=0$ or $\sigma\in\R_-\backslash \{0\}$. If $\sigma=0$, we cannot say anything about the stability/instability of the steady state in the continuum setting, see also similar discussions for the stability/instability of Delta-rings in the continuum setting in \cite{simionethesis} and the discussion after Theorem 2.1 in \cite{Bertozzi2015}. In particular,  linear stability for any  $N\in\N$ is not sufficient to conclude stability in the continuum setting.
\end{remark} 

Note that the asymmetry in the definition of the eigenvalues \eqref{eq:eigenvaluesreallinegeneral} is due to the asymmetric  steady states in \eqref{eq:straightlineansatz}. For $f=f_s=f_l$ the total force in \eqref{eq:totalforcehom} simplifies to $F(d)=f(|d|)d$ for $d=(d_1,d_2)\in [-0.5,0.5]^2$. In this case, the gradient of $F=(F_1,F_2)$ is a symmetric matrix, compare \eqref{eq:straightlinederivativeforce}, and hence the eigenvalues of the stability matrix are real. Since 
\begin{align*}
\frac{\partial F_1}{\partial d_2}=\frac{\partial F_2}{\partial d_1}
\end{align*}
there exists a radially symmetric potential $W(d)=w(|d|)$ such that $F=-\nabla W$  on $[-0.5,0.5]^2$. 
Hence, the stability conditions can be derived in terms of the potential $w$ and we have $$\text{trace}(\nabla F(d))=f'(|d|)|d|+2f(|d|)=-\Delta w(|d|)=\lambda_1+\lambda_2$$ for $d\in[-0.5,0.5]^2$ and the periodic extension $\bar{F}$ of $F$ can be considered on $\R^2$. For $f_s=f_l$ and radially symmetric steady states, this leads to identical conditions for both eigenvalues $\lambda_k$, $k=1,2$. For the analysis of these symmetric steady states, however, it is helpful to consider an appropriate coordinate system such as polar coordinates for ring steady states as in \cite{Bertozzi2015}.

Note that the stability conditions for steady states depend on the choice of the coordinate system. Considering derivatives with respect to the coordinate axes as in \eqref{eq:straightlinederivativeforce} seems to be the natural choice for straight line patterns, in contrast to polar coordinates as in \cite{Bertozzi2015}.

In the sequel, we investigate the high-wave number stability of straight line patterns for the particle model \eqref{eq:particlemodelperiodic}, i.e.\ the stability of straight vertical lines as $m\to \infty$. This can be studied by considering the limit $m\to\infty$ of the eigenvalues \eqref{eq:eigenvalueslinegeneral} of the stability matrix \eqref{eq:stabilitymatrix} associated with the dynamical system \eqref{eq:stabilitysystem}.
\begin{proposition}\label{prop:highwavestraightline}
	Suppose that the coefficient functions $f_s$ and $f_l$ are continuously differentiable on $[0,+\infty)$ with  $f_s(|d|)=f_l(|d|)=0$ for $|d|\geq R_c$
	and  $f_s\geq 0$. The condition  
	\begin{align}\label{eq:condvertlinehighwave}
	\begin{split}
	\int_0^{R_c} f_l(s)\di s \leq 0 \quad\text{and}\quad f_s(R_c)= 0
	\end{split}
	\end{align}
	is necessary for the  high-wave number stability of the single straight vertical line \eqref{eq:straightlineansatz}, i.e.\ \eqref{eq:condvertlinehighwave} is necessary for the stability of the straight vertical line for any $N\in\N$ and in the continuum limit $N\to\infty$.

\end{proposition}
\begin{proof}
	The eigenvalues \eqref{eq:eigenvalueslinegeneral} of the stability matrix \eqref{eq:stabilitymatrix} associated with the equilibrium of a single vertical straight line are of the form
	\begin{align*}
	\lambda(m)&=\int_{-R_c}^{R_c} f(|s|)(1-\exp(-2\pi i ms))\di s\\
	&=2\int_0^{R_c} f(s)\di s-\frac{1}{2\pi i m}\int_{-R_c}^{R_c} f'(|s|)\exp(-2\pi i ms)\di s\\&\quad+\frac{1}{2\pi im}f(R_c)\bl \exp(-2\pi imR_c)-\exp(2\pi imR_c)\br
	\end{align*}
	for a function $f\colon\R_+\to\R$ with $f(|d|)=0$ for $|d|\geq R_c$.
	For high-wave number stability we require 
	\begin{align*}
	\int_0^{R_c} f(s)\di s \leq 0 \qquad\text{and}\qquad |f'| \text{ is integrable on }[0,R_c].
	\end{align*}
	Then, using the definition of the eigenvalues \eqref{eq:eigenvalueslinegeneral} this leads to the conditions 
	\begin{align}\label{eq:condvertlinehighwavehelp}
	\begin{split}
	\int_0^{R_c} f_l(s)\di s\leq 0\qquad \text{and}\qquad \int_0^{R_c}   f_s(s)+ f_s'(s)s\di s\leq 0.
	\end{split}
	\end{align}
	Integration by parts of the second condition in \eqref{eq:condvertlinehighwavehelp} leads to $f_s(R_c)\leq 0$ and the conditions in \eqref{eq:condvertlinehighwave} result from $f_s$ being repulsive, i.e.\ $f_s\geq 0$.
\end{proof}

\begin{remark}
	The necessary condition $f_s(R_c)=0$ in \eqref{eq:condvertlinehighwave} for a stable straight vertical line  is equivalent to the eigenvalue associated with $f_s$ to be equal to zero in the high-wave number limit. Hence, stability/instability of straight vertical lines cannot be concluded in the continuum limit $N\to\infty$ from the linear stability analysis.
\end{remark}

The first condition in \eqref{eq:condvertlinehighwave} implies that the total attractive force over its entire range is larger than the total repulsive force along $l$. 
The second condition in \eqref{eq:condvertlinehighwave} implies that for high-wave stability we require the total force at the cutoff radius $R_c$ should not be repulsive along $s$ which is identical to the assumptions on the cutoff in \eqref{eq:forcecutoff}.

In comparison with the high-wave number conditions in \eqref{eq:condvertlinehighwavehelp} in the proof of Proposition~\ref{prop:highwavestraightline} the integrands for the stability conditions are multiplied by a factor 
\begin{align*}
\Re (1-~\exp(-2\pi ims))=1-\cos(2\pi ms)\in[0,2].
\end{align*}
Even if the necessary  conditions for high-wave number stability \eqref{eq:condvertlinehighwave} are satisfied, this does not guarantee that  $\Re(\lambda_1(m)),\Re(\lambda_2(m))\leq 0$ for all $m\in\N$ and hence necessary  stability conditions for  the single vertical straight line might not be satisfied for all $m\in\N$. 

The general stability  conditions for straight vertical lines can be obtained from the real parts of the eigenvalues \eqref{eq:eigenvalueslinegeneral} of the stability matrix \eqref{eq:stabilitymatrix}. The conditions \eqref{eq:condvertlinehighwave} suggest that stability of the straight line is possible for particular force coefficient choices. This will be investigated in Section \ref{sec:verticallines}.

\begin{remark}
	Note that differentiable approximations $f_R^\epsilon,f_A^\epsilon$ of the force coefficients $f_R$ and $f_A$ in the K\"ucken-Champod model are defined in \eqref{eq:repulsionforcemodel} and \eqref{eq:attractionforcemodel}, respectively. Setting $f_l^\epsilon:=f_A^\epsilon+f_R^\epsilon$ and $f_s^\epsilon:=\chi f_A^\epsilon+f_R^\epsilon$ for some $0<\epsilon\ll R_c$ and a parameter $\chi\in[0,1]$ such that $f_s \geq 0$ on $[0,R_c)$ we consider $f_l^\epsilon,f_s^\epsilon$ instead of $f_l,f_s$ in the definition of the real parts of the eigenvalues  \eqref{eq:eigenvaluesreallinegeneral}. We	
	obtain for the real parts of the eigenvalues of the stability matrix  \eqref{eq:stabilitymatrix} in the K\"ucken-Champod model with total force \eqref{eq:totalforcenewhom} and the spatially homogeneous tensor field $T$ in \eqref{eq:tensorfieldhom}:
	\begin{align*}
	\begin{split}
	\Re(\lambda_1(m))&=2\int_{0}^{R_c} \bl f_A^\epsilon(s)+f_R^\epsilon(s)\br \bl 1-\cos\bl - 2\pi ms\br\br\di s,\\
	\Re(\lambda_2(m))&=2\int_{0}^{R_c}
	\bl  \chi f_A^\epsilon(s)+f_R^\epsilon(s)+ \chi s(f^\epsilon_A)'(s)+s(f_R^\epsilon)'(s)
	\br\bl 1-\cos\bl -2\pi ms\br\br\di s.
	\end{split}
	\end{align*} 
	The necessary stability condition \eqref{eq:condvertlinehighwave} implies that ${f_s^\epsilon(R_c)=\chi f_A^\epsilon(R_c)+f_R^\epsilon(R_c)=0}$, consistent with the definition of the force coefficients \eqref{eq:repulsionforcemodel} and \eqref{eq:attractionforcemodel} in the K\"ucken-Champod model. 
	Hence, the necessary condition \eqref{eq:condvertlinehighwave} for high-wave number stability of a straight vertical line is satisfied in this case.
\end{remark}

\subsection{Instability of a single  horizontal straight line}
In this section we investigate the stability of a single horizontal straight line given by the ansatz \eqref{eq:straightlinehorizontalansatz}  which follows from \eqref{eq:straightlineansatzgenerallength} with $\theta=0$.  
\begin{theorem}
	For $N\in\N$ sufficiently large and in the continuum limit $N\to\infty$, the single horizontal straight line \eqref{eq:straightlinehorizontalansatz}  is an unstable steady state to the particle model \eqref{eq:particlemodelperiodic} for any choice of force coefficients $f_s$ and $f_l$ of the total force \eqref{eq:totalforcehom}, provided the total force is purely repulsive along $s$ on $[0,R_c)$.
\end{theorem}
\begin{proof}
	For  a single horizontal straight line, we have
	$$d_{jk}=\bar{x}_j-\bar{x}_k=\begin{pmatrix}
	\frac{j-k}{N}\\ 0
	\end{pmatrix}$$ 
	and the  derivatives of the total force are given by
	\begin{align*}
	\frac{\partial}{\partial d_1} \bar{F}(d_{jk})&=\begin{pmatrix}
	f_l(|d_{jk}|)+f_l'(|d_{jk}|)|d_{jk}|\\0
	\end{pmatrix},\\
	\frac{\partial}{\partial d_2} \bar{F}(d_{jk})&=\begin{pmatrix}
	0\\f_s(|d_{jk}|)
	\end{pmatrix}
	\end{align*}
	for $d_{jk}\in [-0.5,0.5]^2$.
	Similar as in Section \ref{sec:verticalstraightline} one can show that the  eigenvalues $\lambda_1=\lambda_1(m)$, $\lambda_2=\lambda_2(m)$ of the stability matrix \eqref{eq:stabilitymatrix} are given by
	\begin{align*}
	\begin{split}
	\lambda_1(m)&=2\int_0^{R_c} \bl  f_l(s)+ f_l'(s)s\br\bl 1-\exp\bl - 2\pi ims\br\br\di s,\\
	\lambda_2(m)&=2\int_0^{R_c}
	f_s(s)
	\bl 1-\exp\bl -2\pi ims\br\br\di s
	\end{split}
	\end{align*}
	for a   cutoff radius $R_c\in(0,0.5]$. For high-wave number stability we require 
	\begin{align*}
	\begin{split}
	f_l(R_c)=\frac{1}{R_c}\int_0^{R_c}   f_l(s)+ f_l'(s)s\di s\leq 0 \qquad \text{and}\qquad
	\int_0^{R_c}    f_s(s)\di s\leq 0.
	\end{split}
	\end{align*}
	The forces are assumed to be purely repulsive along $s$ up to the cutoff $R_c$, i.e.\ $f_s > 0$ on $[0,R_c)$, implying $$\int_0^{R_c}   f_s(s)\di s>0.$$ Hence, the single horizontal straight line is high-wave unstable.
\end{proof}

\subsection{Instability of rotated straight line patterns}

In this section we consider the ansatz \eqref{eq:straightlineansatzgenerallength} where the angle of rotation $\theta$ satisfies \eqref{eq:straightlineanglerotation}, resulting in rotated straight line patterns. The entries of  the  stability matrix \eqref{eq:stabilitymatrix} are given by
\begin{align*}
\begin{split}
I_1(m)&=2\int_0^{R_c}\frac{\partial \bar{F} }{\partial d_1}\bl \bl s\cos\bl\theta\br,s \sin\bl \theta\br\br\br\bl 1-\exp\bl - 2\pi ims\br\br\di s,\\
I_2(m)&=2\int_0^{R_c}\frac{\partial \bar{F}}{\partial d_2}\bl \bl s\cos\bl\theta \br,s \sin\bl\theta\br\br\br\bl 1-\exp\bl -2\pi ims\br\br\di s,
\end{split}
\end{align*}
where the derivatives of the total force can easily be determined by

\begin{align*}
\frac{\partial}{\partial d_1} \bar{F}(d)=\begin{pmatrix}
f_l(|d|)+f_l'(|d|)\frac{d_1^2}{|d|}\\ f_s'(|d|)\frac{d_1d_2}{|d|}
\end{pmatrix},\qquad
\frac{\partial}{\partial d_2} \bar{F}(d)=\begin{pmatrix}
f_l'(|d|)\frac{d_1d_2}{|d|}\\ f_s(|d|)+ f_s'(|d|)\frac{d_2^2}{|d|}
\end{pmatrix}
\end{align*}
$d\in[-0.5,0.5]^2$ with the cutoff radius $R_c\in(0,0.5]$. In particular,  the  stability matrix \eqref{eq:stabilitymatrix} is no longer a diagonal matrix in general. To show that the rotated straight line pattern is unstable for $\theta\in(0,\pi)\setminus  [\phi,\pi-\phi]$ for  some $\phi\in(0,\frac{\pi}{2})$ and $N\in \N$ sufficiently large and in the continuum limit $N\to \infty$, it is sufficient to consider the high frequency wave limit and show high-wave number instability. Denoting the entries of $I_k$ by $I_{k1}$ and $I_{k2}$ for $k=1,2$ with $M=(I_1,I_2)$  the high-frequency limit leads to
\begin{align}\label{eq:rotatestabilityentries}
\begin{split}
I_{11}=2\int_0^{R_c} f_l(s)\di s+2\int_0^{R_c} f_l'(s)s\cos^2\bl\theta\br \di s,\qquad
I_{12}=2\int_0^{R_c} f_s'(s)s\sin\bl\theta\br\cos\bl \theta\br \di s,\\
I_{21}=2\int_0^{R_c} f_l'(s)s\sin\bl\theta\br\cos\bl \theta\br \di s,\qquad
I_{22}=2\int_0^{R_c} f_s(s)\di s+2\int_0^{R_c} f_s'(s)s\sin^2\bl\theta\br\di s.
\end{split}
\end{align}
Here, $I_{12}=I_{21}=0$ for $\theta=0$ and $\theta=\frac{\pi}{2}$, i.e.\ for the straight horizontal and the straight vertical line, respectively. Hence, the eigenvalues of the stability matrix are given by $I_{11}$ and $I_{22}$ in this case whose value are given by $$I_{11}=2R_cf_l(R_c), \quad I_{22}=2\int_0^{R_c} f_s(s)\di s$$ for $\theta=0$ and $$I_{11}=2\int_0^{R_c} f_l(s)\di s,\quad I_{22}=2R_c f_s(R_c)$$ for $\theta=\frac{\pi}{2}$. This leads to the necessary conditions for high-wave number stability  for $\theta=\frac{\pi}{2}$ in \eqref{eq:condvertlinehighwave}, while due to Assumption \ref{ass:newforcegeneral} we obtain instability of the straight horizontal line. 

Note that for any $\theta\in[0,\pi)$ the eigenvalues $\lambda_k, k=1,2,$  are either real or complex conjugated and thus the sum and the product of $\lambda_k$ are real. 
The condition $\Re(\lambda_k)\leq 0, k=1,2,$ is equivalent to $\text{trace}(M)=\lambda_1+\lambda_2\leq 0$ and $\text{det}(M)=\lambda_1 \lambda_2\geq 0$. Hence, we require for the stability of the rotated straight line:
\begin{align}\label{eq:rotatedlinecond}
I_{11}+I_{22}\leq 0\quad\text{and}\quad I_{11}I_{22}-I_{12}I_{21}\geq 0.
\end{align} 

For showing the instability of the rotated straight line with angle of rotation $\theta\in (0,\pi)\setminus [\phi,\pi-\phi]$ for some $\phi\in(0,\frac{\pi}{2})$ we show that the two conditions in \eqref{eq:rotatedlinecond} cannot be satisfied simultaneously in this case.

\begin{theorem}
	For $N\in \N$ sufficiently large and in the continuum limit $N\to\infty$, the single  straight line  \eqref{eq:straightlineansatzgenerallength} where the angle of rotation $\theta\in(0,\pi)\setminus [\phi,\pi-\phi]$ for some $\phi\in(0,\frac{\pi}{2})$ satisfies \eqref{eq:straightlineanglerotation} is an unstable steady state to the particle model \eqref{eq:particlemodelperiodic} for any  force coefficients $f_s$ and $f_l$  satisfying the general conditions for force coefficients in Assumption \ref{ass:newforcegeneral} and the conditions in \eqref{eq:condvertlinehighwave}.
\end{theorem}
\begin{proof}
	Note that  we have
	\begin{align*}
	\int_0^{R_c} f_s'(s)s\sin^2\bl\theta\br\di s&=\sin^2\bl\theta\br\bl f_s(R_c)R_c-\int_0^{R_c} f_s(s)\di s\br
	\end{align*}
	by integration by parts. For $\theta=0$ and $f_l(R_c)=0$ we have
	\begin{align*}
	I_{11}+I_{22}=2R_c f_l(R_c)+2\int_0^{R_c} f_s(s) \di s>0,
	\end{align*}
	while for $\theta=\frac{\pi}{2}$ we have
	\begin{align*}
	I_{11}+I_{22}=2R_c f_s(R_c)+2\int_0^{R_c} f_l(s) \di s\leq 0
	\end{align*}
	by \eqref{eq:condvertlinehighwave}.
	Hence, there exists $\phi\in (0,\frac{\pi}{2})$ such that $I_{11}+I_{22}>0$ on $(0,\phi)$. Since $\cos^2(\theta)=\cos^2(\pi-\theta)$ and $\sin^2(\theta)=\sin^2(\pi-\theta)$ we have $I_{11}+I_{22}>0$ on $(\pi-\phi,\pi)$, implying that stability may only be possible on $[\phi,\pi-\phi]$.
\end{proof}

\section{Stability of vertical lines for particular force coefficients}\label{sec:verticallines}

We have investigated the high-wave number stability for $m\to \infty$ in Section \ref{sec:stationarystates}. Since only vertical straight lines for the considered spatially homogeneous tensor field $T$ in \eqref{eq:tensorfieldhom} can lead to stable steady states for any $N\in\N$ we restrict ourselves to vertical straight lines in the sequel. As a next step towards proving stability we consider the stability for fixed modes $m\in \N$ now. 

Due to the  form of the eigenvalues  in \eqref{eq:eigenvalueslinegeneral} no general stability result for the single straight vertical line for the particle system \eqref{eq:particlemodelperiodic} with arbitrary force coefficients $f_s$ and $f_l$ satisfying Assumption \ref{ass:newforcegeneral} can be derived. Thus, additional assumptions on the force coefficients are necessary. 

\subsection{Linear force coefficients}\label{sec:linearcase}

To study the stability of the single straight vertical line for any $N\in \N$, we consider linear force coefficients satisfying Assumption \ref{ass:newforcegeneral}. To guarantee that the force coefficient is differentiable, required for using the results from Section \ref{sec:stationarystates} we consider the differentiable adaptation \eqref{eq:forcecutoff}  for a  given linear force coefficient, leading to a linear force coefficient on $[0,R_c-\epsilon]$ for some $\epsilon>0$, a cubic polynomial on $(R_c-\epsilon,R_c)$ and the constant zero function for $|d|\geq R_c$. This leads to the following conditions:
\begin{assumption}\label{ass:linearforcecoeff}
	For any $\epsilon>0$ with $\epsilon\ll R_c$, we assume that the force coefficients are linear on $[0,R_c-\epsilon]$, i.e.\
	\begin{align}\label{eq:linearcoefffunctions}
	\begin{split}
	f_l^\epsilon(|d|)&:=\begin{cases}a_l|d|+b_l, & |d|\in[0,R_c-\epsilon],\\(2 b_l + 2 R_c a_l - a_l \epsilon)\frac{(|d|-R_c)^3}{\epsilon^3}+(3 b_l+ 3 R_c a_l - 2 a_l\epsilon) \frac{\left(|d|-R_c \right)^2}{\epsilon^2},  & |d| \in (R_c-\epsilon,R_c), \\ 0, & |d| \geq R_c, \end{cases}\\ 
	f_s^\epsilon(|d|)&:=\begin{cases}a_s|d|+b_s, & |d|\in[0,R_c-\epsilon],\\(2 b_s + 2 R_c a_s - a_s \epsilon)\frac{(|d|-R_c)^3}{\epsilon^3}+(3 b_s+ 3 R_c a_s - 2 a_s\epsilon) \frac{\left(|d|-R_c \right)^2}{\epsilon^2},  & |d| \in (R_c-\epsilon,R_c), \\ 0, & |d| \geq R_c, \end{cases}
	\end{split}
	\end{align}
	for constants $a_l, a_s, b_l,b_s$. Since $f_l^\epsilon$ and $f_s^\epsilon$ are short-range repulsive, we require 
	\begin{align*}	
	b_l>0,\quad b_s>0.
	\end{align*}
	Besides, for physically realistic force coefficients the absolute values of $f_l^\epsilon$ and $f_s^\epsilon$ are decaying, i.e.\
	\begin{align*}
	a_l<0,\quad a_s<0.
	\end{align*}
\end{assumption}
Note that for the short-range repulsive, long-range attractive force coefficient $f_l$, we have $a_l R_c +b_l<0$ and in particular $a_l R_c +b_l$ is of order $\mathcal{O}(1)$. Hence, the adaptation  $f_l^\epsilon$ of $f_l$ for $f_l$ linear is not negligible. However, due to the concentration of particles along a straight vertical line the adaptation $f_l^\epsilon$ acting along the vertical axis does not influence the overall dynamics provided $0< \epsilon\ll R_c$. For the force coefficient $f_s$, the adaption $f_s^\epsilon$ of $f_s$ is negligible if  $a_s R_c +b_s$ is of order $\mathcal{O}(\epsilon)$ and also results in the same stability/instability properties numerically, see Section \ref{sec:numericalresultsstability}. If $a_s R_c +b_s$ is of order $\mathcal{O}(1)$, then the adaptation is not negligible, but the numerical results in Section \ref{sec:numericalresultsstability} illustrate that we obtain the same stability/instability results for $f_s^\epsilon$ and $f_s$. 

\begin{remark}\label{rem:linforcecoeff}
	Note that the modelling assumptions in Assumption \ref{ass:newforcegeneral} and Assumption~\ref{ass:linearforcecoeff} can be applied to  linear repulsive and attractive force coefficients $f_R^\epsilon$ and $f_A^\epsilon$ as in \eqref{eq:linearcoefffunctions}, where  the total force of the form  \eqref{eq:totalforcenewhom} consists of repulsion and attraction forces. That is, for $\epsilon>0$ we define
	\begin{align}\label{eq:linearforcecoeffkc}
	\begin{split}
	f_R^\epsilon(|d|)&:=\begin{cases}a_R|d|+b_R, & |d|\in[0,R_c-\epsilon],\\(2 b_R + 2 R_c a_R - a_R \epsilon)\frac{(|d|-R_c)^3}{\epsilon^3}+(3 b_R+ 3 R_c a_R - 2 a_R\epsilon) \frac{\left(|d|-R_c \right)^2}{\epsilon^2},  & |d| \in (R_c-\epsilon,R_c), \\ 0, & |d| \geq R_c, \end{cases}\\ 
	f_A^\epsilon(|d|)&:=\begin{cases}a_A|d|+b_A, & |d|\in[0,R_c-\epsilon],\\(2 b_A + 2 R_c a_A - a_A \epsilon)\frac{(|d|-R_c)^3}{\epsilon^3}+(3 b_A+ 3 R_c a_A - 2 a_A\epsilon) \frac{\left(|d|-R_c \right)^2}{\epsilon^2},  & |d| \in (R_c-\epsilon,R_c), \\ 0, & |d| \geq R_c, \end{cases}
	\end{split}
	\end{align}
	for constants $a_A, a_R, b_A,b_R$ and  we require 
	\begin{align*}
	f_R^\epsilon\geq 0\quad \text{and} \quad f_A^\epsilon \leq 0
	\end{align*}
	for all $\epsilon>0$ with $\epsilon \ll R_c$,
	implying
	\begin{align}\label{eq:linearcoeffrepattr}
	a_Rs+b_R\geq 0\quad \text{and} \quad a_As+b_A\leq 0 \quad \text{for} \quad s\in [0,R_c],
	\end{align}
	and in particular
	\begin{align}\label{eq:linearcoeffb}
	b_R>0\quad \text{and}\quad b_A<0.
	\end{align}
	For realistic interaction force coefficients $f_R^\epsilon$ and $f_A^\epsilon$ we assume that their absolute values decrease as the distance between the particles increases, implying
	\begin{align}\label{eq:linearcoeffa}
	a_R<0\quad \text{and} \quad a_A>0
	\end{align}
	by the definition of $f_R^\epsilon$ and $f_A^\epsilon$ in \eqref{eq:linearforcecoeffkc} and by the condition for $b_R$ and $b_A$ in \eqref{eq:linearcoeffb}. Combining the assumptions on $a_A,a_R$ in \eqref{eq:linearcoeffa} and $b_A,b_R$ in \eqref{eq:linearcoeffb} condition \eqref{eq:linearcoeffrepattr} reduces to 
	\begin{align*}
	a_RR_c+b_R\geq 0\quad \text{and} \quad a_AR_c+b_A \leq 0.
	\end{align*}
	Further we assume that $f_A^\epsilon+f_R^\epsilon$ is short-range repulsive, long-range attractive for any $\epsilon>0$ with $\epsilon\ll R_c$, i.e.\ $$(f_A^\epsilon+f_R^\epsilon)(0)=b_A+b_R>0,\qquad (f_A^\epsilon+f_R^\epsilon)(R_c-\epsilon)=(a_A+a_R)(R_c-\epsilon)+b_A+b_R<0$$ for all $0<\epsilon\ll R_c$  implying
	\begin{align}\label{eq:linearcoeffsum}
	b_A+b_R>0\quad \text{and}\quad a_A+a_R<0.
	\end{align}
	For any $\epsilon>0$, the force coefficient $\chi f_A^\epsilon+f_R^\epsilon$ is purely repulsive along $s$ on $[0,R_c-\epsilon]$ if $\chi\in[0,1]$ is sufficiently small since $f_R^\epsilon$ is  repulsive. Note that \eqref{eq:linearcoeffsum} implies
	\begin{align*}
	\chi a_A+a_R<0, \quad\chi b_A+b_R>0 \quad\text{for all}\quad\chi\in[0,1]
	\end{align*}
	by the positivity of $b_R$ and by the negativity of $a_R$ in \eqref{eq:linearcoeffb} and \eqref{eq:linearcoeffa}, respectively. Since 
	\begin{align*}
	f_l^\epsilon(|d|)= f_A^\epsilon(|d|)+f_R^\epsilon(|d|)&=\bl a_A+a_R\br|d|+b_A+b_R
	\end{align*}
	and
	\begin{align*}
	f_s^\epsilon(|d|)=  \chi f_A^\epsilon(|d|)+f_R^\epsilon(|d|)=\bl \chi a_A+a_R\br|d|+\chi b_A+b_R
	\end{align*}
	for $|d|\in [0,R_c-\epsilon],$ we  have 
	\begin{align*}
	a_l=a_A+a_R<0,\quad a_s= \chi a_A+a_R<0,\quad b_l=b_A+b_R>0,\quad b_s=\chi b_A+b_R>0
	\end{align*}
	as in Assumption \ref{ass:linearforcecoeff}.
\end{remark}

For investigating the stability of the straight line for any $N\in\N$, we consider the real parts of the eigenvalues in \eqref{eq:eigenvaluesreallinegeneral}, i.e.\
\begin{align*}
\Re(\lambda_1(m))&=2\int_{0}^{R_c}  f_l^\epsilon(s)\bl 1-\cos\bl - 2\pi ms\br\br\di s,\\
\Re(\lambda_2(m))&=2\int_{0}^{R_c}
\bl  f_s^\epsilon(s)+ s(f_s^\epsilon)'(s)
\br\bl 1-\cos\bl -2\pi ms\br\br\di s.
\end{align*}
Note that the coefficient functions of the integrands in the definition of the eigenvalues  are given by
\begin{align*}
f_l^\epsilon(s)= a_ls+b_l,\quad f_s^\epsilon(s)+s(f_s^\epsilon)'(s)= 2a_s s+b_s
\end{align*}
for $s\in [0,R_c-\epsilon]$ with $a_l,a_s<0,~b_l,b_s>0$ and 
\begin{align*}
f_l^\epsilon(s)&= (2 b_l + 2 R_c a_l - a_l \epsilon)\frac{(s-R_c)^3}{\epsilon^3}+(3 b_l+ 3 R_c a_l - 2 a_l\epsilon) \frac{\left(s-R_c \right)^2}{\epsilon^2},\\ f_s^\epsilon(s)+s(f_s^\epsilon)'(s)&= (2 b_s + 2 R_c a_s - a_s \epsilon)\frac{(s-R_c)^3}{\epsilon^3}+(3 b_s+ 3 R_c a_s - 2 a_s\epsilon) \frac{\left(s-R_c \right)^2}{\epsilon^2}\\&\qquad +3(2 b_s + 2 R_c a_s - a_s \epsilon)\frac{s(s-R_c)^2}{\epsilon^3}+2(3 b_s+ 3 R_c a_s - 2 a_s\epsilon) \frac{s \left(s-R_c \right)}{\epsilon^2},
\end{align*}
for $s \in [R_c-\epsilon, R_c]$
by Assumption \ref{ass:linearforcecoeff}. Since $f_l^\epsilon(s)$ and $f_s^\epsilon(s)+s(f_s^\epsilon)'(s)$ are bounded on $[R_c-\epsilon,R_c]$, we obtain
\begin{align}\label{eq:realeigenvalues}
\begin{split}
\Re(\lambda_1(m))&=2\int_{0}^{R_c-\epsilon}  (a_ls+b_l)\bl 1-\cos\bl - 2\pi ms\br\br\di s+\mathcal{O}(\epsilon),\\
\Re(\lambda_2(m))&= 2\int_{0}^{R_c-\epsilon}
	\bl  2a_s s+b_s
	\br\bl 1-\cos\bl -2\pi ms\br\br\di s \\&\quad  +\frac{12(b_s +  R_c a_s)}{\epsilon} \int_{R_c-\epsilon}^{R_c} \left(\frac{s(s-R_c)^2}{\epsilon^2}+ \frac{s \left(s-R_c \right)}{\epsilon}\right)\bl 1-\cos\bl - 2\pi ms\br\br \di s+\mathcal{O}(\epsilon).
\end{split}
\end{align}
Note that $s (f_s^{\epsilon})'$ is of order $\mathcal{O}(1/\epsilon)$ on $[R_c-\epsilon,R_c]$ and hence the integral over $[R_c-\epsilon,R_c]$ also contributes to the leading order term for $ \Re(\lambda_2(m))$.
Here, $f_l^\epsilon(s)$ and $f_s^\epsilon(s)+s(f_s^\epsilon)'(s)$ are linear functions on $[0,R_c-\epsilon]$ of the form $f|_{[0,R_c-\epsilon]}\to\R, s\mapsto as+b$ for constants $a<0$ and $b>0$. In particular,  $\Re(\lambda_1)$ and the first term in $\Re(\lambda_2)$ are of the form
	\begin{align}\label{eq:eigenvalueslinearforceparameter}
	2\int_0^{R_c-\epsilon} \bl a_ks+b_k\br\bl 1-\cos\bl 2\pi ms\br\br\di s
	\end{align}
where
\begin{align}\label{eq:linearforceparameter}
a_1=a_l,\quad a_2=2a_s,\quad b_1=b_l,\quad b_2= b_s.
\end{align}
For ease of notation we drop the indices of $a_k$ and $b_k$ in the sequel. Note that 
\begin{align}\label{eq:straightlinestabilitylinear}
\begin{split}
&\int_0^{R_c-\epsilon} \bl as+b\br\bl 1-\cos\bl 2\pi ms\br\br\di s\\
&=\frac{2\pi m\bl \pi m (R_c-\epsilon)\bl a (R_c-\epsilon)+2b\br-\bl a (R_c-\epsilon)+b\br\sin\bl 2\pi m (R_c-\epsilon)\br\br}{4\pi^2 m^2}\\&\quad +\frac{a -a\cos\bl 2\pi m(R_c-\epsilon)\br}{4\pi^2 m^2}.
\end{split}
\end{align}
In the limit $m\to\infty$, all  terms in the second line of \eqref{eq:straightlinestabilitylinear} vanish except for the first term. Since $R_c>0$, we require
\begin{align*}
a \leq -b\frac{2}{R_c -\epsilon}
\end{align*}
for high-wave number stability for any $\epsilon>0$ with $\epsilon\ll R_c$. In particular, this condition is consistent  with the necessary condition for high-wave number stability  in Proposition \ref{prop:highwavestraightline} for arbitrary force coefficients $f_s^\epsilon$ and $f_l^\epsilon$ satisfying Assumption \ref{ass:newforcegeneral}. In the limit $\epsilon\to 0$, it reduces to
\begin{align}\label{eq:straightlinehighwaveansatz}
a \leq -b\frac{2}{R_c }
\end{align}
Since $R_c\in(0,0.5]$ and $b>0$, \eqref{eq:straightlinehighwaveansatz} implies that $a<0$ is necessary for high-wave number stability. Hence, we can assume
\begin{align*}
a<0\quad\text{and}\quad b>0
\end{align*}
in the sequel.

\begin{lemma}\label{lem:stabilitystraightvertline}
	Let  $b>0$ and $R_c\in(0,0.5]$. For $\epsilon>0$, set
	\begin{align}\label{eq:straightlinestabilityhelpfunctions}
	\begin{split}
	g_\epsilon(m)&:=2\pi m \bl\pi m (R_c-\epsilon)^2-(R_c-\epsilon)\sin(2\pi m(R_c-\epsilon))\br+1-\cos(2\pi m (R_c-\epsilon)),\\
	h_\epsilon(m)&:=2\pi m \bl2\pi m (R_c-\epsilon)-\sin(2\pi m (R_c-\epsilon))\br.
	\end{split}
	\end{align}
	Then, 
	\begin{align}\label{eq:straightlinestabilitylinearhelp}
	\int_0^{R_c-\epsilon} \bl as+b\br\bl 1-\cos\bl 2\pi ms\br\br\di s\leq 0
	\end{align}
	is satisfied for all $m\in\N$ and all $\epsilon>0$ with $\epsilon\ll R_c$ if and only if $a\leq a_0$ with
	\begin{align}\label{eq:defa0}
	a_0:=-b\max_{m\in\N}  \frac{h_\epsilon(m)}{g_\epsilon(m)} \leq -\frac{2b}{R_c}\leq 0.
	\end{align}
\end{lemma}
\begin{proof}
	Note that the numerator of  \eqref{eq:straightlinestabilitylinear} is of the form $ag_\epsilon(m)+bh_\epsilon(m)$ for functions $g_\epsilon$ and $h_\epsilon$, defined in \eqref{eq:straightlinestabilityhelpfunctions}. Condition \eqref{eq:straightlinestabilitylinearhelp} is equivalent to $$a\leq -b\frac{h_\epsilon(m)}{g_\epsilon(m)}\quad \text{for all~}m\in\N.$$
	Herein, $h_\epsilon(m)\geq 0$ for all $m\geq 0$ since $h_\epsilon$ is an increasing function. Further note that 
	\begin{align*}
	g_\epsilon'(m)&=2\pi\bl \pi m (R_c-\epsilon)^2-(R_c-\epsilon)\sin(2\pi m (R_c-\epsilon))\br\\&\quad+2\pi m \bl \pi (R_c-\epsilon)^2-2\pi (R_c-\epsilon)^2\cos(2\pi m (R_c-\epsilon))\br+2\pi (R_c-\epsilon)\sin(2\pi m (R_c-\epsilon))\\
	&=4\pi^2 m(R_c-\epsilon)^2\bl 1-\cos(2\pi m (R_c-\epsilon))\br
	\end{align*}
	is nonnegative implying that $g_\epsilon$ is an increasing function with $g_\epsilon(0)=0$. In particular, $g_\epsilon$ and $h_\epsilon$ are nonnegative functions for all $m\in\N$. Hence, \eqref{eq:straightlinestabilitylinearhelp} is satisfied for all $m\in\N$ if and only if $a<a_0$. Note that 
	\begin{align*}
	\lim_{m\to \infty} \frac{h_\epsilon(m)}{g_\epsilon(m)}=\frac{2}{R_c-\epsilon},
	\end{align*}
	implying that $$\sup_{m\in\N} \frac{h_\epsilon(m)}{g_\epsilon(m)}\in\R$$ by the nonnegativity and continuity of $g_\epsilon$ and $h_\epsilon$. 
	
	Let $R_c\in(0,0.5]$ and $\epsilon>0$. We have
	\begin{align*}
	\max_{m\in\N} \frac{h_\epsilon(m)}{g_\epsilon(m)}&\geq 
	\frac{2\pi m \bl2\pi m (R_c-\epsilon)-\sin(2\pi m (R_c-\epsilon))\br}{2\pi m \bl\pi m (R_c-\epsilon)^2-(R_c-\epsilon)\sin(2\pi m(R_c-\epsilon))\br+2}\\&=\frac{2}{R_c-\epsilon}\bl 1+\frac{\pi m\sin(2\pi m(R_c-\epsilon))-\frac{2}{R_c-\epsilon}}{2\pi m \bl\pi m (R_c-\epsilon)-\sin(2\pi m(R_c-\epsilon))\br+\frac{2}{R_c-\epsilon}}\br
	\end{align*}
	for all $m\in\N$. 
	Since $R_c-\epsilon\in(0,0.5)$ there exists $m\in \N$ such that $\pi m\sin(2\pi m(R_c-\epsilon))-\frac{2}{R_c-\epsilon}>0$ and hence
	\begin{align*}
	\max_{m\in\N} \frac{h_\epsilon(m)}{g_\epsilon(m)}>\frac{2}{R_c-\epsilon}>\frac{2}{R_c}
	\end{align*}
	for $\epsilon>0$ with $\epsilon\ll R_c$. For $R_c \in (0,0.5)$, we obtain
	\begin{align*}
	\lim_{\epsilon\to 0}	\max_{m\in\N} \frac{h_\epsilon(m)}{g_\epsilon(m)}>\frac{2}{R_c}.
	\end{align*}
	For $R_c=0.5$, we have
	\begin{align*}
	\lim_{\epsilon\to 0}\frac{h_\epsilon(m)}{g_\epsilon(m)}&=\begin{cases}\frac{2}{R_c}, &\quad m \text{~even},\\
	\frac{4\pi^2  R_c}{2\pi^2  R_c^2+\frac{2}{m^2}}<\frac{2}{R_c}, &\quad m\text{~odd},\end{cases}
	\end{align*}
	implying that 
	\begin{align*}
	\lim_{\epsilon\to 0}\max_{m\in\N} \frac{h_\epsilon(m)}{g_\epsilon(m)}= \frac{2}{R_c}.
	\end{align*}
	Hence, $a\leq a_0$  is equivalent to the necessary condition \eqref{eq:straightlinehighwaveansatz} for high-wave number stability for $R_c=0.5$. 
\end{proof}

\begin{remark}\label{rem:stabilitystraightvertline}
	For the stability of line patterns with force coefficients $f_s^\epsilon,f_l^\epsilon$ of the form \eqref{eq:linearcoefffunctions}, we require $\Re(\lambda_k(m))\leq 0$ for $k=1,2$ for the real parts of the eigenvalues $\Re(\lambda_k(m))$, $k=1,2$, in \eqref{eq:realeigenvalues}. Note that the nonnegativity of the leading-order term of $\Re(\lambda_1)$  which is of the form \eqref{eq:eigenvalueslinearforceparameter} is equivalent to condition \eqref{eq:straightlinestabilitylinearhelp} in Lemma \ref{lem:stabilitystraightvertline}. Similarly, the nonnegativity of the first term in \eqref{eq:realeigenvalues} which is also of the form \eqref{eq:eigenvalueslinearforceparameter} is equivalent to condition \eqref{eq:straightlinestabilitylinearhelp} in Lemma \ref{lem:stabilitystraightvertline}. 
	
	From the proof of Lemma \ref{lem:stabilitystraightvertline} it follows that
	\begin{align}\label{eq:maxquotientapprox}
	\frac{2}{R_c}\leq \max_{m\in\N} \frac{h_\epsilon(m)}{g_\epsilon(m)}.
	\end{align}
	The inequality in \eqref{eq:maxquotientapprox} is strict for $R_c -\epsilon\in(0,0.5)$, i.e.\ a necessary condition for 
		\eqref{eq:straightlinestabilitylinearhelp} to hold for $R_c -\epsilon \in(0,0.5)$ is given by
		\begin{align*}
		a<-\frac{2b}{R_c}.
		\end{align*}
	For $R_c=0.5$ and $\epsilon\to 0$,  condition \eqref{eq:straightlinestabilitylinearhelp} holds for any $a<0$ satisfying the necessary condition \eqref{eq:straightlinehighwaveansatz} for high-wave number stability.  If the necessary condition \eqref{eq:straightlinehighwaveansatz} for high-wave number stability is satisfied with equality, i.e.\ $a=-\frac{2b}{R_c}$, the leading order term of the left-hand side of \eqref{eq:realeigenvalues} vanishes for $\epsilon\to 0$ in the high-wave limit and lower order terms have to be considered.
\end{remark}

\begin{figure}[htbp]
	\centering
	\subfloat[$\frac{h_0(m)}{g_0(m)}$]{\includegraphics[width=0.45\textwidth]{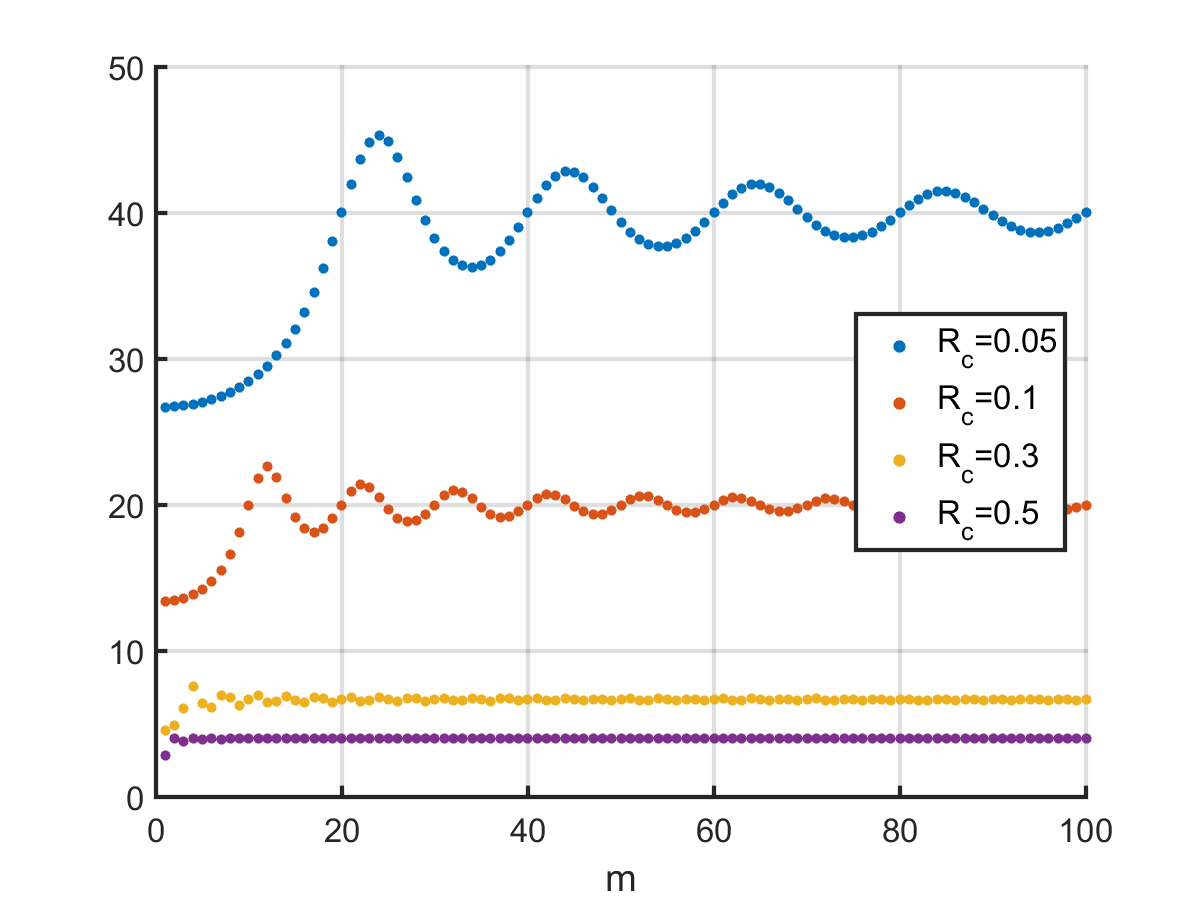}\label{fig:linearcoeffquotientmdep}}
	\subfloat[$\max_{m\in\N}\frac{h_0(m)}{g_0(m)}$]{
		\includegraphics[width=0.45\textwidth]{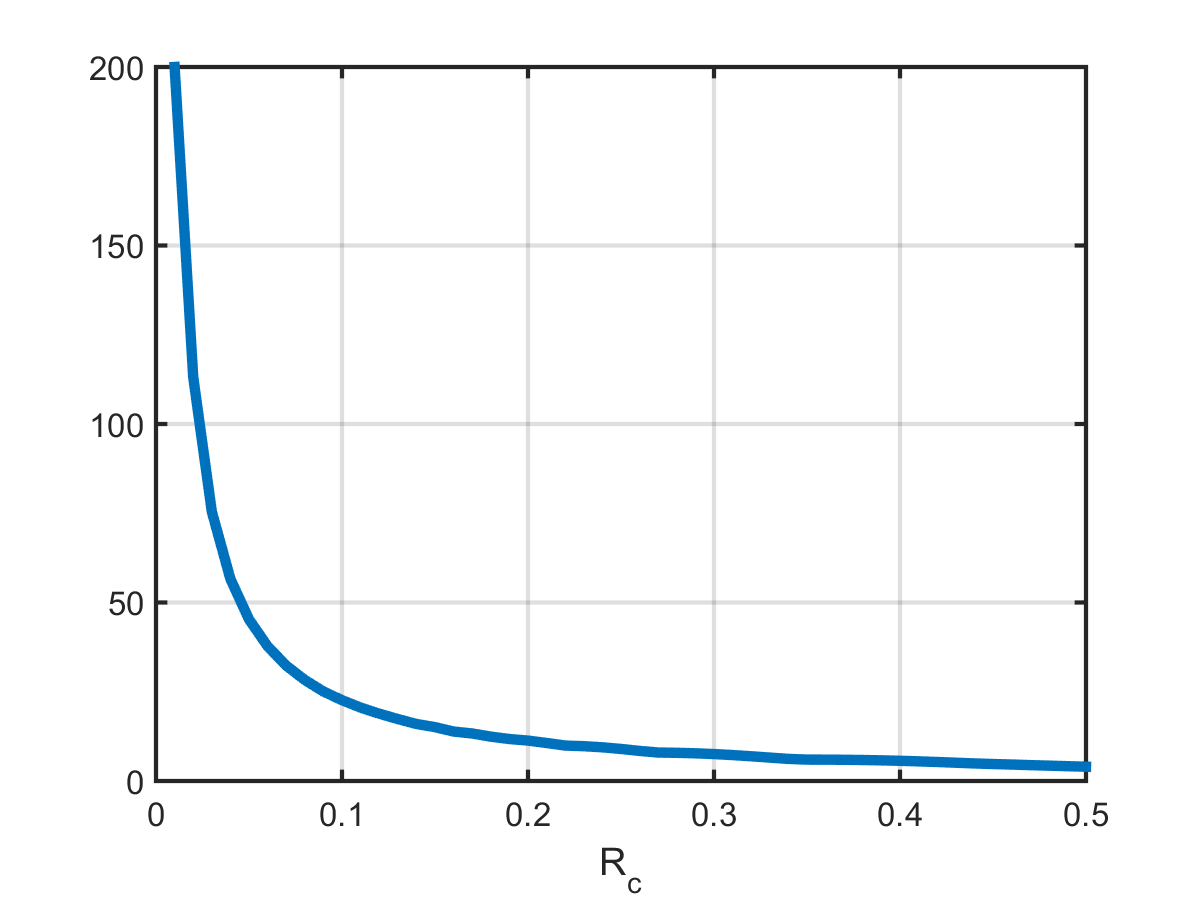}\label{fig:linearcoeffquotientmax}}\\
	\subfloat[$R_c\frac{h_0(m)}{g_0(m)}$]{\includegraphics[width=0.45\textwidth]{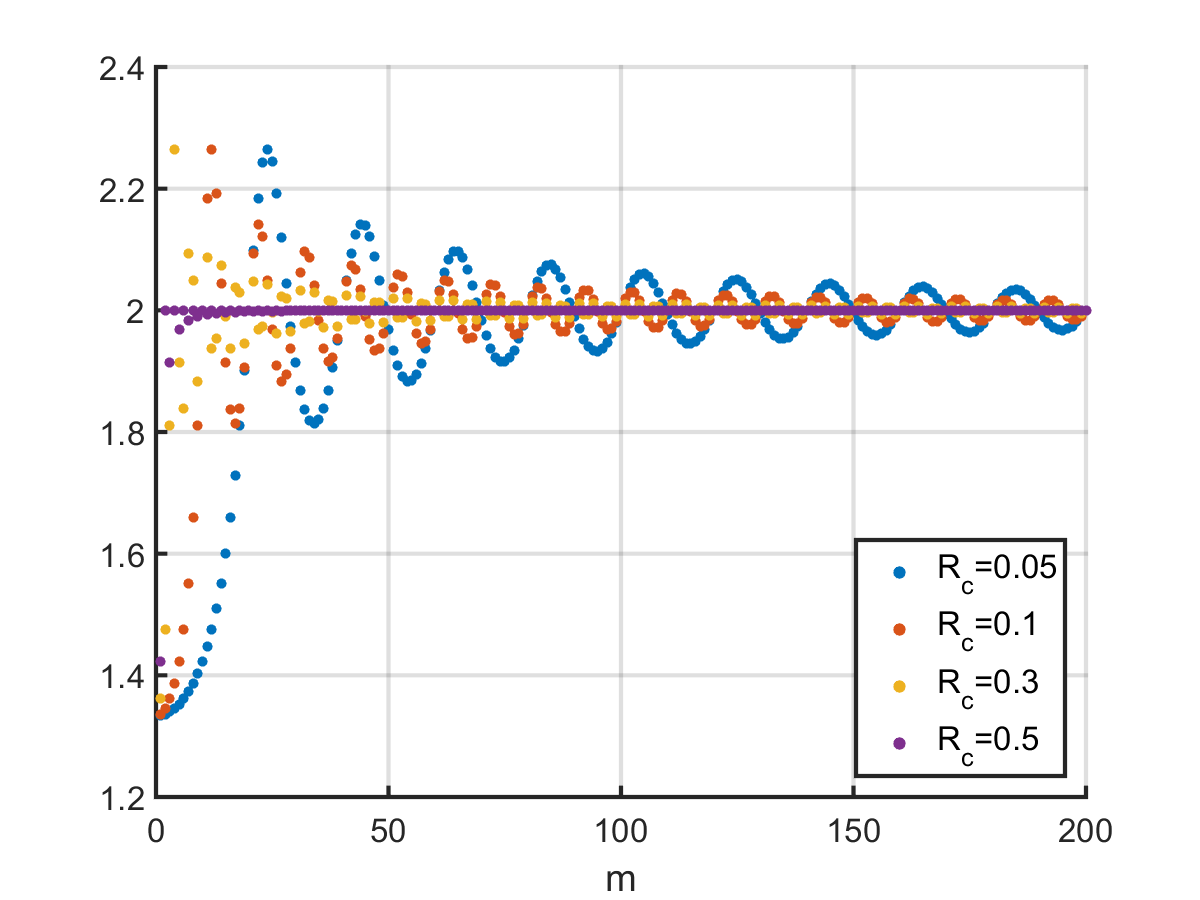}\label{fig:linearcoeffquotientmdeprescalerc}}
	\subfloat[$R_c \max_{m\in\N}\frac{h_0(m)}{g_0(m)}$]{
		\includegraphics[width=0.45\textwidth]{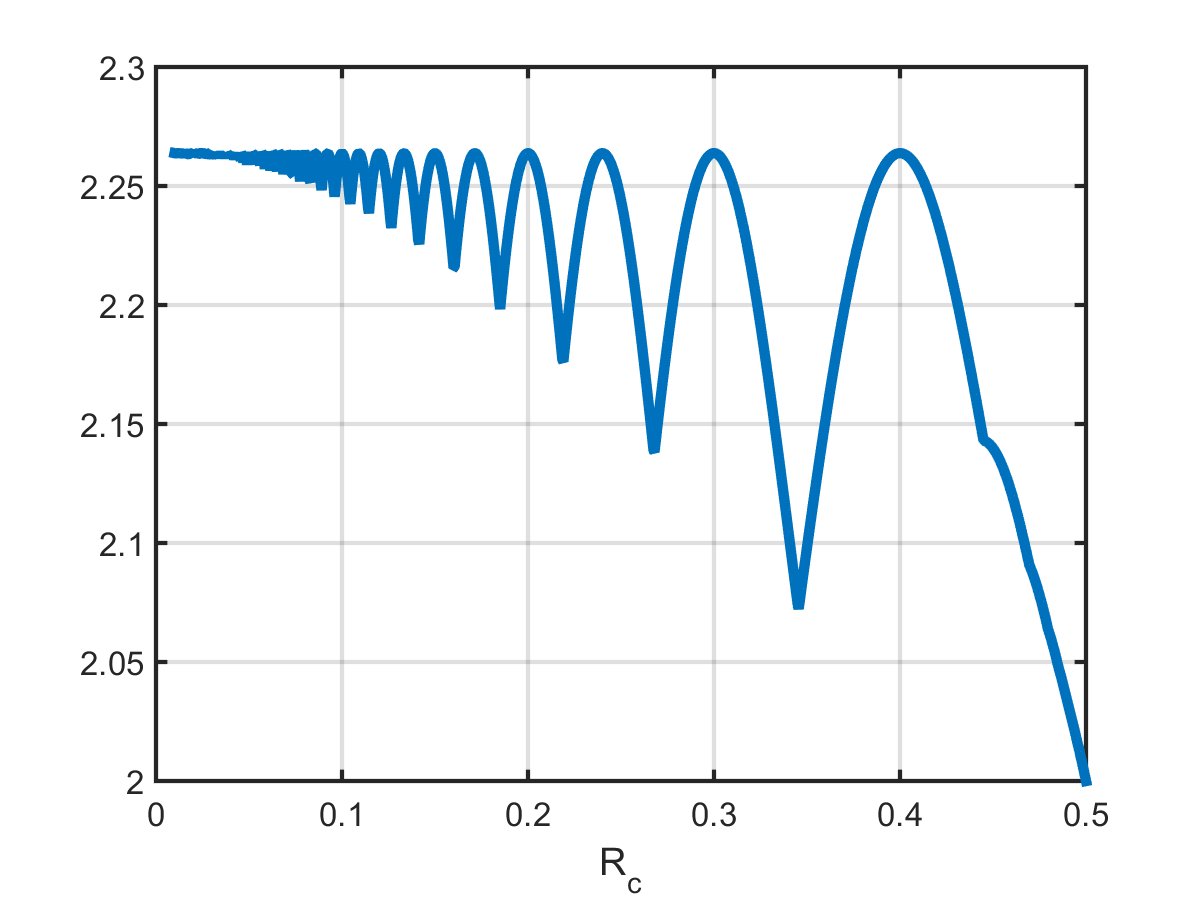}\label{fig:linearcoeffquotientmaxrescalerc}}        
	\caption{Scaling factor of $a_0$ in \eqref{eq:defa0} as a functions of $R_c$ where $g_0,h_0$ are defined in \eqref{eq:straightlinestabilityhelpfunctions}.}\label{fig:linearcoeffgeneral}
\end{figure}

In Figure \ref{fig:linearcoeffgeneral}, we investigate the scaling factor of $a_0$, defined in \eqref{eq:defa0}, numerically. In Figure~\subref*{fig:linearcoeffquotientmdep} the quotient  $h_0/g_0$ is shown as a function of $m\in\N$ for different values of the cutoff radius $R_c$. Note that for smaller values of $R_c$,  the  maximum of  $h_0/g_0$ gets larger as shown in Figure~\subref*{fig:linearcoeffquotientmax}. In Figures \subref*{fig:linearcoeffquotientmdeprescalerc} and \subref*{fig:linearcoeffquotientmaxrescalerc} we consider the quotient  $h_0/g_0$ scaled by $R_c$. Figure \subref*{fig:linearcoeffquotientmdeprescalerc} shows that $R_c h_0/g_0\to 2$ as $m\to\infty$, independently of the value of $R_c$, and that the maximum of $R_c h_0/g_0$ is obtained for smaller values of $m\in\N$ in general. The value of $$R_c \max_{m\in\N}\frac{h_0(m)}{g_0(m)}$$
is shown in Figure \subref*{fig:linearcoeffquotientmaxrescalerc} as a function of $R_c$. In particular, the scaled maximum is larger than 2 if and only if $R_c\in(0,0.5)$ and is equal to $2$ for $R_c=0.5$. Hence, this numerical investigation is consistent with the results in Lemma \ref{lem:stabilitystraightvertline}.

Applying  Lemma \ref{lem:stabilitystraightvertline} to the specific form of the stability conditions for a single straight vertical line leads to the following stability results for the linear force coefficients \eqref{eq:linearcoefffunctions}:

\begin{proposition}\label{prop:stabilityverticalline}
	For $R_c\in(0,0.5)$, the single straight vertical line is an unstable steady state of \eqref{eq:particlemodelperiodic} for any $N\in \N$ sufficiently large and in the continuum limit $N\to \infty$, where the forces are of the form \eqref{eq:totalforcehom} for any linear coefficient functions $f_s^\epsilon,f_l^\epsilon$ with $0<\epsilon\ll R_c$  such that Assumption~\ref{ass:linearforcecoeff} is satisfied. In particular, the single straight vertical line is an unstable steady state for force coefficients $f_s^\epsilon,f_l^\epsilon$ for $R_c\in(0,0.5)$ in the limit $\epsilon\to 0$. 
\end{proposition}
	\begin{proof}
		Note that the leading order term of $\Re(\lambda_1(m))$ and the first term of $\Re(\lambda_2(m))$ in \eqref{eq:realeigenvalues} are of the form \eqref{eq:eigenvalueslinearforceparameter} with parameters \eqref{eq:linearforceparameter}. For stability we require $\Re(\lambda_k)\leq 0$ for $k=1,2$. 	\\
		Let us consider the nonnegativity of $\Re(\lambda_2(m))$ in \eqref{eq:realeigenvalues} first. 
		Note that the second leading order term of $\Re(\lambda_2(m))$ in \eqref{eq:realeigenvalues} can be rewritten as
		\begin{align}\label{eq:secondleadingterm}
		\begin{split}
		& 12(b_s+R_ca_s) \int_{-1}^0  ((\epsilon s +R_c)s^2+(\epsilon s +R_c)s)\bl 1-\cos\bl 2\pi m(\epsilon s +R_c)\br\br  \di s \\
		&=-2(b_s+R_ca_s)R_c \bl 1-\cos\bl 2\pi m R_c\br\br +\mathcal{O}(\epsilon).
		\end{split}
		\end{align}	
		Hence, $\Re(\lambda_2(m))$ is of the form
		\begin{align*}
		\Re(\lambda_2(m))=\frac{2}{4\pi^2m^2}(2a_s g_0(m)+b_s h_0(m))-2(b_s+R_ca_s)R_c \bl 1-\cos\bl 2\pi m R_c\br\br+\mathcal{O}(\epsilon)
		\end{align*}
		by \eqref{eq:straightlinestabilitylinear} where $g_0,h_0$ are defined in \eqref{eq:straightlinestabilityhelpfunctions}. For the nonnegativity of the leading order term of $\Re(\lambda_2(m))$ we require
		\begin{align*}
		a_s \bl\frac{2g_0(m)}{4\pi^2m^2}-R_c^2\bl 1-\cos\bl 2\pi m R_c\br\br \br+b_s \bl \frac{h_0(m)}{4\pi^2m^2}-R_c \bl 1-\cos\bl 2\pi m R_c\br\br\br\leq 0,
		\end{align*}
		which can be rewritten as 
		\begin{align*}
		a_s \tilde{g}_0(m)+b_s \tilde{h}_0(m)\leq 0
		\end{align*}
		where 
		\begin{align*}
		\begin{split}
		\tilde{g}_0(m)&:=2\pi m \bl 2\pi m R_c^2 \cos(2\pi m R_c)-2R_c\sin(2\pi mR_c)\br+2-2\cos(2\pi m R_c),\\
		\tilde{h}_0(m)&:=2\pi m \bl2\pi m R_c\cos(2\pi m R_c)-\sin(2\pi m R_c)\br.
		\end{split}
		\end{align*}
		For $m$ sufficiently large, we have
		\begin{align*}
		a_s \tilde{g}_0(m)+b_s \tilde{h}_0(m)=4\pi^2m^2 R_c^2\cos(2\pi m R_c) a_s + 4\pi^2m^2 R_c\cos(2\pi m R_c) b_s+\mathcal{O}(m)
		\end{align*}
		and by only considering the leading order term we obtain the condition  
		\begin{align*}
		R_c\cos(2\pi m R_c) a_s +  \cos(2\pi m R_c) b_s \leq 0.
		\end{align*}
		Note that there exist infinitely many $m\in\N$ such that $\cos(2\pi m R_c)> 0$ and such that $\cos(2\pi m R_c)<0$, independently of the choice of $R_c \in (0,0.5)$. Hence, we can conclude $a_sR_c+b_s=0$. In this case, the second leading order term of $\Re(\lambda_2(m))$ vanishes by \eqref{eq:secondleadingterm} and thus, it is sufficient to consider $\Re(\lambda_1(m))$ and the first term of $\Re(\lambda_2(m))$ in \eqref{eq:realeigenvalues}.
		Applying Lemma \ref{lem:stabilitystraightvertline} together with Remark \ref{rem:stabilitystraightvertline} for $R_c\in(0,0.5)$ to the linear force coefficients $f_l^\epsilon, f_s^\epsilon$ in \eqref{eq:linearcoefffunctions} results in the stability conditions
		\begin{align}\label{eq:condstabilitylincoeffrcsmall}
		a_l<-\frac{2b_l}{R_c}\quad\text{and}\quad a_s<-\frac{b_s}{R_c}
		\end{align}
		for any $\epsilon>0$ which are necessary for the nonnegativity of $\Re(\lambda_1(m))$ and the (first) leading order term of $\Re(\lambda_2(m))$. 
		Hence, the single straight vertical line is unstable for ${R_c\in(0,0.5)}$ and $0<\epsilon\ll R_c$, both in the continuum limit $N\to \infty$ and for any $N\in \N$ sufficiently large. Similarly, we obtain instability of straight vertical line patterns for force coefficients $f_s^\epsilon,f_l^\epsilon$ for $R_c \in (0,0.5)$ in the limit $\epsilon\to 0$.
	\end{proof}

\begin{remark}
		For $R_c=0.5$, we cannot conclude stability/instability of the straight vertical line for the linear force coefficients in \eqref{eq:linearcoefffunctions} with $a_s=-\frac{b_s}{R_c}$, while we can conclude instability for $a_s \neq -\frac{b_s}{R_c}$. 
		To see this, note that for $R_c=0.5$ the calculations in the proof of Proposition~\ref{prop:stabilityverticalline} imply $a_s=-\frac{b_s}{R_c}$ as a necessary condition for stability, from Lemma~\ref{lem:stabilitystraightvertline} we obtain 
		\begin{align*}
		a_l\leq -\frac{2 b_l}{R_c}\quad\text{and}\quad  a_s\leq -\frac{ b_s}{R_c},
		\end{align*}
		and together with the condition that $f_s^\epsilon$ is purely repulsive we get the necessary conditions
		\begin{align}\label{eq:linearcoeffcond}
		a_l\leq -\frac{2 b_l}{R_c}\quad\text{and}\quad a_s= -\frac{b_s}{R_c}
		\end{align}
		for the stability of the straight vertical line.
		Note that the conditions \eqref{eq:linearcoeffcond} are consistent with each other since $a_l,a_s<0$ and $b_l,b_s>0$ by Assumption~\ref{ass:linearforcecoeff} and it is possible to choose the parameters $a_l,a_s,b_l,b_s$ in such a way that both \eqref{eq:linearcoeffcond} and the assumptions on the force coefficients $f_s^\epsilon,f_l^\epsilon$ in Assumption~\ref{ass:linearforcecoeff} are satisfied. In this case, we have
		\begin{align*}
		f_s^\epsilon(s)+s(f_s^\epsilon)'(s)= a_s(2s-0.5)
		\end{align*}
		for $s\in [0,0.5-\epsilon]$ with $a_s<0$ and 
		\begin{align*}
		f_s^\epsilon(s)+s(f_s^\epsilon)'(s)&=  - a_s \frac{(s-R_c)^3}{\epsilon^2}- 2 a_s \frac{\left(s-R_c \right)^2}{\epsilon} -3 a_s \frac{s(s-R_c)^2}{\epsilon^2} -4 a_s \frac{s \left(s-R_c \right)}{\epsilon}
		\end{align*}
		for $s \in [0.5-\epsilon, 0.5]$
		by Assumption \ref{ass:linearforcecoeff}. 
		Clearly, the leading order term of $\Re(\lambda_2(m))$ vanishes in the high-wave limit $m\to\infty$ and lower order terms in $\epsilon$ have to be considered. 
		An easy computation reveals that 
		\begin{align}\label{eq:highwaverc05}
		\lim_{m\to \infty}\Re(\lambda_2(m))=2\int_{0}^{0.5}
		f_s^\epsilon(s)+ s(f_s^\epsilon)'(s)
		=0
		\end{align}
		for any $\epsilon>0$ and as $\epsilon\to 0$.
		Further note that using \eqref{eq:highwaverc05}, $\Re(\lambda_2(m))$ reduces to
		\begin{align*}
		\Re(\lambda_2(m))=-2\int_{0}^{0.5}
		\bl  f_s^\epsilon(s)+ s(f_s^\epsilon)'(s)
		\br\cos\bl 2\pi ms\br\di s.
		\end{align*}
		We obtain
		\begin{align*}
		\int_0^{0.5-\epsilon} a_s(2s-0.5)\cos\bl 2\pi ms\br\di s=\begin{cases} -0.5 a_s \epsilon + a_s \epsilon^2+ \mathcal{O}(\epsilon^3), &  m\in 2\N, \\ -\frac{a_s}{\pi^2 m^2} +0.5 a_s \epsilon - a_s \epsilon^2+ \mathcal{O}(\epsilon^3), &  m\in 2\N+1,\end{cases}
		\end{align*}
		\begin{align*}
		-a_s\int_{0.5-\epsilon}^{0.5}\bl  \frac{(s-R_c)^3}{\epsilon^2}+ 2 \frac{\left(s-R_c \right)^2}{\epsilon}\br\cos\bl 2\pi ms\br\di s=\begin{cases} -\frac{5}{12} a_s \epsilon^2+ \mathcal{O}(\epsilon^3), &  m\in 2\N, \\  \frac{5}{12} a_s \epsilon^2+ \mathcal{O}(\epsilon^3), &  m\in 2\N+1,\end{cases}
		\end{align*}
		\begin{align*}
		-3a_s\int_{0.5-\epsilon}^{0.5} \frac{s \left(s-R_c \right)^2}{\epsilon^2}\cos\bl 2\pi ms\br\di s=\begin{cases} -0.5 a_s  \epsilon + \frac{3}{4} a_s\epsilon^2+ \mathcal{O}(\epsilon^3), &  m\in 2\N, \\0.5 a_s \epsilon -\frac{3}{4}a_s \epsilon^2+ \mathcal{O}(\epsilon^3), &  m\in 2\N+1,\end{cases}
		\end{align*}
		and
		\begin{align*}
		-4a_s\int_{0.5-\epsilon}^{0.5} \frac{s \left(s-R_c \right)}{\epsilon}\cos\bl 2\pi ms\br\di s=\begin{cases} a_s  \epsilon - \frac{4}{3} a_s\epsilon^2+ \mathcal{O}(\epsilon^3), &  m\in 2\N, \\ -a_s \epsilon + \frac{4}{3}a_s \epsilon^2+ \mathcal{O}(\epsilon^3), &  m\in 2\N+1,\end{cases}
		\end{align*}
		implying that
		\begin{align*}
		\int_{0}^{0.5}
		\bl  f_s^\epsilon(s)+ s(f_s^\epsilon)'(s)
		\br\cos\bl 2\pi ms\br\di s=\begin{cases}  \mathcal{O}(\epsilon^3), &  m\in 2\N, \\-\frac{a_s}{\pi^2 m^2} +  \mathcal{O}(\epsilon^3), &  m\in 2\N+1.\end{cases}
		\end{align*}
		Since the real part of the largest eigenvalue $\Re(\lambda_2(m))$ is zero in the high-wave number limit and it vanishes in the limit $\epsilon\to 0$ for any $m\in\N$, we cannot conclude stability/instability of the straight vertical line for $R_c=0.5$ and $\epsilon>0$ or $\epsilon\to 0$ in the continuum limit $N\to \infty$ or any $N\in \N$ sufficiently large. However, the numerical results in Section \ref{sec:numericalresultsstability} suggest instability for $\epsilon>0$ and in the limit $\epsilon\to 0$.
\end{remark}

Since we have the relations $f_l^\epsilon=f_A^\epsilon+f_R^\epsilon$ and $f_s^\epsilon=\chi f_A^\epsilon+f_R^\epsilon$ between the force coefficients $f_l^\epsilon,f_s^\epsilon$ in the general force formulation \eqref{eq:totalforcehom} and the total force \eqref{eq:totalforcenewhom} in the K\"ucken-Champod model  with repulsive and attractive force coefficients $f_R^\epsilon$ and $f_A^\epsilon$, respectively, we  have 
\begin{align*}
a_l=a_A+a_R,\quad a_s= \chi a_A+a_R,\quad b_l=b_A+b_R,\quad b_s=\chi b_A+b_R.
\end{align*}
Hence, Proposition \ref{prop:stabilityverticalline} leads to a similar statement for the forces in the K\"ucken-Champod model:
\begin{corollary}
		For $R_c\in(0,0.5)$ the single straight vertical line is an unstable steady state of \eqref{eq:particlemodelperiodic}  for any $N\in \N$ sufficiently large and for the continuum limit $N\to\infty$, where the forces are of the form \eqref{eq:totalforcehom} for any choice of parameters in the definition of the linear coefficient functions $f_R^\epsilon,f_A^\epsilon$ in \eqref{eq:linearforcecoeffkc} with $0<\epsilon\ll R_c$ or $\epsilon\to 0$. For $R_c=0.5$, the condition
		\begin{align*}
		a_A+a_R\leq -\frac{2\bl b_A+b_R\br}{R_c}\quad\text{and}\quad \chi a_A+a_R= -\frac{\chi b_A+b_R}{R_c}
		\end{align*}
		in addition to the assumptions on $a_A,a_R,b_A,b_R$ in Remark \ref{rem:linforcecoeff} is necessary for the stability of the single straight vertical line for force coefficients $f_R^\epsilon,f_A^\epsilon$ where $0<\epsilon\ll R_c$ or $\epsilon\to 0$. 
		This does not guarantee the stability/instability of the straight vertical line for force coefficients $f_R^\epsilon,f_A^\epsilon$ with $0<\epsilon\ll R_c$ or $\epsilon\to 0$ for $R_c=0.5$ and $N\in\N$ sufficiently large or in the continuum limit $N\to\infty$. 
\end{corollary}

\subsection{Algebraically decaying force coefficients}
Since the straight vertical line is unstable for $N\in\N$ sufficiently large and for $N\to\infty$ for the differentiable force coefficient $f_s^\epsilon$, defined in \eqref{eq:linearcoefffunctions} along $s$, which is linear  on $[0,R_c-\epsilon]$ for $R_c\in(0,0.5)$ and $\epsilon>0$, we consider faster decaying force coefficients along $s$ in the sequel. In this section we consider 
\begin{align*}
f_s(|d|)=\frac{c}{\bl 1+a|d|\br^b}
\end{align*}
for $a> 0$, $b>0$ and $c>0$. To obtain a differentiable force coefficient $f_s^\epsilon$ on $(0,\infty)$  we consider the modification in \eqref{eq:forcecutoff}, i.e.
\begin{align*}
f_s^\epsilon(|d|)=\begin{cases} \frac{c}{\bl 1+a|d|\br^b}, & |d|\in[0,R_c-\epsilon],\\
-\frac{abc}{(1+a(R_c-\epsilon))^{b+1}}\bl \frac{(|d|-R_c)^3}{\epsilon^2}+  \frac{\left(|d|-R_c \right)^2}{\epsilon}\br &\\\quad+  \frac{c}{\bl 1+a(R_c-\epsilon)\br^b} \bl 2\frac{(|d|-R_c)^3}{\epsilon^3}+  3\frac{\left(|d|-R_c \right)^2}{\epsilon^2} \br,  & |d| \in (R_c-\epsilon,R_c), 
\\ 0, & |d| \geq R_c
\end{cases}
\end{align*}
where $R_c \in (0,0.5]$.
Note that for this algebraically decaying force coefficient $f_s^\epsilon$, the necessary condition
${f_s^\epsilon(R_c)=0}$ in \eqref{eq:condvertlinehighwave} for high-wave number stability of a straight vertical line is satisfied.
To guarantee that the term $a|d|$ for $|d|\in[0,R_c]$ dominates the denominator and to avoid too large jumps we require $a\gg1$ additionally. The assumption $a\gg 1$ also guarantees that $f_s^\epsilon(R_c-\epsilon)\ll 1$. In this case, differences between the adaptation $f_s^\epsilon$ and the algebraically decaying force coefficient $f_s$, and their derivatives $(f_s^\epsilon)'$ and $f_s'$, are small. Without loss of generality we can assume that $c=1$ since this positive multiplicative constant leads to a rescaled stability condition but is not relevant for change of sign of the eigenvalues. Hence, we consider the algebraically decaying force coefficient
\begin{align}\label{eq:algebraicdecaycoeff}
f_s^\epsilon(|d|)=\begin{cases} \frac{1}{\bl 1+a|d|\br^b}, & |d|\in[0,R_c-\epsilon],\\
-\frac{ab}{(1+a(R_c-\epsilon))^{b+1}}\bl \frac{(|d|-R_c)^3}{\epsilon^2}+  \frac{\left(|d|-R_c \right)^2}{\epsilon}\br &\\\quad+  \frac{1}{\bl 1+a(R_c-\epsilon)\br^b} \bl 2\frac{(|d|-R_c)^3}{\epsilon^3}+  3\frac{\left(|d|-R_c \right)^2}{\epsilon^2} \br,  & |d| \in (R_c-\epsilon,R_c), 
\\ 0, & |d| \geq R_c
\end{cases}
\end{align}
in the sequel.

\begin{proposition}\label{prop:algebraicdecay}
	For the single straight vertical line to be a stable steady state of \eqref{eq:particlemodelperiodic} with forces of the form \eqref{eq:totalforce} for any $n\in\N$ sufficiently large and for the continuum limit $N\to \infty$ with algebraically decaying force coefficients $f_s^\epsilon$ of the form \eqref{eq:algebraicdecaycoeff} it is necessary that
	\begin{align*}
	b>1 \quad \text{and}\quad \frac{2}{a(b-1)}<R_c.
	\end{align*}
\end{proposition}
\begin{proof}
	Because of the definition of the eigenvalues \eqref{eq:eigenvaluesreallinegeneral} we consider
	\begin{align}\label{eq:algebraicdecayintegralhelp}
	\begin{split}
	&\int_{0}^{R_c}\bl  f_s^\epsilon(s)+ s(f_s^\epsilon)'(s) \br\bl 1-\cos\bl -2\pi ms\br\br\di s\\
	&= \int_{0}^{R_c-\epsilon}  \frac{1}{\bl 1+as\br^{b+1}}\bl  1+as(1-b)  \br\bl 1-\cos\bl -2\pi ms\br\br\di s 
	+\mathcal{O}(\epsilon).
	\end{split}
	\end{align}
	The linear function $ 1+as(1-b)$ is positive for $s\in(0,s_0)$ and negative for $s\in(s_0,R_c-\epsilon)$ for all $\epsilon>0$ where 
	\begin{align*}
	s_0=\frac{1}{a(b-1)}\in(0,R_c)
	\end{align*}
	implying $b>1$. 
	Note that the integral on the right-hand side of \eqref{eq:algebraicdecayintegralhelp} can be rewritten as
	\begin{align*}
	&\int_{0}^{R_c-\epsilon} g(s)\di s=\int_{0}^{s_0} g(s)\di s+\int_{s_0}^{R_c-\epsilon} g(s)\di s
	\end{align*}
	for any $\epsilon>0$
	where 
	\begin{align*}
	g(s)=\frac{1}{\bl 1+as\br^{b+1}}\bl  1+as(1-b)\br\bl 1-\cos\bl -2\pi ms\br\br
	\end{align*}
	is nonnegative on $[0,s_0]$ and not positive on $[s_0,R_c-\epsilon]$ for any $\epsilon>0$ by the definition of $s_0$ and the fact that $1-\cos\bl -2\pi ms\br\in(0,2)$. Setting 
	\begin{align*}
	h(s)=\frac{1}{\bl 1+as\br^{b+1}}
	\end{align*}
	note that  $h(R_c)<h(s_0)<h(0)=1$.
	A lower bound of the integral can be obtained by estimating $h(s)$  on $[0,s_0]$ by $h(s_0)$ due to the nonnegativity of the integrand and since the integrand changes sign at $s_0$ the factor $h(s)$   can be replaced by its maximum on $[s_0,R_c-\epsilon]$ for $\epsilon>0$, i.e.\ by $h(s_0)$, for obtaining a lower bound of the integral. Hence, a lower bound of the integral in \eqref{eq:algebraicdecayintegralhelp} is given by
	\begin{align*}
	&\frac{1}{\bl 1+as_0\br^{b+1}}\int_{0}^{R_c-\epsilon} \bl  1+as(1-b)\br\bl 1-\cos\bl -2\pi ms\br\br\di s
	+\mathcal{O}(\epsilon)\\
	&=\frac{1}{\bl 1+as_0\br^{b+1}}\bl \frac{2\pi m\bl \pi m (R_c-\epsilon)\bl p (R_c-\epsilon)+2q\br-\bl p (R_c-\epsilon)+q\br\sin\bl 2\pi m (R_c-\epsilon)\br\br }{4\pi^2 m^2}\right.\\&\quad \left.+\frac{p-p\cos\bl 2\pi m(R_c-\epsilon)\br}{4 \pi^2m^2}\br
	+\mathcal{O}(\epsilon)
	\end{align*}
	with $p=a(1-b)$ and $q=1$ where the explicit computation is analogous to the discussion of the linear force coefficients in \eqref{eq:straightlinestabilitylinear}. For large values of $m$ the first term of the above right-hand side dominates and we require 
	\begin{align*}
	pR_c+2q=a(1-b)R_c+2<0
	\end{align*}
	for all $\epsilon>0$. 
	This concludes the proof.
\end{proof}

In the sequel, we can restrict ourselves to algebraically decaying force coefficients \eqref{eq:algebraicdecaycoeff} with $a>0,~b>1$ due to  Proposition \ref{prop:algebraicdecay}. We show that the straight vertical line \eqref{eq:straightlineansatz} is an unstable steady state for any $N\in\N$ sufficiently large and in the continuum limit $N\to \infty$ in this case.  Due to the definition of the eigenvalues in \eqref{eq:eigenvaluesreallinegeneral} in Theorem \ref{th:stabilitystraightline} it is sufficient to show that there exists $m\in\N$ such that
\begin{align*}
0<\int_{0}^{R_c-\epsilon}\bl  f_s^\epsilon(s)+ s(f_s^\epsilon)'(s) \br\bl 1-\cos\bl -2\pi ms\br\br\di s
\end{align*}
for all $0<\epsilon\ll R_c$.
This is equivalent to showing that there exists $m\in \N$ such that 
\begin{align}\label{eq:algebraicdecaycomputation}
\begin{split}
0&<\lim_{\epsilon\to 0}\int_{0}^{R_c-\epsilon} \frac{1}{\bl 1+as\br^{b+1}}\bl  1+as(1-b)\br\bl 1-\cos\bl -2\pi ms\br\br\di s.
\end{split}
\end{align}
is satisfied.

\begin{lemma}
	For any $a>0,~b>1$ and $R_c\in(0,0.5]$ there exists $m\in \N$ such that 
	\eqref{eq:algebraicdecaycomputation} is satisfied.
\end{lemma}
\begin{proof}
	We denote the incomplete Gamma function by 
	\begin{align*}
	\Gamma(y,z)=\int_z^{\infty} s^{y-1}\exp(-s)\di s
	\end{align*}
	for $y\in\R$ and $z\in\C$.
	Then the right-hand side of \eqref{eq:algebraicdecaycomputation} can be written  as 
	\begin{align*}
	&-\frac{1}{4a\pi m\bl 1+aR_c\br^b} \bl \vphantom{Re \left[
		\sin\bl \frac{2\pi m}{a}\br   m^{b+1} \bl c_1\Gamma\bl -b,\frac{2i\pi m}{a}\br\br\right.} 2a\sin(-2\pi mR_c)\right.
	\\&\quad+\Re \left[
	\sin\bl \frac{2\pi m}{a}\br   m^{b+1} \bl c_1\Gamma\bl -b,\frac{2i\pi m}{a}\br+c_2\Gamma\bl -b,\frac{2i\pi (1+aR_c)m}{a}\br\br\right.
	\\&\quad+\sin\bl \frac{2\pi m}{a}\br   m^{b} \bl c_3\Gamma\bl 1-b,\frac{2i\pi m}{a}\br+c_4\Gamma\bl 1-b,\frac{2i\pi (1+aR_c)m}{a}\br\br
	\\&\quad+\cos\bl \frac{2\pi m}{a}\br   m^{b+1} \bl c_5\Gamma\bl -b,\frac{2i\pi m}{a}\br+c_6\Gamma\bl -b,\frac{2i\pi (1+aR_c)m}{a}\br\br
	\\&\quad+\left.\left.\cos\bl \frac{2\pi m}{a}\br   m^{b} \bl c_7\Gamma\bl 1-b,\frac{2i\pi m}{a}\br+c_8\Gamma\bl 1-b,\frac{2i\pi (1+aR_c)m}{a}\br\br\right]\right)
	\end{align*}
	for constants $c_i\in\C,i=1,\ldots,8$, depending on $a,b$ and $R_c$, but independent of $m$ where not all constants $c_i$ are equal to zero. Note that
	all incomplete Gamma functions above are of the form $\Gamma(-y,iz)$ for $y,z\in\R$ with $y,z>0$.  Integration by parts yields
		\begin{align*}
		\Gamma(-y,iz)=(iz)^{-y-1}\exp(-iz)+(-y-1)\Gamma(-y-1,iz)
		\end{align*}	
		where 
		\begin{align*}
		|\Gamma(-y-1,iz)|=\left|\int_{iz}^\infty s^{-y-2} \exp(-s)\di s \right|\leq |(iz)^{-y-2} \exp(-iz)|. 
		\end{align*}	
		In particular, we have
		\begin{align*}
		\Gamma(-y,iz)=(iz)^{-y-1} \exp(-iz)\left(1+\mathcal{O}((iz)^{-1})\right),
		\end{align*}
		implying
		\begin{align*}
		\Re (\Gamma(-y,iz))=\tilde{c} z^{-y-1} \left(1+\mathcal{O}(z^{-1})\right)
		\end{align*}	
		where $\tilde{c}=\Re(i^{-y-1}\exp(-iz))\in\R$. This leads to the approximation 
		\begin{align*}
		&	\Re \left[
		\sin\bl \frac{2\pi m}{a}\br   m^{b+1} \bl c_1\Gamma\bl -b,\frac{2i\pi m}{a}\br+c_2\Gamma\bl -b,\frac{2i\pi (1+aR_c)m}{a}\br\br\right.
		\\&\quad \left.+\sin\bl \frac{2\pi m}{a}\br   m^{b} \bl c_3\Gamma\bl 1-b,\frac{2i\pi m}{a}\br+c_4\Gamma\bl 1-b,\frac{2i\pi (1+aR_c)m}{a}\br\br\right]\\&=\tilde{c}_1\sin\bl \frac{2\pi m}{a}\br  \left(1+\mathcal{O}(m^{-1})\right)
		\end{align*}		
		for some constant $\tilde{c}_1\in\R$. The other terms of the right-hand side of \eqref{eq:algebraicdecaycomputation} can be rewritten in a similar way, resulting in
		\begin{align*}
		&-\frac{1}{4a\pi m\bl 1+aR_c\br^b} \bl  -2a\sin(2\pi mR_c)+\tilde{c}_1
		\sin\bl \frac{2\pi m}{a}\br +\tilde{c}_2\cos\bl \frac{2\pi m}{a}\br +\mathcal{O}(m^{-1})  \br
		\end{align*}
	for constants $\tilde{c}_1,\tilde{c}_2\in\R$, independent of $m$. Note that there exist infinitely many $m\in\N$ such that $\tilde{c}_1
		\sin\bl \frac{2\pi m}{a}\br +\tilde{c}_2\cos\bl \frac{2\pi m}{a}\br >0$ and such that $\tilde{c}_1
		\sin\bl \frac{2\pi m}{a}\br +\tilde{c}_2\cos\bl \frac{2\pi m}{a}\br <0$. If $R_c=\frac{1}{a}$, the second factor consists of the sum of a sine and a cosine function of the same period length and hence for $R_c \in (0,0.5]$ given, there exists $m\in\N$ such that the second factor is negative and the leading order term of \eqref{eq:algebraicdecaycomputation} is positive. If the first term in the second factor is of different period length as the second and third summand, there also exists $m\in\N$ such that the second factor is negative.  
	In particular, this implies that there exists an $m\in\N$ such that  \eqref{eq:algebraicdecaycomputation} is satisfied.
\end{proof}

\begin{corollary}\label{th:stabilityalgebraic}
	For any cutoff radius $R_c\in(0,0.5]$ the single straight vertical line is an unstable steady state of \eqref{eq:particlemodelperiodic} for any $N\in\N$ sufficiently large and for the continuum limit $N\to \infty$ with forces of the form \eqref{eq:totalforce} with algebraically decaying force coefficients $f_s^\epsilon$ of the form \eqref{eq:algebraicdecaycoeff} with $b>0$ and for any $\epsilon>0$ or in the limit $\epsilon\to 0$.
\end{corollary}

\subsection{Exponential force coefficients}\label{sec:expforcecoeff}
In this section we consider exponentially decaying force coefficient  along $s$ and short-range repulsive, long-range attractive forces along $l$ such that the necessary condition \eqref{eq:condvertlinehighwave} for high-wave number stability  is satisfied.

To express the  force coefficient along $l$ in terms of exponentially decaying functions we consider
\begin{align}\label{eq:expforcel}
\begin{split}
f_l^\epsilon(|d|)=\begin{cases} c_{l_1}\exp(-e_{l_1} |d|)+c_{l_2}\exp(-e_{l_2} |d|), & |d|\in[0,R_c-\epsilon],\\
\sum_{j=1}^2(-\epsilon c_{l_j}e_{l_j}+2c_{l_j}) \exp(-e_{l_j}(R_c-\epsilon))\frac{(|d|-R_c)^3}{\epsilon^3}&\\\quad+\sum_{j=1}^2(-\epsilon c_{l_j}e_{l_j}+3c_{l_j}) \exp(-e_{l_j}(R_c-\epsilon)) \frac{\left(|d|-R_c \right)^2}{\epsilon^2},  & |d| \in (R_c-\epsilon,R_c), \\ 0, & |d| \geq R_c,
\end{cases}
\end{split}
\end{align}
for parameters $c_{l_1}, c_{l_2}, e_{l_1}$ and $e_{l_2}$ with $e_{l_1}>0$ and $e_{l_2}>0$.
Note that exponentially decaying functions are either purely repulsive or purely attractive, depending on the sign of the multiplicative parameter. Since we require $f_l^\epsilon$ to be short-range repulsive, long-range attractive we  consider the sum of two exponentially decaying functions here. Without loss of generality we assume that the first summand in \eqref{eq:expforcel} is repulsive and the second one is attractive, i.e.\ $c_{l_1}>0>c_{l_2}$. To guarantee that $f_l^\epsilon$ is short-range repulsive we require $c_{l_1}>|c_{l_2}|$. For long-range attractive forces we require that the second term decays slower, i.e.\ $e_{l_1}>e_{l_2}$.
These assumptions lead to the parameter choice:
\begin{align}\label{eq:expforcelconstant}
c_{l_1}>0>c_{l_2}, \quad c_{l_1}>|c_{l_2}|\quad \text{and}  \quad e_{l_1}>e_{l_2}>0.
\end{align}
Note that we have 
\begin{align*}
\int_0^{R_c} f_l^\epsilon(s)\bl 1-\cos\bl -2\pi ms\br\br\di s=	\int_0^{R_c-\epsilon} f_l^\epsilon(s)\bl 1-\cos\bl -2\pi ms\br\br\di s+\mathcal{O}(\epsilon)
\end{align*}
due to the boundedness of $f_l^\epsilon$ on $[R_c-\epsilon,R_c]$ and hence it is sufficient to consider the integral on $[0,R_c-\epsilon]$ for $\epsilon>0$ sufficiently small and in the limit $\epsilon\to 0$.
Further note that for constants $c,e_l\in\R$ we obtain:
\begin{align*}
\int_0^{R_c-\epsilon}c\exp(-e_l s)\di s=\bl 1-\exp(-e_l (R_c-\epsilon))\br\frac{c}{e_l}.
\end{align*}
Hence, we require 
\begin{align*}
\bl 1-\exp(-e_{l_1} (R_c-\epsilon))\br\frac{c_{l_1}}{e_{l_1}}+\bl 1-\exp(-e_{l_2} (R_c-\epsilon))\br\frac{c_{l_2}}{e_{l_2}}\leq 0
\end{align*}
for all $\epsilon>0$
as in the necessary condition for high-wave number stability, implying 
\begin{align*}
\frac{c_{l_1}}{e_{l_1}}\leq \frac{\left|c_{l_2}\right|}{e_{l_2}}.
\end{align*}
Since
\begin{align*}
\int_0^{R_c-\epsilon}c_{l_1}\exp(-e_{l_1} s)\bl 1-\cos\bl -2\pi ms\br\br\di s &>0\,,\\
\int_0^{R_c-\epsilon}c_{l_2}\exp(-e_{l_2} s)\bl 1-\cos\bl -2\pi ms\br\br\di s &<0\,,
\end{align*}
for all $\epsilon>0$ and $m\in\N$ the parameters $c_{l_1},c_{l_2},e_{l_1},c_{l_2}$ in \eqref{eq:expforcelconstant} can clearly be chosen in such a way that 
\begin{align}\label{eq:condstabilityl}
\int_0^{R_c-\epsilon} f_l^\epsilon(s)\bl 1-\cos\bl -2\pi ms\br\br\di s\leq 0
\end{align}
is satisfied for all $m\in \N$ and $0<\epsilon\ll R_c$ where $f_l^\epsilon$ is defined in \eqref{eq:expforcel} with a cutoff radius $R_c\in(0,0.5]$. Note that the adaptation $f_l^\epsilon$ of $f_l$ is not negligible. However, due to the concentration of the particles along a straight vertical axis, this adaptation does not change the overall dynamics.

For the purely repulsive force coefficient $f_s^\epsilon$ we may consider a force coefficient of the form
\begin{align*}
\begin{split}
&f_s^\epsilon\colon \R_+\to \R,\\
&f_s^\epsilon(|d|)=\begin{cases} c\exp(-e_s |d|), & |d|\in[0,R_c-\epsilon],\\
-ce_s \exp(-e_s(R_c-\epsilon))\bl\frac{(|d|-R_c)^3}{\epsilon^2}+ \frac{\left(|d|-R_c \right)^2}{\epsilon}\br&\\\quad+ c \exp(-e_s(R_c-\epsilon))\bl 2\frac{(|d|-R_c)^3}{\epsilon^3}+  3\frac{\left(|d|-R_c \right)^2}{\epsilon^2} \br,   & |d| \in (R_c-\epsilon,R_c), \\ 0, & |d| \geq R_c,
\end{cases}
\end{split}
\end{align*}
by considering \eqref{eq:forcecutoff} for exponentially decaying force coefficients. Since 
\begin{align*}
\Re(\lambda_2(m))&=2\int_{0}^{R_c-\epsilon}
\bl  f_s^\epsilon(s)+ (f_s^\epsilon)'(s)s
\br\bl 1-\cos\bl -2\pi ms\br\br\di s+\mathcal{O}(\epsilon).
\end{align*} 
we require the non-positivity of $\Re(\lambda_2(m))$. Note that 
\begin{align*}
&\int_{0}^{R_c-\epsilon}
\bl  f_s^\epsilon(s)+ (f_s^\epsilon)'(s)s
\br\bl 1-\cos\bl -2\pi ms\br\br\di s 
\\&=(R_c-\epsilon)\exp(-e_s (R_c-\epsilon))-\frac{8e_s \pi^2 m^2  }{(4 \pi^2 m^2 + e_s^2)^2}	\\&\quad-\frac{\exp(-e_s (R_c-\epsilon))}{(4 \pi^2 m^2 + e_s^2)^2} \left[ e_s (4 \pi^2 m^2 (e_s (R_c-\epsilon) - 2) + e_s^3 (R_c-\epsilon)) \cos(2 \pi m (R_c-\epsilon)) \right.
\\&\qquad  \qquad
- 2 \pi m (4 \pi^2 m^2 (e_s (R_c-\epsilon) - 1) + e_s^2 (e_s (R_c-\epsilon) + 1)) \sin(2 \pi m (R_c-\epsilon))
\\&\qquad \qquad\left. - (R_c-\epsilon) (4 \pi^2 m^2 + e_s^2)^2\right],  
\end{align*}
implying that we have $\int_{0}^{R_c-\epsilon}
\bl  f_s^\epsilon(s)+ (f_s^\epsilon)'(s)s
\br\bl 1-\cos\bl -2\pi ms\br\br\di s >0$ for any $\epsilon>0$ and $m\in\N$ sufficiently large, i.e.\ high-wave stability cannot be achieved. However, note that $\exp(-e_s R_c)$ can be assumed to be very small for $e_s>0$ sufficiently large. This motivates to 
consider a force coefficient function of the form
\begin{align}\label{eq:expforcedef}
\begin{split}
&f_s^\epsilon\colon \R_+\to \R,\\
&f_s^\epsilon(|d|)=\begin{cases} c\exp(-e_s |d|)-c\exp(-e_s (R_c-\epsilon)), & |d|\in[0,R_c-\epsilon],\\
-ce_s \exp(-e_s(R_c-\epsilon))\bl\frac{(|d|-R_c)^3}{\epsilon^2}+ \frac{\left(|d|-R_c \right)^2}{\epsilon}\br,  & |d| \in (R_c-\epsilon,R_c), \\ 0, & |d| \geq R_c,
\end{cases}
\end{split}
\end{align}
with $c>0$ and $e_s>0$. Here, the first term in \eqref{eq:expforcedef} represents the exponential decay of the force coefficient. To approximate the high-wave number stability condition, we require $f_s(R_c-\epsilon)=0$ which can be guaranteed by subtracting the constant  $\exp(-e_s(R_c-\epsilon))$. Note  that we can choose $e_s\gg 1$ such that $\exp(-e_s (R_c-\epsilon))$ is a small positive number. Subtracting the constant $c\exp(-e_s (R_c-\epsilon))$ as in \eqref{eq:expforcedef} leads to $f_s^\epsilon(R_c-\epsilon)=0$. This additional constant only changes  the force coefficient $f_s^\epsilon$ slightly and does not change its derivative $(f_s^\epsilon)'$ on $[0,R_c-\epsilon]$, i.e.\ $f_s'=(f_s^\epsilon)'$ on $[0,R_c-\epsilon]$. Note that the differences between $f_s^\epsilon$ and $f_s$, and $(f_s^\epsilon)'$ and $f_s'$ on $[R_c-\epsilon,R_c]$ are negligible provided $e_s>0$ is chosen sufficiently large such that $e_s \exp(-e_s(R_c-\epsilon))\ll 1$. 
Thus, we make the following assumption in the sequel:
\begin{assumption}\label{ass:exponentialdecaycoeff}
	We assume that the purely repulsive, exponentially decaying force coefficient $f_s$  along $s$ is given by \eqref{eq:expforcedef}, i.e.\
	\begin{align*}
	&f_s^\epsilon\colon \R_+\to \R,\\
	&f_s^\epsilon(|d|)=\begin{cases} c\exp(-e_s |d|)-c\exp(-e_s (R_c-\epsilon)), & |d|\in[0,R_c-\epsilon],\\
	-ce_s \exp(-e_s(R_c-\epsilon))\bl\frac{(|d|-R_c)^3}{\epsilon^2}+ \frac{\left(|d|-R_c \right)^2}{\epsilon}\br,  & |d| \in (R_c-\epsilon,R_c), \\ 0, & |d| \geq R_c,
	\end{cases}
	\end{align*}
	where $c>0$  and $e_s\gg 1$. 
	For the forces along $l$ we either consider linear or exponentially decaying force coefficients. For a linear force coefficient we consider \eqref{eq:linearcoeffcond}, i.e.
	\begin{align*}
	f_l^\epsilon(|d|)&:=\begin{cases}a_l|d|+b_l, & |d|\in[0,R_c-\epsilon],\\(2 b_l + 2 R_c a_l - a_l \epsilon)\frac{(|d|-R_c)^3}{\epsilon^3}+(3 b_l+ 3 R_c a_l - 2 a_l\epsilon) \frac{\left(|d|-R_c \right)^2}{\epsilon^2},  & |d| \in (R_c-\epsilon,R_c), \\ 0, & |d| \geq R_c, \end{cases}
	\end{align*}
	where we assume that the parameters $a_l,b_l$ satisfy the sign conditions $a_l<0,~b_l>0$ in Assumption~\ref{ass:linearforcecoeff} as well as the necessary stability condition along $l$ in \eqref{eq:linearcoeffcond}. For an exponentially decaying force coefficient $f_l^\epsilon$ we assume that $f_l^\epsilon$ is of the form \eqref{eq:expforcel}, i.e.
	\begin{align*}
	\begin{split}
	f_l^\epsilon(|d|)=\begin{cases} c_{l_1}\exp(-e_{l_1} |d|)+c_{l_2}\exp(-e_{l_2} |d|), & |d|\in[0,R_c-\epsilon],\\
	\sum_{j=1}^2(-\epsilon c_{l_j}e_{l_j}+2c_{l_j}) \exp(-e_{l_j}(R_c-\epsilon))\frac{(|d|-R_c)^3}{\epsilon^3}&\\\quad+\sum_{j=1}^2(-\epsilon c_{l_j}e_{l_j}+3c_{l_j}) \exp(-e_{l_j}(R_c-\epsilon)) \frac{\left(|d|-R_c \right)^2}{\epsilon^2},  & |d| \in (R_c-\epsilon,R_c), \\ 0, & |d| \geq R_c,
	\end{cases}
	\end{split}
	\end{align*}
	for parameters 
	\begin{align*}
	c_{l_1}>0>c_{l_2}, \quad c_{l_1}>|c_{l_2}|\quad \text{and}  \quad e_{l_1}>e_{l_2}>0
	\end{align*} 
	as in \eqref{eq:expforcel}--\eqref{eq:expforcelconstant}  such that  the necessary stability condition  \eqref{eq:condstabilityl} for a straight vertical line is satisfied for all $m\in\N$ and $0<\epsilon\ll R_c$. 
\end{assumption}

\begin{theorem}\label{th:expstability}
	For the cutoff radius $R_c=0.5$, the straight vertical line is stable for the particle model \eqref{eq:particlemodelperiodic} for any $N\in \N$ sufficiently large  with the exponentially decaying force coefficient $f_s^\epsilon$ in \eqref{eq:expforcedef} along $s$ and a linear or exponentially decaying force coefficient $f_l^\epsilon$ as in Assumption \ref{ass:exponentialdecaycoeff} along $l$  in the limit $\epsilon\to 0$.
	For $R_c\in(0,0.5)$ the straight vertical line is an unstable steady state to \eqref{eq:particlemodelperiodic} for any $N\in\N$ sufficiently large and for the continuum limit $N\to \infty$ for any exponential decay $e_s>0$ in the limit $\epsilon\to 0$. For any $0<\epsilon\ll R_c$, the straight vertical lime is an unstable steady state for any $R_c \in (0,0.5]$. 
\end{theorem}
\begin{proof}
	Due to the assumptions on $f_l^\epsilon$ in Assumption \ref{ass:exponentialdecaycoeff} the real part for the first eigenvalue in \eqref{eq:eigenvaluesreallinegeneral}, given by
	\begin{align*}
	\Re(\lambda_1(m))&=2\int_{0}^{R_c-\epsilon}  f_l^\epsilon(s)\bl 1-\cos\bl - 2\pi ms\br\br\di s+\mathcal{O}(\epsilon),
	\end{align*}
	is not positive for any $m\in\N$ and any $0<\epsilon\ll R_c$ sufficiently small. 
	The real part of the second eigenvalue  \eqref{eq:eigenvaluesreallinegeneral} is given by 
	\begin{align*}
	\Re(\lambda_2(m))&=2\int_{0}^{R_c-\epsilon}
	\bl  f_s^\epsilon(s)+ (f_s^\epsilon)'(s)s
	\br\bl 1-\cos\bl -2\pi ms\br\br\di s+\mathcal{O}(\epsilon).
	\end{align*} 
	For the non-positivity of $\Re(\lambda_2(m))$
	it is sufficient to  require
	\begin{align}\label{eq:expcondestimate}
	\int_0^{R_c-\epsilon} \bl f_s^\epsilon(s)+s(f_s^\epsilon)'(s)\br (1-\cos(2\pi m s))\di s\leq 0,
	\end{align}
	for any $\epsilon>0$ sufficiently small where the left-hand side is given by
	\begin{align}\label{eq:expcondcomputation}
	\begin{split}
	&c\int_0^{R_c-\epsilon}\bl \exp(-e_s s)(1-e_ss)-\exp(-e_s (R_c-\epsilon))\br(1-\cos(2\pi m s))\di s\\&=-\frac{c e_s \exp\bl-e_s (R_c-\epsilon)\br }{2 \pi  m \left(e_s^2+4 \pi ^2 m^2\right)^2} \left[2 \pi  m \left(e_s^3 (R_c-\epsilon)+4 \pi ^2 m^2 (e_s (R_c-\epsilon)-2)\right) \cos (2 \pi  m (R_c-\epsilon))\right.\\&\qquad\left. -\left(e_s^3+4 \pi ^2 e_s^2 m^2 (R_c-\epsilon)+12 \pi ^2 e_s m^2+16 \pi ^4 m^4 (R_c-\epsilon)\right) \sin (2 \pi  m (R_c-\epsilon))\right.\\&\qquad\left.+16 \pi ^3 m^3 \exp\bl e_s (R_c-\epsilon)\br\right].
	\end{split}
	\end{align}	
	For $R_c=0.5$ we have $\lim_{\epsilon\to 0}\sin (2 \pi  m (R_c-\epsilon))=0$ and the right-hand side of \eqref{eq:expcondcomputation} simplifies to $g_\epsilon(m)h_\epsilon(m)$
	where 
	\begin{align*}
	g_\epsilon(m)&=-\frac{c e_s \exp\bl-e_s (R_c-\epsilon)\br }{ \left(e_s^2+4 \pi ^2 m^2\right)^2},\\
	h_\epsilon(m)&=  \left(e_s^3 (R_c-\epsilon)+4 \pi ^2 m^2 (e_s (R_c-\epsilon)-2)\right) \cos (2 \pi  m (R_c-\epsilon))+8 \pi ^2 m^2 \exp\bl e_s (R_c-\epsilon)\br.
	\end{align*}
	\\	
	For determining the limit $m\to\infty$ of $g_\epsilon(m)h_\epsilon(m)$ note that the leading order term of $g_\epsilon$ is $m^{-4}$ while the highest order term of $h_\epsilon$ is $m^2$, implying that the product $g_\epsilon(m)h_\epsilon(m)$, i.e.\  the right-hand side of \eqref{eq:expcondcomputation}, goes to zero as $m\to \infty$. 
	\\
	Note that for $R_c=0.5$ we have 
	\begin{align*}
	\lim_{\epsilon\to 0}\cos(2\pi m(R_c-\epsilon))=\begin{cases}
	1, & \text{m even},\\
	-1, & \text{m odd}.
	\end{cases}
	\end{align*}
	Let us consider $e_s>0$ with $e_s\leq 4$ first, i.e.\ $\lim_{\epsilon\to 0}e_s(R_c-\epsilon)\leq 2$. Then,
	\begin{align*}
	\lim_{\epsilon\to 0}h_\epsilon(m)= \begin{cases} e_s^3 R_c+4 \pi ^2 m^2 \bl e_s R_c-2+2\exp\bl e_s R_c\br\br, & \text{m even},\\
	-e_s^3 R_c+4 \pi ^2 m^2 \bl -e_s R_c+2+2\exp\bl e_s R_c\br\br, & \text{m odd}.\end{cases}
	\end{align*}
	Note that $\lim_{\epsilon\to 0}g_\epsilon(m)<0$ for all $m\in \N$ and $\lim_{\epsilon\to 0} h_\epsilon(m)>0$ for all even $m$ since $2\exp(e_sR_c)>2$. For $m$ odd, note that the term in brackets is positive and  a lower bound of $\lim_{\epsilon\to 0}h_\epsilon$ is given by
	\begin{align*}
	-16 e_sR_c+4 \pi ^2  \bl -e_s R_c+2+2\exp\bl e_s R_c\br\br\geq 8 \pi ^2  \bl -e_s R_c+1+\exp\bl e_s R_c\br\br,
	\end{align*}
	which is clearly positive. Hence, $\lim_{\epsilon\to 0}h_\epsilon(m)$ is positive for all $m\in\N$ and thus, we obtain $\lim_{\epsilon\to 0}g_\epsilon(m)h_\epsilon(m)<0$, provided $e_s\leq 4$ and $R_c=0.5$. This implies that \eqref{eq:expcondestimate} is satisfied for all $m\in\N$ in this case.
	\\	
	Let us now consider $\lim_{\epsilon\to 0}e_s(R_c-\epsilon)>2$ with $R_c=0.5$. Note that a lower bound of $\lim_{\epsilon\to 0}h_\epsilon$ is obtained from $\lim_{\epsilon\to 0}\cos(2\pi m (R_c-\epsilon))\geq-1$, leading to the upper bound 
	\begin{align*}
	\lim_{\epsilon\to 0} g_\epsilon(m) \left[ -\left(e_s^3 R_c+4 \pi ^2 m^2 (e_s R_c-2)\right) +8 \pi ^2 m^2 \exp\bl e_s R_c\br\right]
	\end{align*}
	of $\lim_{\epsilon\to 0}g_\epsilon(m)h_\epsilon(m)$ since $g_\epsilon(m)<0$ for all $\epsilon>0$. This upper bound can be rewritten as 
	\begin{align*}
	\lim_{\epsilon\to 0} g_\epsilon(m) \left[ -e_s^3 R_c+4 \pi ^2 m^2\bl- e_s R_c+2 +2 \exp\bl e_s R_c\br\br\right].
	\end{align*}
	Note that $- e_s R_c+2 +2 \exp\bl e_s R_c\br>0$. Besides, 
	\begin{align*}
	\frac{e_s^3 R_c}{4 \pi ^2 \bl- e_s R_c+2 +2 \exp\bl e_s R_c\br\br}<1
	\end{align*}
	is satisfied for all $e_s>4$, implying
	\begin{align*}
	-e_s^3 R_c+4 \pi ^2 m^2\bl- e_s R_c+2 +2 \exp\bl e_s R_c\br\br>0
	\end{align*}
	for all $m\in\N$. Hence,  the right-hand side of \eqref{eq:expcondcomputation}, i.e.\ $g_\epsilon(m)h_\epsilon(m)$, is negative for all $m\in\N$ in the limit $\epsilon\to 0$. In particular, this shows that condition \eqref{eq:expcondestimate} is satisfied for all $m\in \N$ for $R_c=0.5$.
	\\	
	For $R_c\in(0,0.5]$ and $\epsilon>0$ we have $\sin(2\pi m (R_c-\epsilon))> 0$ for countably many $m\in\N$. In particular, there exists $\delta>0$ and a countably infinite set $\mathcal{N}\subset \N$ such that $\sin(2\pi m (R_c-\epsilon))> \delta$ for all $m\in\mathcal{N}$.   Hence the second term in \eqref{eq:expcondcomputation} is negative with upper bound
	\begin{align*}
	-\left(e_s^3+4 \pi ^2 e_s^2 m^2 (R_c-\epsilon)+12 \pi ^2 e_s m^2+16 \pi ^4 m^4 (R_c-\epsilon)\right) \delta<0
	\end{align*}
	for all $m\in\mathcal{N}$. 
	This implies that the right-hand side of \eqref{eq:expcondcomputation}, i.e. $g_\epsilon(m),h_\epsilon(m)$, can be estimated from below by
	\begin{align*}
	&\frac{g_\epsilon(m)}{2\pi m}\left[2 \pi  m \max\{e_s^3 (R_c-\epsilon)+4 \pi ^2 m^2 (e_s (R_c-\epsilon)-2),-e_s^3 (R_c-\epsilon)-4 \pi ^2 m^2 (e_s (R_c-\epsilon)-2)\}\right.\\&\qquad \left.+16 \pi ^3 m^3 \exp\bl e_s (R_c-\epsilon)\br-\left(e_s^3+4 \pi ^2 e_s^2 m^2 (R_c-\epsilon)+12 \pi ^2 e_s m^2+16 \pi ^4 m^4 (R_c-\epsilon)\right) \delta\right]
	\end{align*}
	for all $m\in\mathcal{N}$ and $0<\epsilon\ll R_c$ since $g_\epsilon(m)<0$ for all $m\in \N$. Thus, there exists $m_0\in\mathcal{N}$ such that the term in square brackets is negative for all $m\in\mathcal{N}$ with $m\geq m_0$ and all $\epsilon>0$ sufficiently small since the highest order term of power $m^4$ in the square brackets dominates for $m$ large enough. In particular,  $g_\epsilon(m_0)<0$ for $\epsilon> 0$ implies that we have found a positive lower bound of the right-hand side in \eqref{eq:expcondcomputation} and one can easily show that this positive lower bound also holds in the limit $\epsilon\to 0$. Hence, stability cannot be achieved in the case $R_c\in(0,0.5]$ and any $\epsilon>0$, as well as for $R_c \in (0,0.5)$ and $\epsilon\to 0$, both for $N\in\N$ sufficiently large and in the continuum limit $N\to\infty$.
\end{proof}

\begin{remark}\label{rem:expstability}
	For $R_c\in (0,0.5)$ and $\epsilon\to 0$, no stability  can be shown analytically. However, note that 
	an upper bound of the integral
	\begin{align}\label{eq:expintegral}
	\int_0^{R_c-\epsilon} \bl f_s^\epsilon(s)+s(f_s^\epsilon)'(s)\br (1-\cos(2\pi m s))\di s
	\end{align}
	in the necessary stability condition \eqref{eq:expcondestimate} is given by 
	\begin{align}\label{eq:expcomputationupperbound}
	\begin{split}
	&\quad-\frac{c e_s \exp\bl-e_s (R_c-\epsilon)\br }{2 \pi  m \left(e_s^2+4 \pi ^2 m^2\right)^2} \left[-2 \pi  m \left(e_s^3 (R_c-\epsilon)+4 \pi ^2 m^2 (e_s (R_c-\epsilon)-2)\right) \right.\\&\qquad\left. -\left(e_s^3+4 \pi ^2 e_s^2 m^2 (R_c-\epsilon)+12 \pi ^2 e_s m^2+16 \pi ^4 m^4 (R_c-\epsilon)\right) +16 \pi ^3 m^3 \exp\bl e_s (R_c-\epsilon)\br\right]\\
	&=-\frac{c e_s \exp\bl-e_s (R_c-\epsilon)\br }{2 \pi  m \left(e_s^2+4 \pi ^2 m^2\right)^2} \left[-e_s^3-2 \pi   e_s^3 (R_c-\epsilon) m-\bl 4 \pi ^2 e_s^2  (R_c-\epsilon)+12 \pi ^2 e_s\br m^2 \right.\\&\qquad\left. +\bl -8 \pi ^3  (e_s (R_c-\epsilon)-2) +16 \pi ^3\exp\bl e_s (R_c-\epsilon)\br\br m^3 -16 \pi ^4 (R_c-\epsilon) m^4 \right]
	\end{split}
	\end{align}
	for any $0<\epsilon\ll R_c$
	due to \eqref{eq:expcondcomputation}. For $\exp(e_sR_c)\gg 1$ there exists $m_0\in\N$ of order $ \exp(e_sR_c)\gg 1$ such that the term $16 \pi ^3 m^3 \exp\bl e_s R_c\br$ is the dominating term in  the upper bound \eqref{eq:expcomputationupperbound} of the integral \eqref{eq:expintegral} for all $m\in\N$ with $m\leq m_0$. Hence negativity of the upper bound \eqref{eq:expcomputationupperbound} and thus of the integral \eqref{eq:expintegral} in the necessary stability condition  can be guaranteed for all $m\leq m_0$. For $m>m_0$, however,  the highest order term of power $m^4$ dominates the sum. Since $m_0\gg 1$, we have stability for $N\in\N$ sufficiently large and for the continuum limit $N\to \infty$ for almost all, but finitely many, Fourier modes for $e_s\gg 1$, $R_c\in (0,0.5)$ and any $\epsilon>0$ sufficiently small or in the limit $\epsilon\to 0$. 
	
	The integral \eqref{eq:expintegral} is explicitly evaluated  in \eqref{eq:expcondcomputation}. For large values of $m\in\N$  the highest order term in \eqref{eq:expcondcomputation} is associated with the summand $16\pi^4m^4(R_c-\epsilon)\sin(2\pi m(R_c-\epsilon))$ and can be written in the form
	\begin{align*}
	\frac{8\pi^3e_s \exp(-e_s(R_c-\epsilon))(R_c-\epsilon) m^3\sin(2\pi m(R_c-\epsilon))}{\left(e_s^2+4 \pi ^2 m^2\right)^2}.
	\end{align*}
	Here, the numerator increases as $m^3$ for large $m$ while the denominator is of order $m^4$, multiplied by a factor $\exp(-e_sR_c)\ll 1$, leading to decaying sinusoidal oscillations around zero as $m$ increases. Since this approximation is only valid for $m>m_0\gg 1$ the absolute value of the right-hand side in \eqref{eq:expcondcomputation} may be so small that it is numerically zero and one may see stable vertical line patterns for exponentially decaying force coefficients $f_s^\epsilon$ along $s$ for $R_c\in (0,0.5)$, $\epsilon>0$ or in the limit $\epsilon\to 0$, and  $N\in\N$ sufficiently large, see the numerical experiment in Figure \subref*{fig:expcoeffnumericssmallrc}. 
\end{remark}

\begin{corollary}\label{col:expsumforce}
	Let $c_1,c_2\in\R$ with $c_1>0$, $c_1>|c_2|$ be given. There exist parameters   ${e_2\geq e_1>0}$ such that the straight vertical line is stable for the particle model \eqref{eq:particlemodelperiodic} for $N\in\N$ sufficiently large for 
	the exponentially decaying force coefficient $f_s^\epsilon$ along $s$ given by 
	$f_s^\epsilon\colon \R_+\to \R$ with
	\begin{align}\label{eq:expsumcoeff}
	f_s^\epsilon(|d|)=\begin{cases} c_1\exp(-e_1 |d|)+c_2\exp(-e_2 |d|)-c, & |d|\in[0,R_c-\epsilon],\\
	(f_s^\epsilon)'(R_c-\epsilon)\bl\frac{(|d|-R_c)^3}{\epsilon^2}+  \frac{\left(|d|-R_c \right)^2}{\epsilon}\br,  & |d| \in (R_c-\epsilon,R_c), 
	\\ 0, & |d| \geq R_c,
	\end{cases}
	\end{align}
	with
	\begin{align*}
	c=c_1\exp(-e_1 (R_c-\epsilon))+c_2\exp(-e_2 (R_c-\epsilon))
	\end{align*}
	and a linear or an exponential  force coefficient $f_l^\epsilon$ along $l$ as in Assumption \ref{ass:exponentialdecaycoeff} for a cutoff radius $R_c=0.5$. For the continuum limit $N\to \infty$ stability/instability cannot be concluded.
\end{corollary}
\begin{proof}
	For the stability of the straight vertical line for $N\in\N$ sufficiently large we require that the force coefficient $f_s^\epsilon$ in \eqref{eq:expsumcoeff} is purely repulsive for any $\epsilon>0$ and hence at least one of the constants $c_1,c_2$ has to be positive. Since we can assume $c_1>0$ without loss of generality this implies that $c_1$ is a repulsive multiplicative factor, while the sign of $c_2$ is not given by the assumptions. Thus, we require that the first  term in the definition of  $f_s^\epsilon$ in \eqref{eq:expsumcoeff} decays slower than the second one, implying $0<e_1\leq e_2$. Hence, the conditions on the parameters are verified.
	
	As in the proof of Theorem \ref{th:expstability} we evaluate integrals of the form \eqref{eq:expcondcomputation} where the term with factor $\sin(2\pi m R_c)$ vanishes for our choice of $R_c=0.5$.  If $c_2\geq 0$ one can choose $e_1,e_2$ sufficiently large such  that the term $16\pi^3m^3\exp(e_{k}R_c), k=1,2,$ in the square brackets in \eqref{eq:expcondcomputation} dominates as in the proof of Theorem \ref{th:expstability}, leading to the stability of the vertical straight line for $N\in\N$ sufficiently large. For $c_2<0$ one can choose $e_1,e_2$ sufficiently large such  that the term $16\pi^3m^3\exp(e_{k}R_c), k=1,2,$ dominate the square brackets. However, since $c_1>0>c_2$ we require in addition that the term  with multiplicative factor $c_1$ dominates over the term with multiplicative factor $c_2$, leading to the condition
	\begin{align*}
	&\quad-\frac{c_1 e_1 \exp\bl-e_1 R_c\br }{2 \pi  m \left(e_1^2+4 \pi ^2 m^2\right)^2} \left[16 \pi ^3 m^3 \exp\bl e_1 R_c\br\right]+\frac{|c_2| e_2 \exp\bl-e_2 R_c\br }{2 \pi  m \left(e_2^2+4 \pi ^2 m^2\right)^2} \left[16 \pi ^3 m^3 \exp\bl e_2 R_c\br\right]\\
	&=-\frac{16 \pi ^3 m^3c_1 e_1 }{2 \pi  m \left(e_1^2+4 \pi ^2 m^2\right)^2} +\frac{16 \pi ^3 m^3|c_2| e_2  }{2 \pi  m \left(e_2^2+4 \pi ^2 m^2\right)^2}<0
	\end{align*}
	in the limit $\epsilon\to 0$
	which is equivalent to
	\begin{align*}
	-c_1 e_1  \left(e_2^2+4 \pi ^2 m^2\right)^2 +|c_2| e_2 \left(e_1^2+4 \pi ^2 m^2\right)^2 <0.
	\end{align*}
	Since $c_1>|c_2|$ and $e_2\geq e_1>0$  by assumption this condition is satisfied for $e_2>e_1$ sufficiently large. Hence, stability of the straight vertical line can be achieved for $N\in\N$ sufficiently large.
\end{proof}
The force  coefficient $f_s^\epsilon$ of the form \eqref{eq:expsumcoeff} along $s$ is motivated by the force coefficients in the K\"ucken-Champod model. Here, $f_s^\epsilon=\chi f_A^\epsilon+f_R^\epsilon$ for $\chi\in[0,1]$ where, motivated by this section, $f_R^\epsilon,f_A^\epsilon$ are defined as
\begin{align}\label{eq:repulsionforcemodelconstant}
f_R^\epsilon(|d|)=\begin{cases} f_R(|d|)-f_R(R_c-\epsilon), & |d|\in[0,R_c-\epsilon],\\
f_R'(R_c-\epsilon)\bl \frac{(|d|-R_c)^3}{\epsilon^2}+  \frac{\left(|d|-R_c \right)^2}{\epsilon}\br, 
& |d| \in (R_c-\epsilon,R_c), 
\\ 0, & |d| \geq R_c,
\end{cases}
\end{align}
and 
\begin{align}\label{eq:attractionforcemodelconstant}
f_A^\epsilon(|d|)=\begin{cases} f_A(|d|)-f_A(R_c-\epsilon), & |d|\in[0,R_c-\epsilon],\\
f_A'(R_c-\epsilon)\bl \frac{(|d|-R_c)^3}{\epsilon^2}+  \frac{\left(|d|-R_c \right)^2}{\epsilon}\br, 
& |d| \in (R_c-\epsilon,R_c), 
\\ 0, & |d| \geq R_c.
\end{cases}
\end{align}
This corresponds to the sum of an attractive and a repulsive force coefficient as in \eqref{eq:expsumcoeff} for $c_1>0>c_2$ where the repulsive term, i.e.\ $c_1>|c_2|$, dominates. This motivates that we obtain stability of the straight vertical line for the force coefficients in the K\"ucken-Champod model for $N\in\N$ sufficiently large by considering  force coefficients of the form \eqref{eq:repulsionforcemodelconstant}, \eqref{eq:attractionforcemodelconstant}.

\subsection{Kücken-Champod model}

\begin{figure}[htbp]
	\centering
	\includegraphics[width=0.45\textwidth]{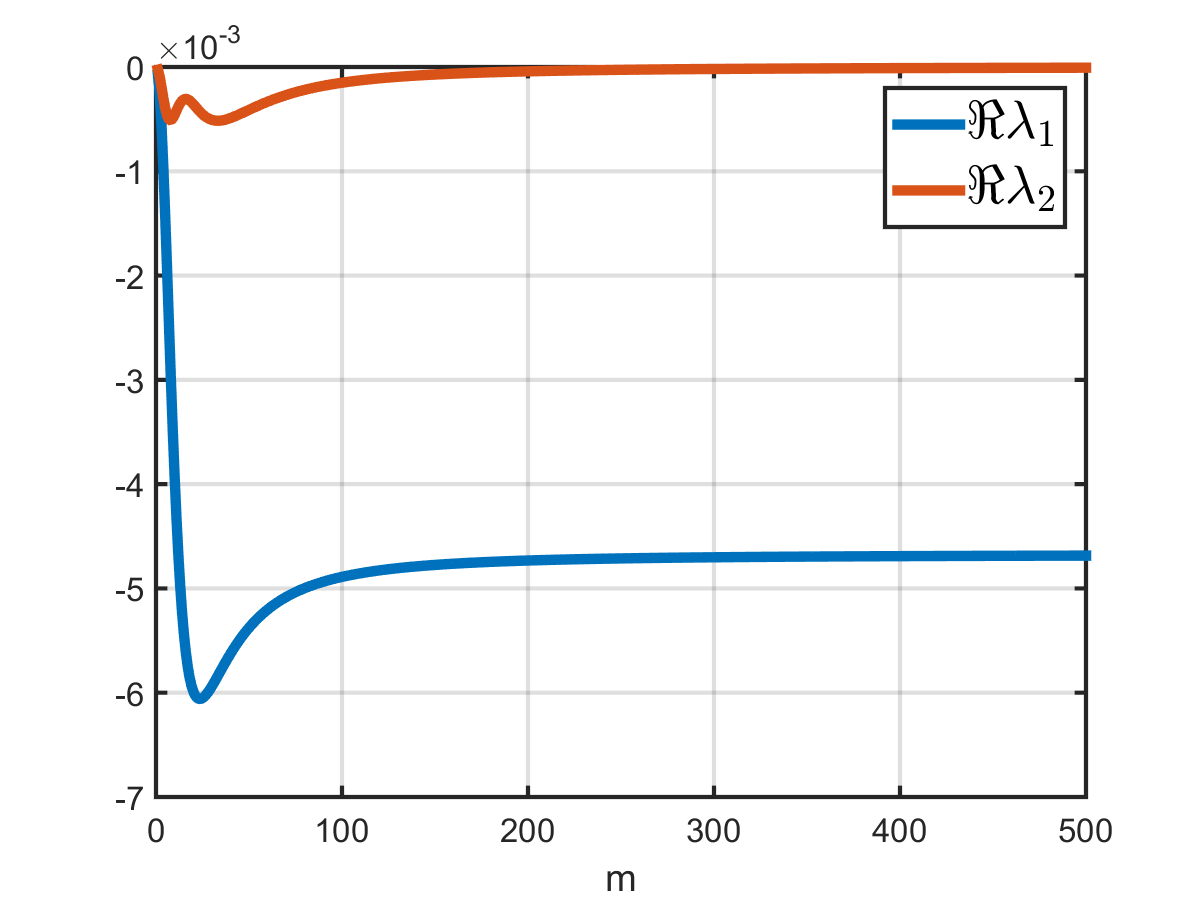}
	\caption{$\Re(\lambda_i)$ in \eqref{eq:eigenvaluesreallinegeneral} as a function of $m$ for the force coefficients $f_R^\epsilon$ in \eqref{eq:repulsionforcemodelconstant} and $f_A^\epsilon$ in \eqref{eq:attractionforcemodelconstant} of repulsion force \eqref{eq:repulsionforce} and attraction force \eqref{eq:attractionforce}, respectively,  for parameter values in  \eqref{eq:parametervaluesRepulsionAttraction} in the limit $\epsilon\to 0$ where $f_l^\epsilon=f_A^\epsilon+f_R^\epsilon$ and ${f_s^\epsilon=\chi f_A^\epsilon+f_R^\epsilon}$ \label{fig:kueckenchampodstability}.}
\end{figure}

For the specific forces in the K\"{u}cken-Champod model, given by the repulsive and attractive force coefficients $f_R^\epsilon$ and $f_A^\epsilon$ in \eqref{eq:repulsionforcemodelconstant} and \eqref{eq:attractionforcemodelconstant}, respectively, we require the non-positivity of the real parts of the eigenvalues $\lambda_k,k=1,2$, given by
\begin{align*}
\Re(\lambda_1(m))&=2\int_{0}^{R_c}  f_l^\epsilon(s)\bl 1-\cos\bl - 2\pi ms\br\br\di s,\\
\Re(\lambda_2(m))&=2\int_{0}^{R_c}
\bl  f_s^\epsilon(s)+ s(f_s^\epsilon)'(s)
\br\bl 1-\cos\bl -2\pi ms\br\br\di s
\end{align*}
in \eqref{eq:eigenvaluesreallinegeneral} where $f_l^\epsilon=f_A^\epsilon+f_R^\epsilon$ and $f_s^\epsilon=\chi f_A^\epsilon+f_R^\epsilon$. 
In Figure \ref{fig:kueckenchampodstability} we evaluate $\Re (\lambda_k)$ numerically for the force coefficients \eqref{eq:repulsionforcemodelconstant} and \eqref{eq:attractionforcemodelconstant} in the K\"{u}cken-Champod model for the parameters in \eqref{eq:parametervaluesRepulsionAttraction} and a cutoff radius $R_c=0.5$ in the limit $\epsilon\to 0$. Clearly, $\Re(\lambda_1)\leq 0$, while $\Re(\lambda_2)$ is negative for small modes $m$ but tends to zero for large modes $m$. 
Investigating the high-wave number stability for the forces in the K\"{u}cken-Champod model can  be done analytically.  For the  general necessary high-wave number condition \eqref{eq:condvertlinehighwave} for $\lambda_1$ we require $$\int_0^{R_c} f_l^\epsilon\di s\leq 0.$$ Note that 
\begin{align*}
&\lim_{\epsilon\to 0}\int_0^{R_c-\epsilon}\exp(-e_Rs)\bl \alpha s^2+\beta\br-\gamma\exp(-e_A s)s\di s
\\&=\frac{\alpha \left(\exp(-e_R R_c) (-e_R R_c(e_R R_c+2)-2)+2\right)}{(e_R)^3}+\frac{\beta-\beta \exp(-e_R R_c)}{e_R}\\&\quad-\frac{\gamma \left(1-\exp(-e_A R_c) (e_A R_c+1)\right)}{(e_A)^2}\\&\approx \frac{2\alpha}{(e_R)^3}+\frac{\beta}{e_R}-\frac{\gamma}{(e_A)^2}
\end{align*}
which is clearly negative for the choice of parameters in \eqref{eq:parametervaluesRepulsionAttraction}. For the high-wave stability we also consider the condition associated with $\lambda_2$, leading to the condition  $${\int_{0}^{R_c}   f_s^\epsilon(s)+ s(f_s^\epsilon)'(s) \di s \leq 0}.$$
We evaluate the integral
\begin{align*}
&\int_0^{R_c-\epsilon} \exp(-e_R s)\bl\alpha \bl 3s^2-e_R s^3\br+\beta(1-e_R s)\br-\chi\gamma\exp(-e_As) s(2-e_A s)\di s
\\&=(R_c-\epsilon) \left[(R_c-\epsilon) \left[\alpha (R_c-\epsilon) \exp(-e_R (R_c-\epsilon))-\chi \gamma \exp(-e_A (R_c-\epsilon))\right]\right.\\&\qquad+\left.\beta \exp(-e_R (R_c-\epsilon))\right]
\\&=(R_c-\epsilon)\bl f_R(R_c-\epsilon)+\chi f_A(R_c-\epsilon)\br
\end{align*}
for $f_R$ and $f_A$ defined in \eqref{eq:repulsionforcekc} and \eqref{eq:attractionforcekc}, respectively, implying that $$\int_0^{R_c-\epsilon} f_s^\epsilon(s) +s(f_s^\epsilon)'(s)\di s=0$$ for any $\epsilon>0$. In particular, the  straight vertical line is high-wave number stable for any $N\in\N$ sufficiently large and in the continuum limit $N\to \infty$ for the K\"ucken-Champod model with force coefficients $f_R^\epsilon$ and $f_A^\epsilon$ in \eqref{eq:repulsionforcemodelconstant} and \eqref{eq:attractionforcemodelconstant}, respectively, the parameters in \eqref{eq:parametervaluesRepulsionAttraction} and $\epsilon\to 0$. By definition of $f_s^\epsilon=\chi f_A^\epsilon+f_R^\epsilon$, we have $f_s^\epsilon(R_c)=0$, i.e.\ the high-wave number stability of the straight vertical line, compare Proposition \ref{prop:highwavestraightline}, is satisfied. Note that $$\lim_{\epsilon\to 0}{\chi f_A(R_c-\epsilon)+f_R(R_c-\epsilon)= 4.8144\cdot 10^{-21}}$$ for $R_c=0.5$, i.e.\  the force coefficient $\chi f_A+f_R$ has only slightly been modified to obtain $\chi f_A^\epsilon+f_R^\epsilon$ with $(\chi f_A^\epsilon+f_R^\epsilon)'\approx (\chi f_A^\epsilon+f_R^\epsilon)'$, provided $e_R \exp(-e_R R_c)\ll 1$ and $e_A \exp(-e_A R_c)\ll 1$.  


Note that it is not possible to analyse the stability of the straight vertical line for all modes $m\in\N$ for the forces $f_R^\epsilon$ and $f_A^\epsilon$ in \eqref{eq:repulsionforcemodelconstant} and \eqref{eq:attractionforcemodelconstant} in the K\"ucken-Champod model analytically for all possible parameter values due to the large number of parameters in the model. Besides, the force coefficients strongly depend on the choice of parameters. In Corollary~\ref{col:expsumforce}, however, we investigated the stability of the straight vertical line for $N\in\N$ sufficiently large where $f_s^\epsilon$, restricted to $[0,R_c-\epsilon]$ for some $\epsilon>0$, is the sum of the positive term $c_1\exp(-e_1|d|)$, the negative term $c_2\exp(-e_2|d|)$ and a constant to guarantee $f_s^\epsilon(R_c-\epsilon)=0$ where $c_1>|c_2|>0$. Besides, we required  $e_1<e_2$ for the positivity of the sum $c_1\exp(-e_1|d|)+c_2\exp(-e_2|d|)$ for $|d|\in [0,R_c-\epsilon]$ and showed stability of the straight vertical line for $N\in\N$ sufficiently large provided the parameters $e_1,e_2>0$ are chosen sufficiently large enough. In Figure \ref{fig:forces} the absolute value of the terms $\chi f_A$ and $f_R$, defined  in  \eqref{eq:repulsionforcekc}--\eqref{eq:attractionforcekc}, are plotted for the parameters in \eqref{eq:parametervaluesRepulsionAttraction}. As in Corollary \ref{col:expsumforce} the positive term always dominates and the terms $\chi f_A$ and $f_R$ have  fast exponential decays. This suggests that the straight vertical line is a stable steady state for the K\"ucken-Champod model for $N\in\N$ sufficiently large with the adopted force coefficient $f_s^\epsilon=\chi f_A^\epsilon+f_R^\epsilon$. Besides, the numerical evaluation of the real part of the eigenvalue $\lambda_2$ for $f_s^\epsilon$  for $\epsilon>0$, i.e.\ a differentiable force coefficient with the additional constant $-\bl \chi f_A(R_c-\epsilon)+f_R(R_c-\epsilon) \br$ for $|d|\in[0,R_c-\epsilon]$ leads to non-positivity of the real part of the eigenvalue $\lambda_2$.

\subsection{Summary}
In this section, we summarize the results from the previous subsections  on the stability of the  straight vertical line \eqref{eq:straightlineansatz} of the particle model \eqref{eq:particlemodelperiodic} with linear, algebraically decaying and exponentially decaying force coefficients for different values of the cutoff radius $R_c\in(0,0.5]$. This summary is shown in Table~\ref{table:summary}.

\begin{table}
	\caption{Stability/Instability of the straight vertical line \eqref{eq:straightlineansatz} for the particle model \eqref{eq:particlemodelperiodic} with  force coefficients $f_s$ along $s$ and different cutoff radii $R_c\in(0,0.5]$.}
	\label{table:summary}
	\begin{tabular}{p{5cm}p{4cm}p{4cm}}
		\hline\noalign{\smallskip}
		Force coefficient $f_s$ along $s$ & $R_c\in(0,0.5)$ & $R_c=0.5$ \\
		\noalign{\smallskip}\hline\noalign{\smallskip}
		Linear force coefficient \eqref{eq:linearcoefffunctions} & Instability for any ${N\in\N}$ sufficiently large and for $N\to \infty$ (see Theorem~\ref{prop:stabilityverticalline}) & Stability or instability  since stability conditions are satisfied with equality (see Corollary~\ref{prop:stabilityverticalline}) \\\noalign{\smallskip}\noalign{\smallskip}
		Algebraically decaying force coefficient \eqref{eq:algebraicdecaycoeff} & Instability for any ${N\in\N}$ sufficiently large and for $N\to \infty$ (see Corollary~\ref{th:stabilityalgebraic}) & Instability for any ${N\in \N}$ sufficiently large and for $N\to\infty$ (see Corollary~\ref{th:stabilityalgebraic}) \\\noalign{\smallskip}\noalign{\smallskip}
		Exponentially decaying force coefficient \eqref{eq:expforcedef} & Instability for any ${N\in\N}$ sufficiently large and for ${N\in\N}$ (see Theorem \ref{th:expstability}), but stability may be seen in numerical simulations (see Remark \ref{rem:expstability}) & Stability for any ${N\in\N}$ sufficiently large (see Theorem~\ref{th:expstability}) \\
		\noalign{\smallskip}\hline\noalign{\smallskip}
	\end{tabular}
\end{table}

\section{Numerical simulations}\label{sec:numerics}

\subsection{Numerical methods}
As in \cite{patternformationanisotropicmodel,anisotropicfingerprint} we consider the unit square with periodic boundary conditions as the domain for our numerical simulations if not stated otherwise. The particle system \eqref{eq:particlemodelperiodic} is solved by either the simple explicit Euler scheme
or higher order methods such as the Runge-Kutta-Dormand-Prince method, all resulting in very
similar simulation results. Note that the time step has to be adjusted depending on the value of the cutoff radius $R_c$. For efficient numerical simulation we consider cell lists as outlined in \cite{anisotropicfingerprint}. 

\subsection{Numerical results}\label{sec:numericalresultsstability}
\begin{figure}[htbp]
	\centering
	\subfloat[Initial data]{\includegraphics[width=0.32\textwidth]{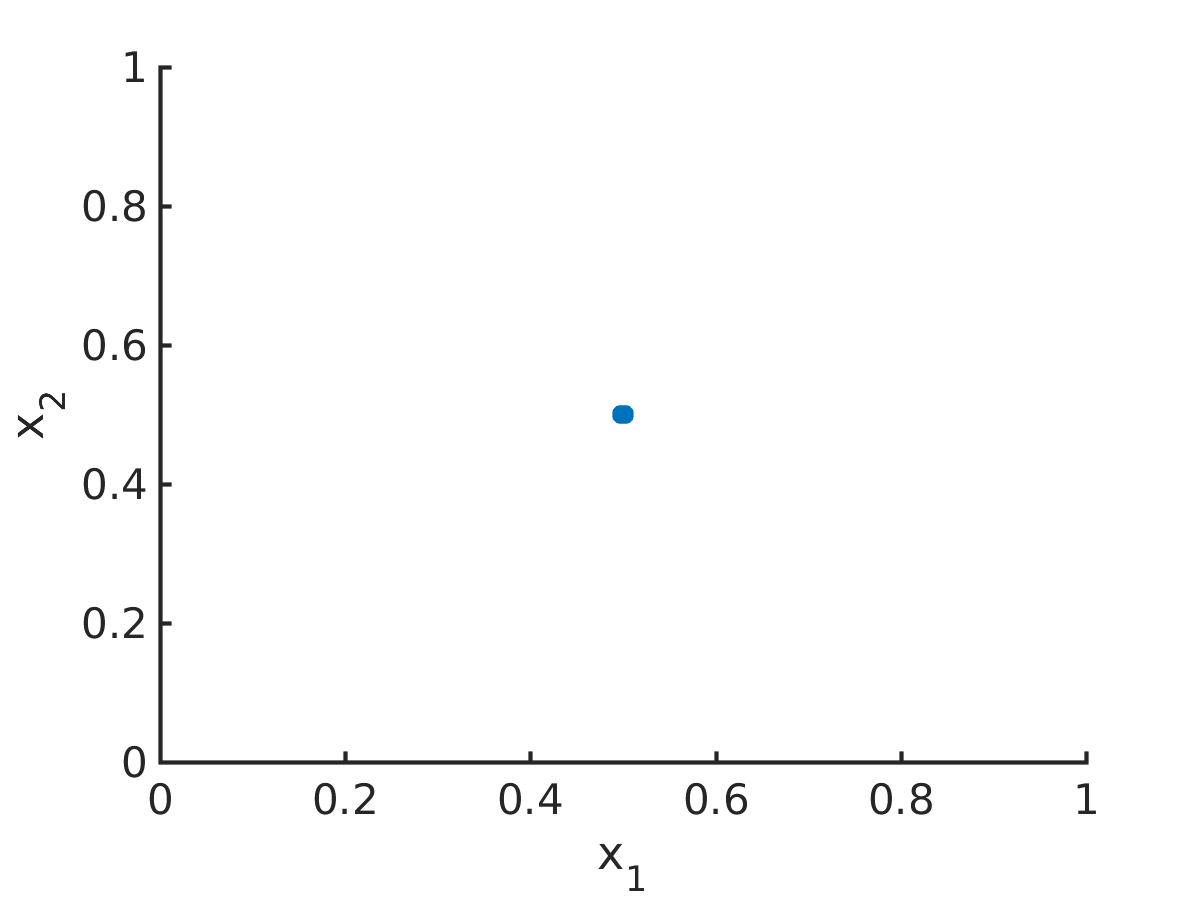}\label{fig:linearcoeffnumericsinitial}}
	\subfloat[${a_s=-0.2}$, ${b_s=0.1}, {R_c=0.3}$]{
		\includegraphics[width=0.32\textwidth]{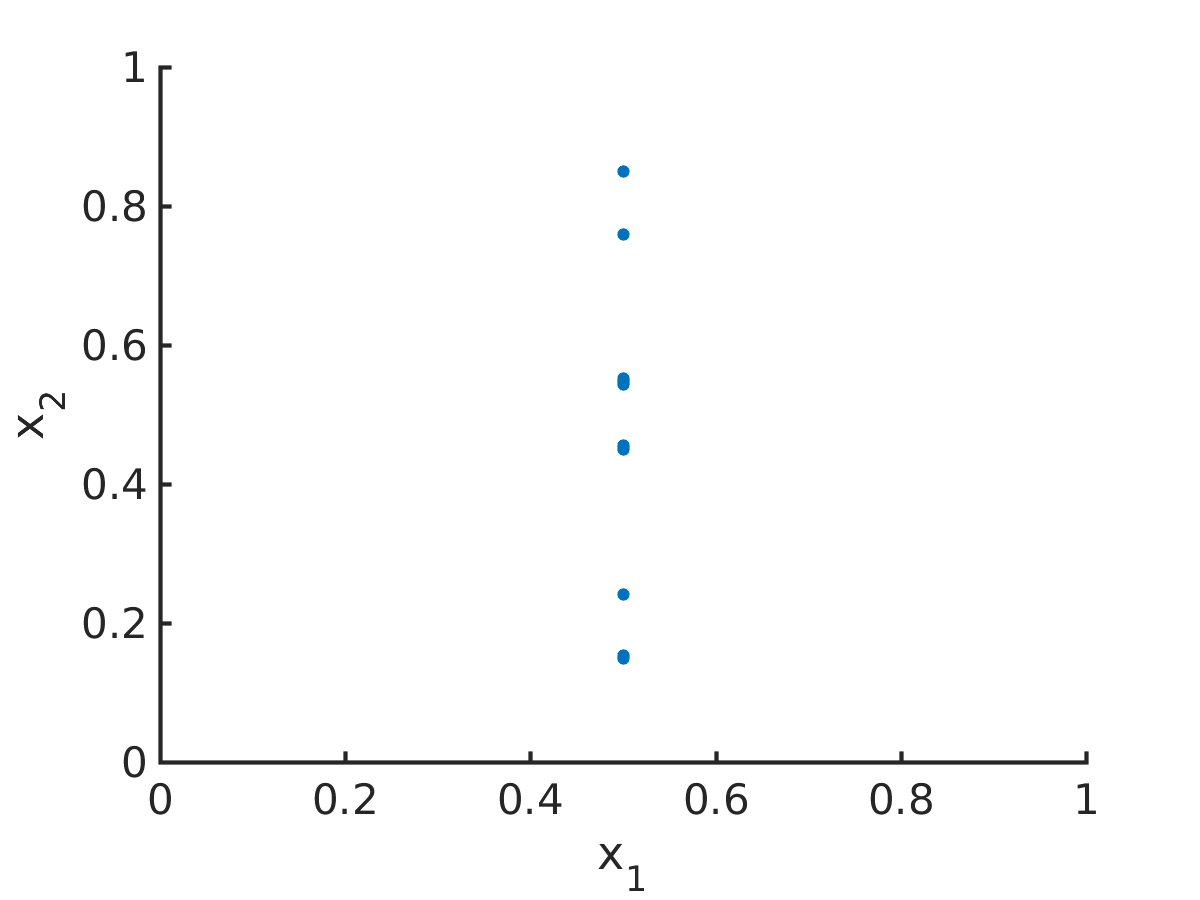}\label{fig:linearcoeffnumericssmallrc}}
	\subfloat[${a_s=-0.2}$, ${b_s=0.1}, {R_c=0.5}$]{\includegraphics[width=0.32\textwidth]{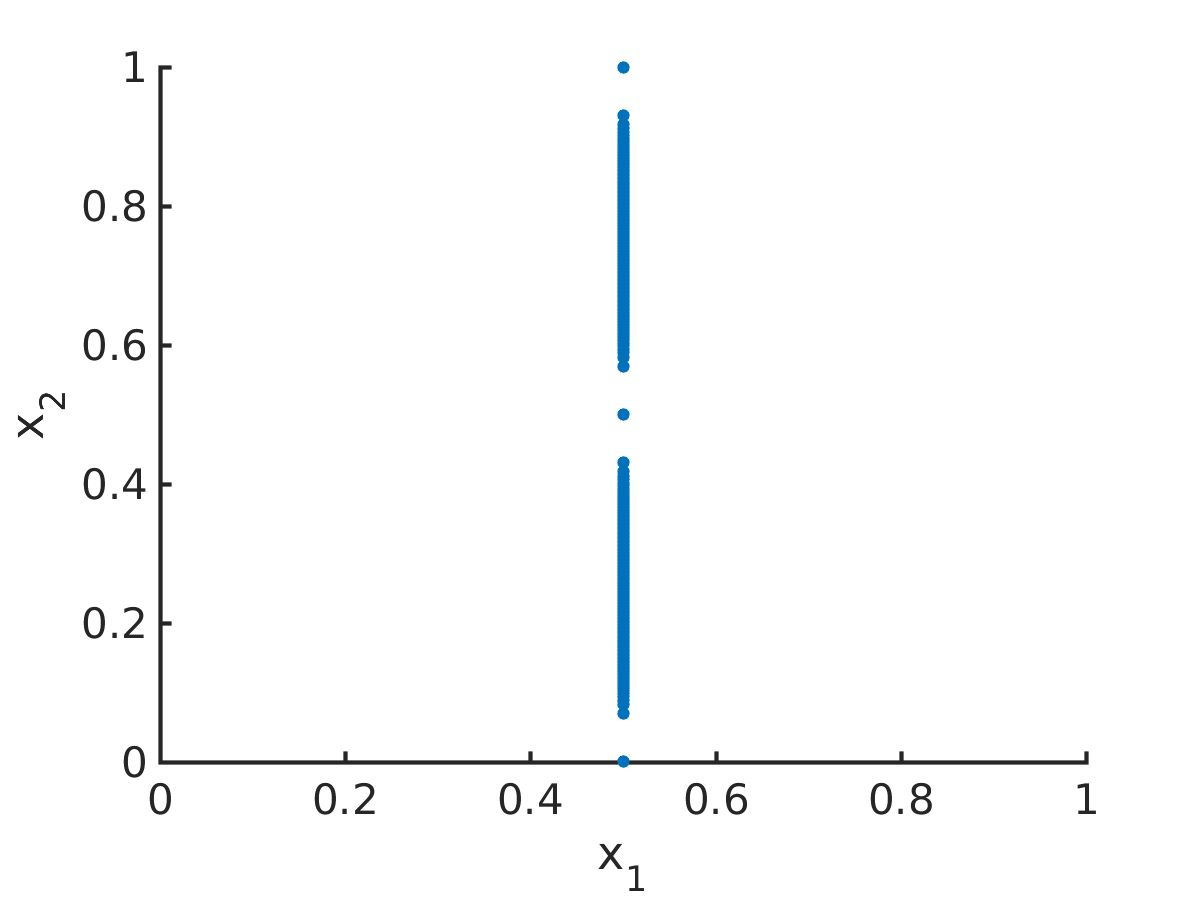}\label{fig:linearcoeffnumericseq1}}\\
	\subfloat[${a_s=-0.02}$, ${b_s=0.01, R_c=0.5}$]{\includegraphics[width=0.32\textwidth]{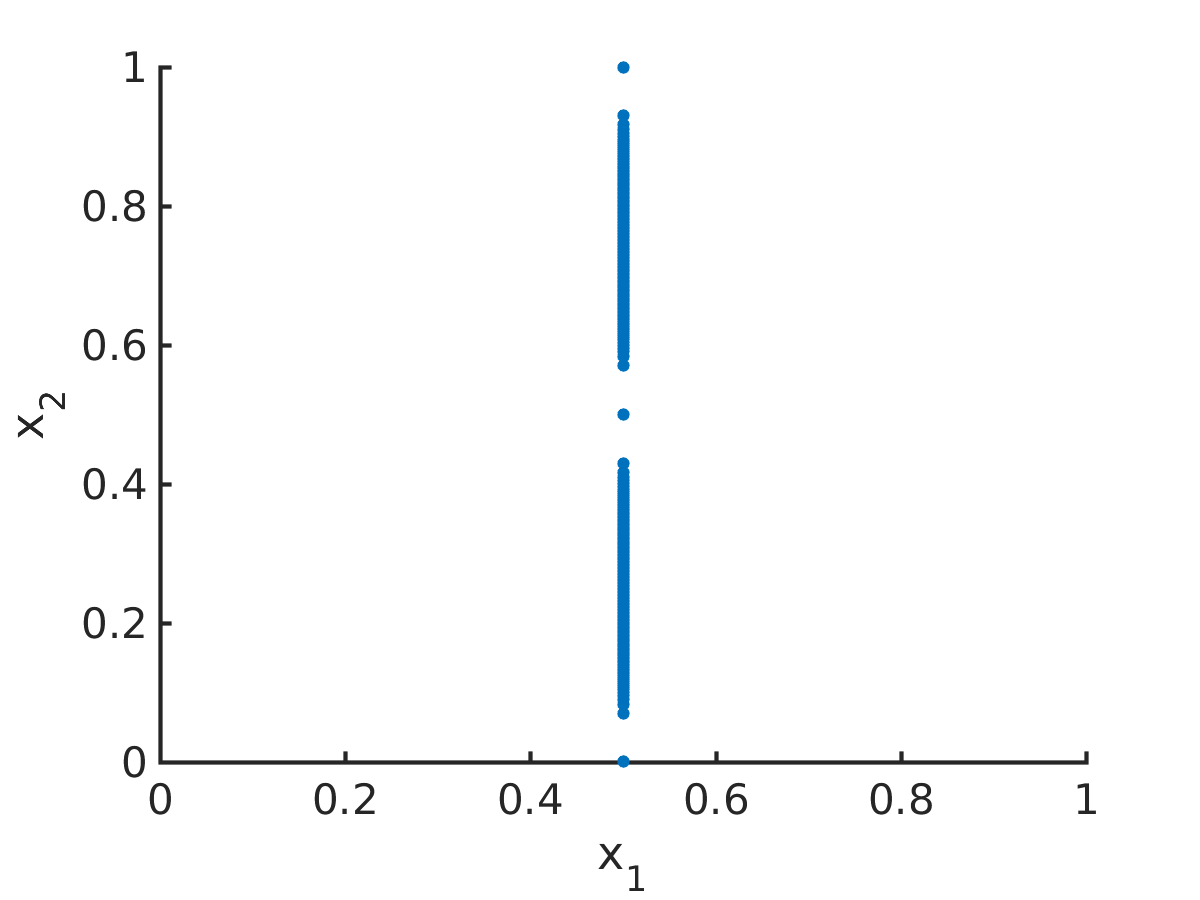}\label{fig:linearcoeffnumericseq2}}    
	\subfloat[${a_s=-0.1}$, ${b_s=0.1}, {R_c=0.5}$]{
		\includegraphics[width=0.32\textwidth]{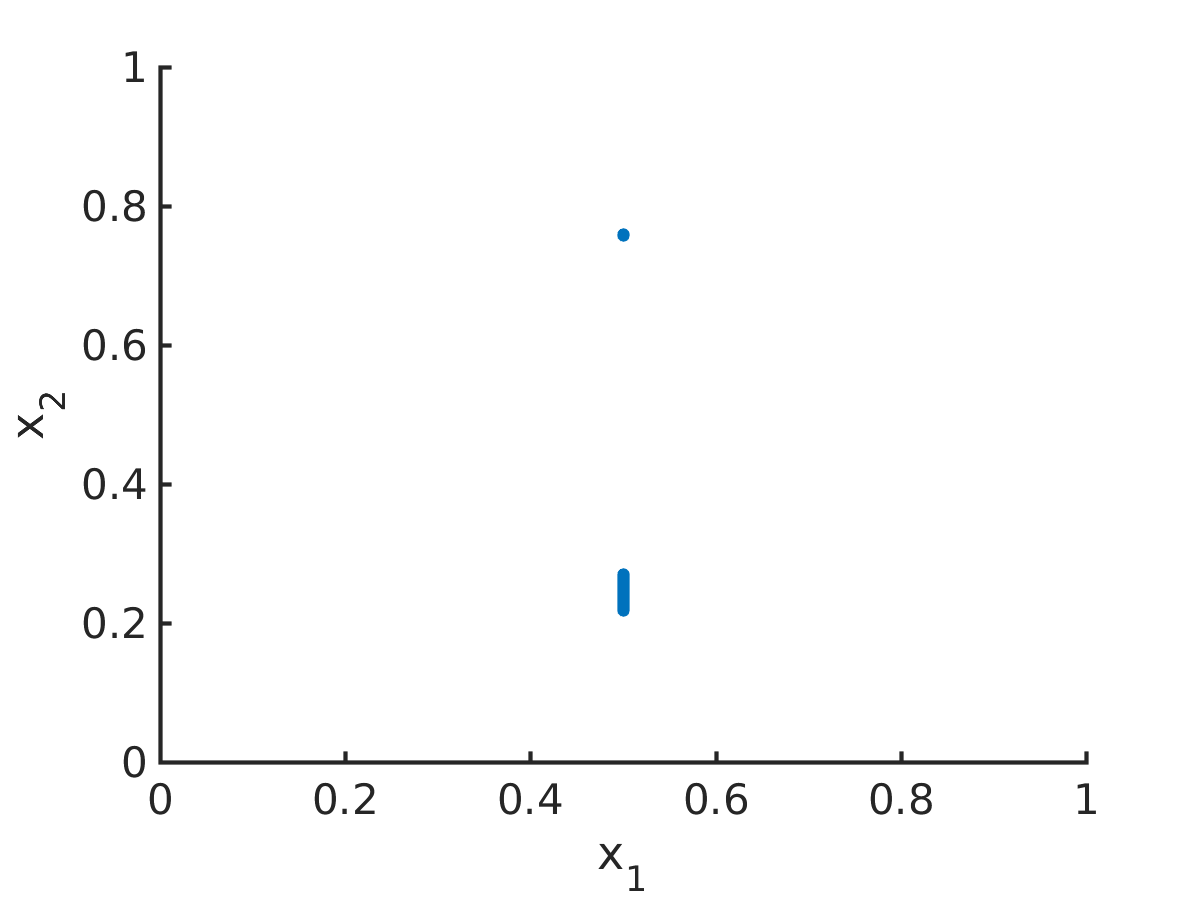}\label{fig:linearcoeffnumericssmaller}}
	\subfloat[${a_s=-0.4}$, ${b_s=0.1}, {R_c=0.5}$]{
		\includegraphics[width=0.32\textwidth]{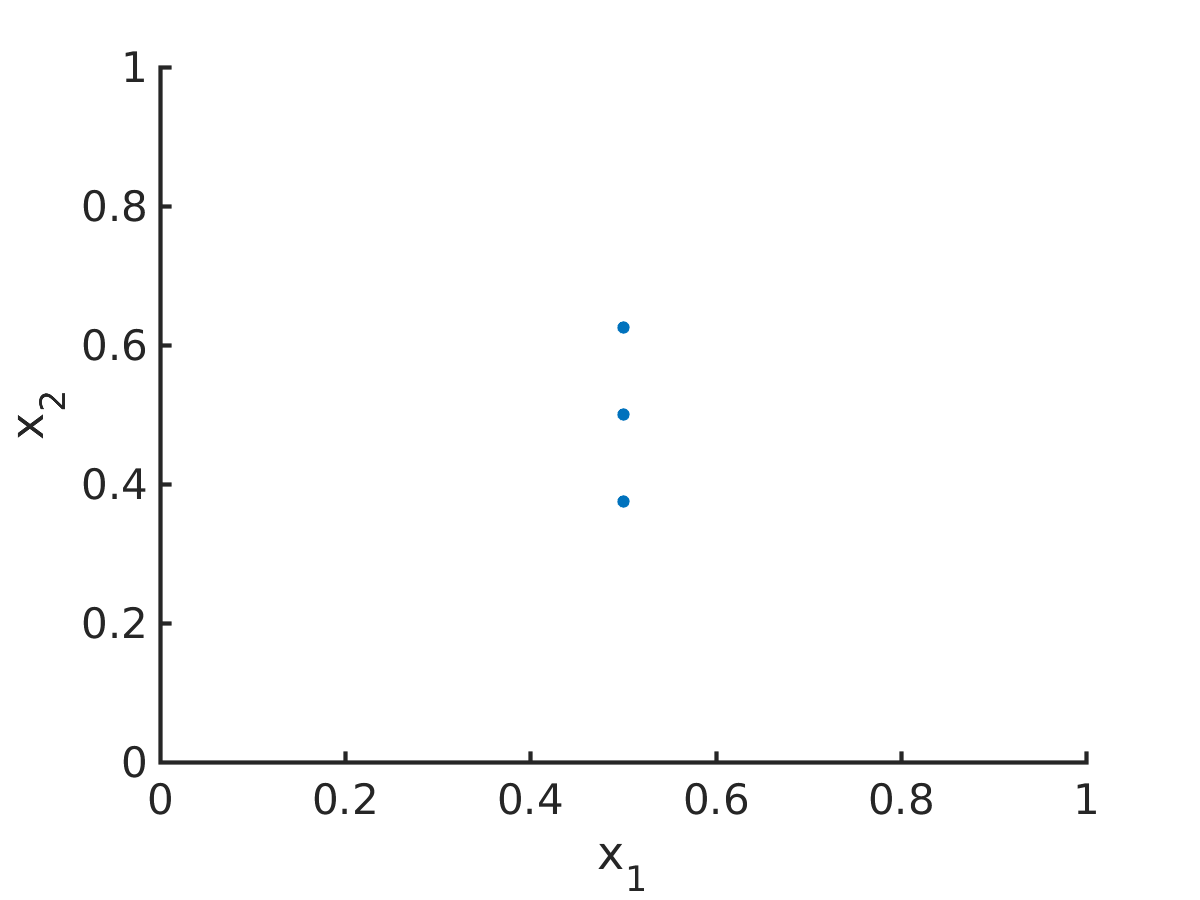}\label{fig:linearcoeffnumericslarger}}    
	\caption{Stationary solution to the model \eqref{eq:particlemodelperiodic} for total force \eqref{eq:totalforcehom} with linear force coefficients $f_l^\epsilon(|d|)=a_l|d|+b_l$, $f_s^\epsilon(|d|)=a_s|d|+b_s$ for $|d|\in[0,R_c-\epsilon]$ in \eqref{eq:linearcoefffunctions} with $a_l=-3,b_l=0.1$ and cutoff radius $R_c$ in the limit $\epsilon\to 0$.}\label{fig:linearcoeffnumerics}
\end{figure}

\begin{figure}[htbp]
	\centering
	\subfloat[${a_s=-0.2}$, ${b_s=0.1}, {R_c=0.3}$]{
		\includegraphics[width=0.32\textwidth]{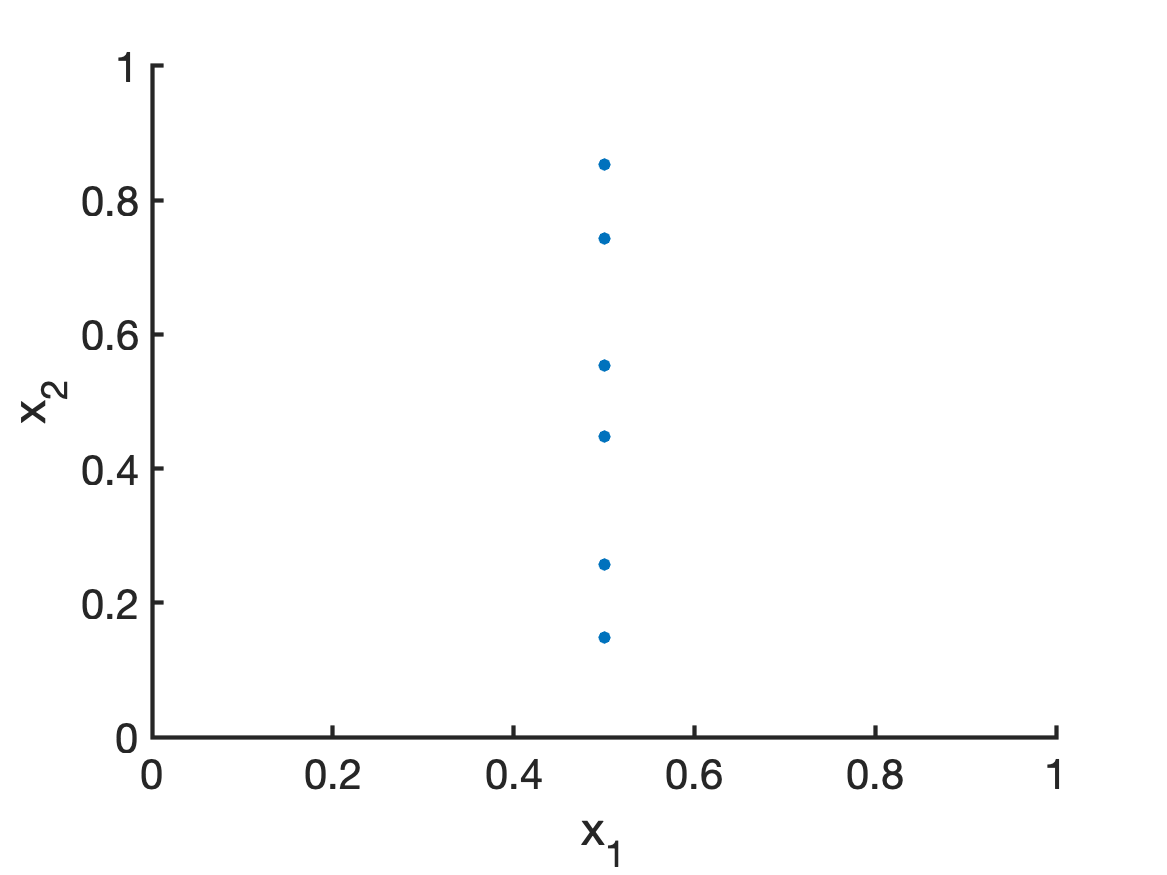}\label{fig:linearcoeffnumericssmallrceps}}
	\subfloat[${a_s=-0.2}$, ${b_s=0.1}, {R_c=0.5}$]{\includegraphics[width=0.32\textwidth]{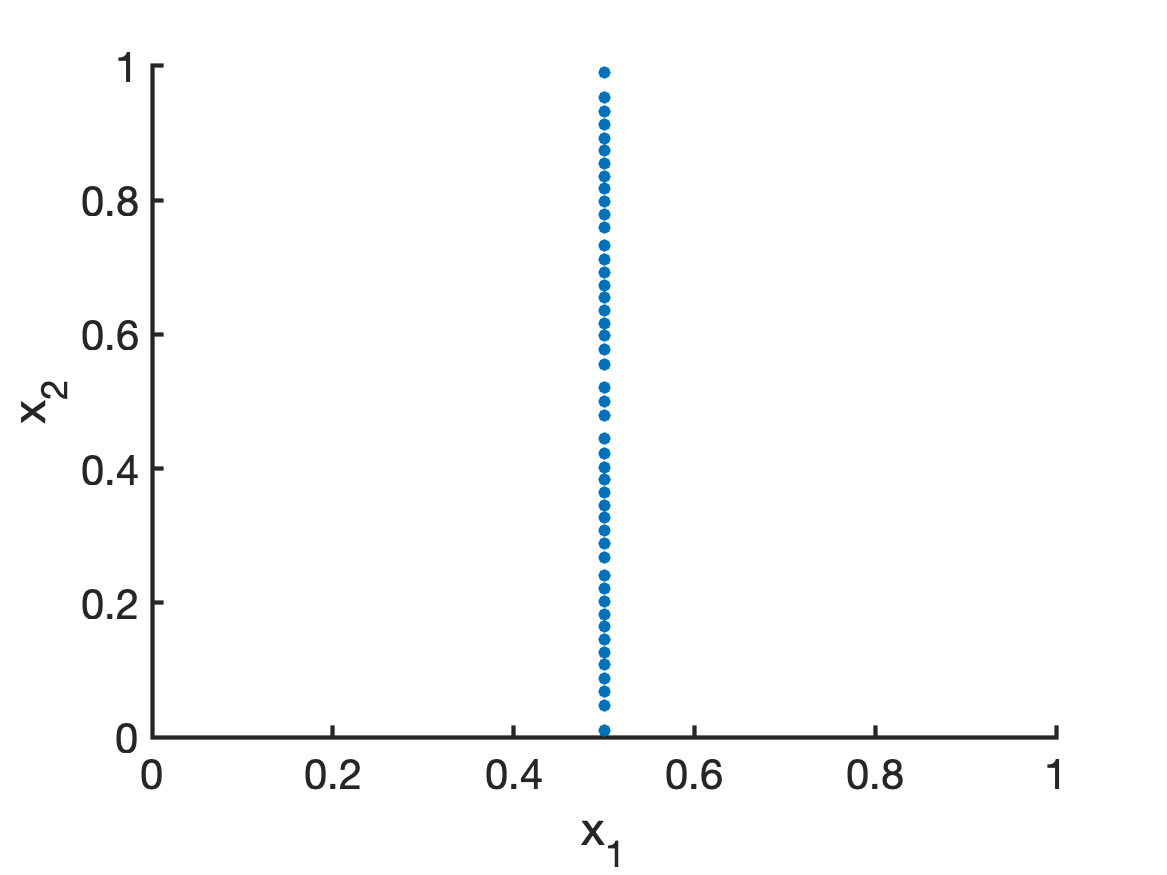}\label{fig:linearcoeffnumericseq1eps}}
	\subfloat[${a_s=-0.1}$, ${b_s=0.1}, {R_c=0.5}$]{
		\includegraphics[width=0.32\textwidth]{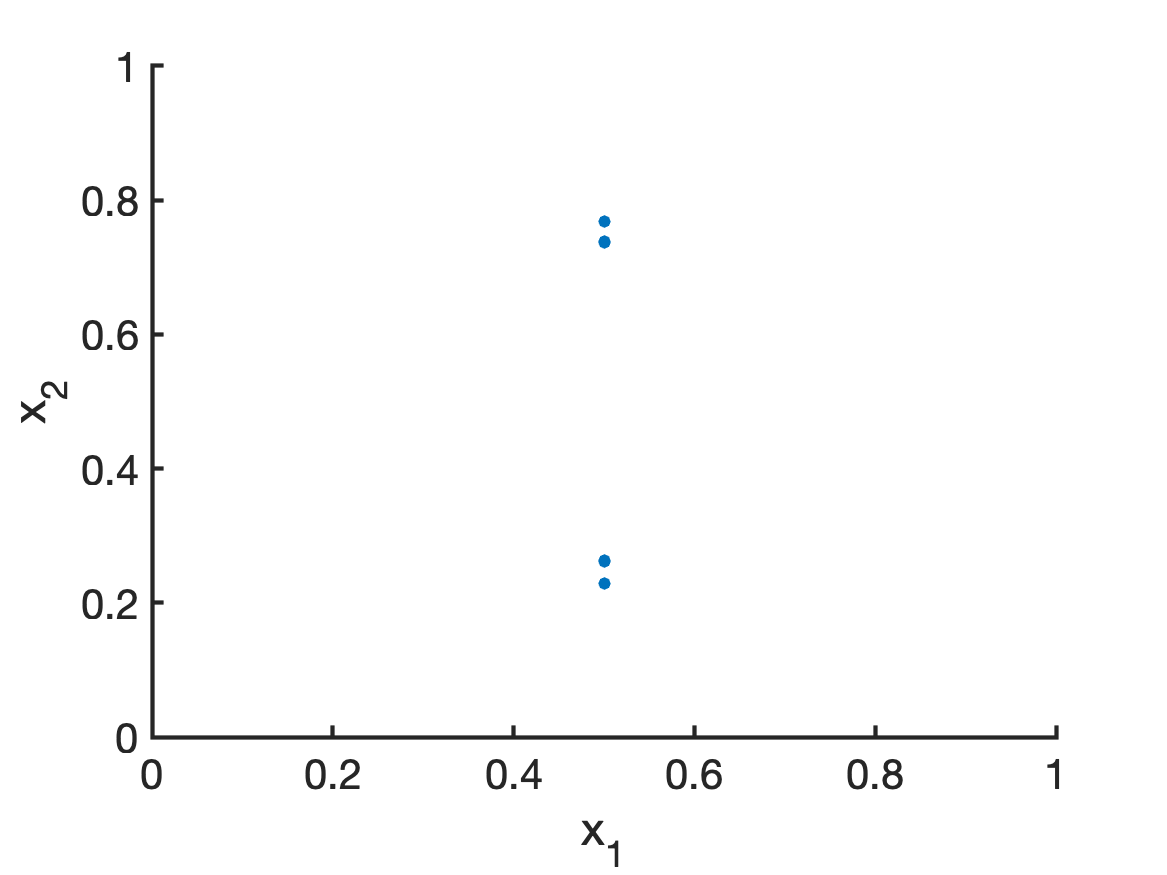}\label{fig:linearcoeffnumericssmallereps}}   
	\caption{Stationary solution to the model \eqref{eq:particlemodelperiodic} for total force \eqref{eq:totalforcehom} with linear force coefficients $f_l^\epsilon(|d|)=a_l|d|+b_l$, $f_s^\epsilon(|d|)=a_s|d|+b_s$ for $|d|\in[0,R_c-\epsilon]$ in \eqref{eq:linearcoefffunctions} with $a_l=-3,b_l=0.1$ and cutoff radius $R_c$ for $\epsilon=0.01$.}\label{fig:linearcoeffnumericseps}
\end{figure}

\begin{figure}[htbp]
	\centering
	\subfloat[Initial data]{\includegraphics[width=0.32\textwidth]{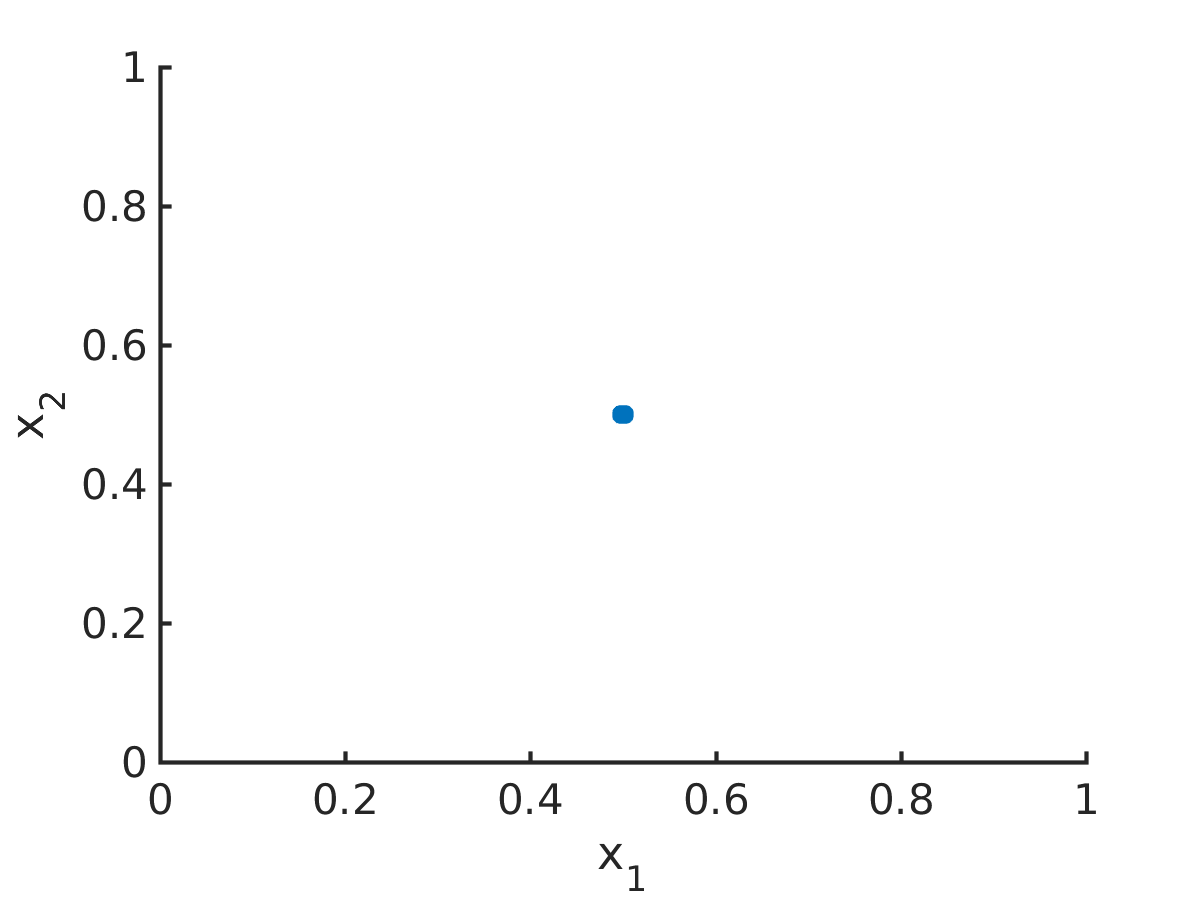}\label{fig:expcoeffnumericsinitial}}
	\subfloat[${e_s=10}$, ${R_c=0.1}$, {$f_l$ linear}]{
		\includegraphics[width=0.32\textwidth]{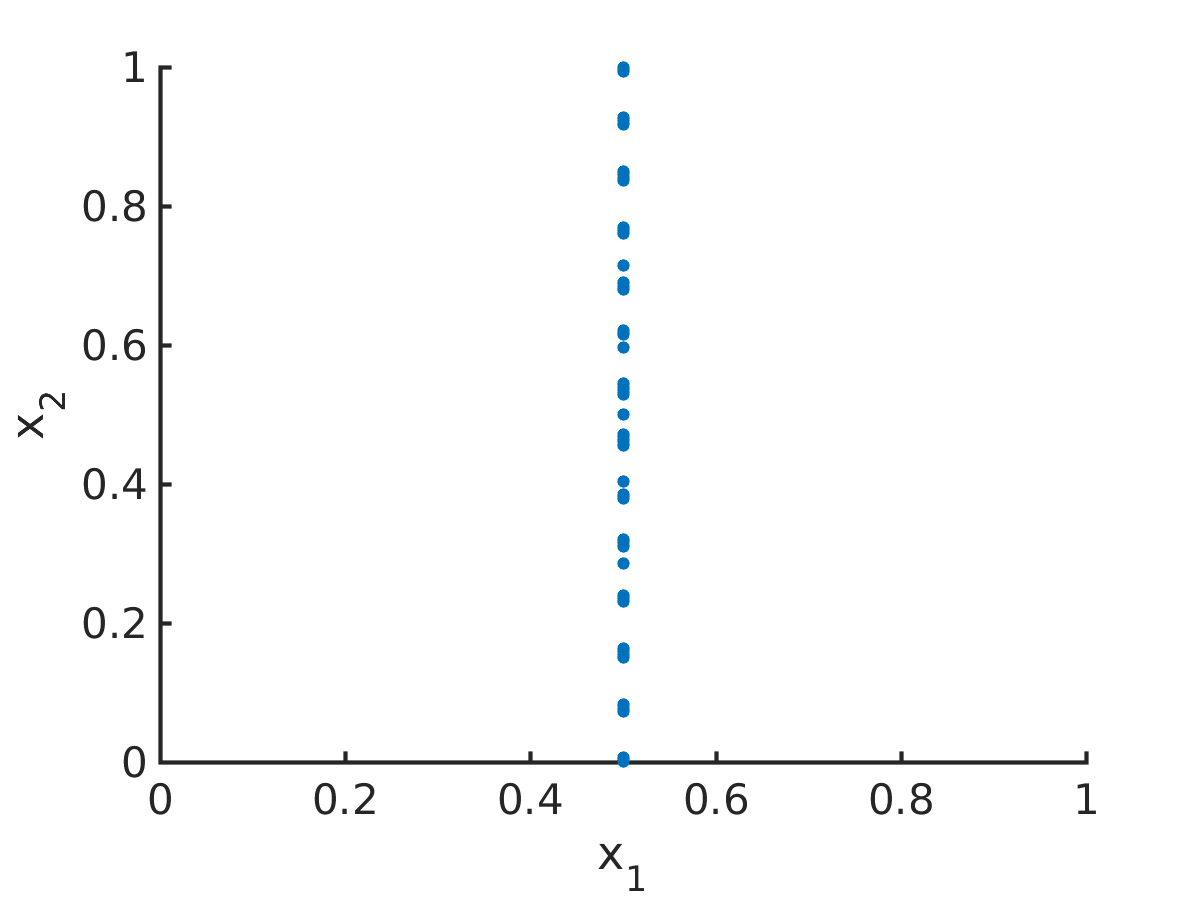}\label{fig:expcoeffnumericssmallessmallrcaddconst}}
	\subfloat[${e_s=10}$, ${R_c=0.5}$, {$f_l$ linear}]{
		\includegraphics[width=0.32\textwidth]{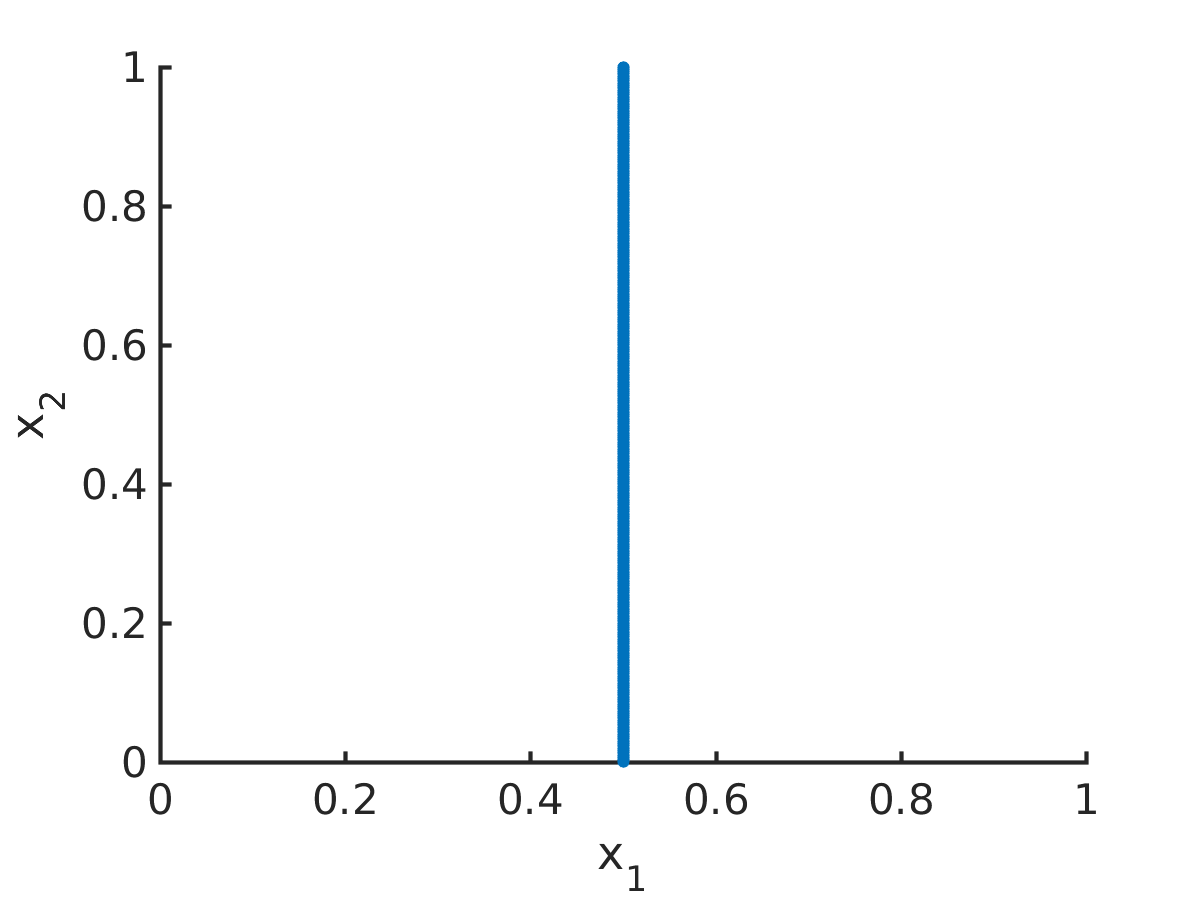}\label{fig:expcoeffnumericssmalleslargercaddconst}}\\
	\subfloat[${e_s=10}$, ${R_c=0.5}$, {$f_l$ linear, no add. constant for $f_s$}]{
		\includegraphics[width=0.32\textwidth]{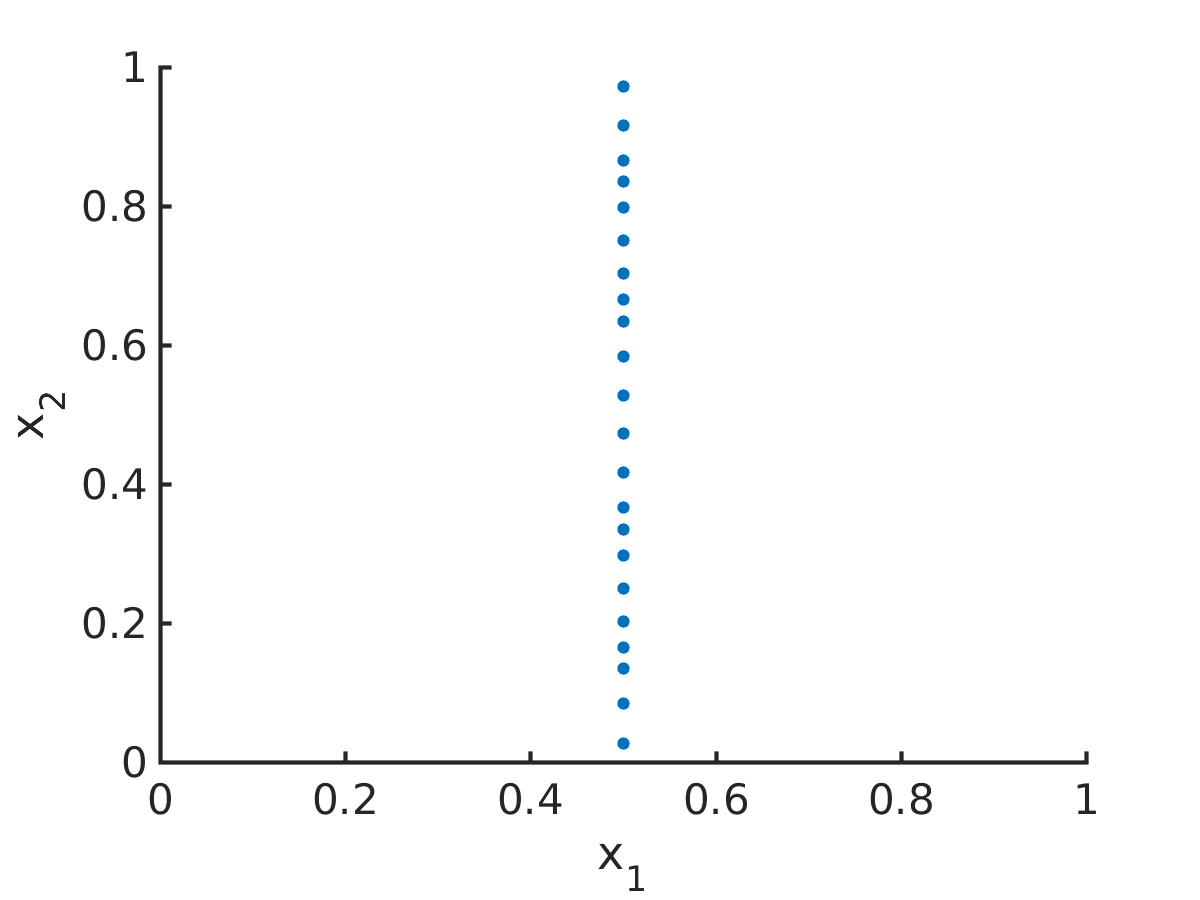}\label{fig:expcoeffnumericssmalleslargerc}}
	\subfloat[${e_s=100}$, $R_c=0.1$, {$f_l$ linear}]{\includegraphics[width=0.32\textwidth]{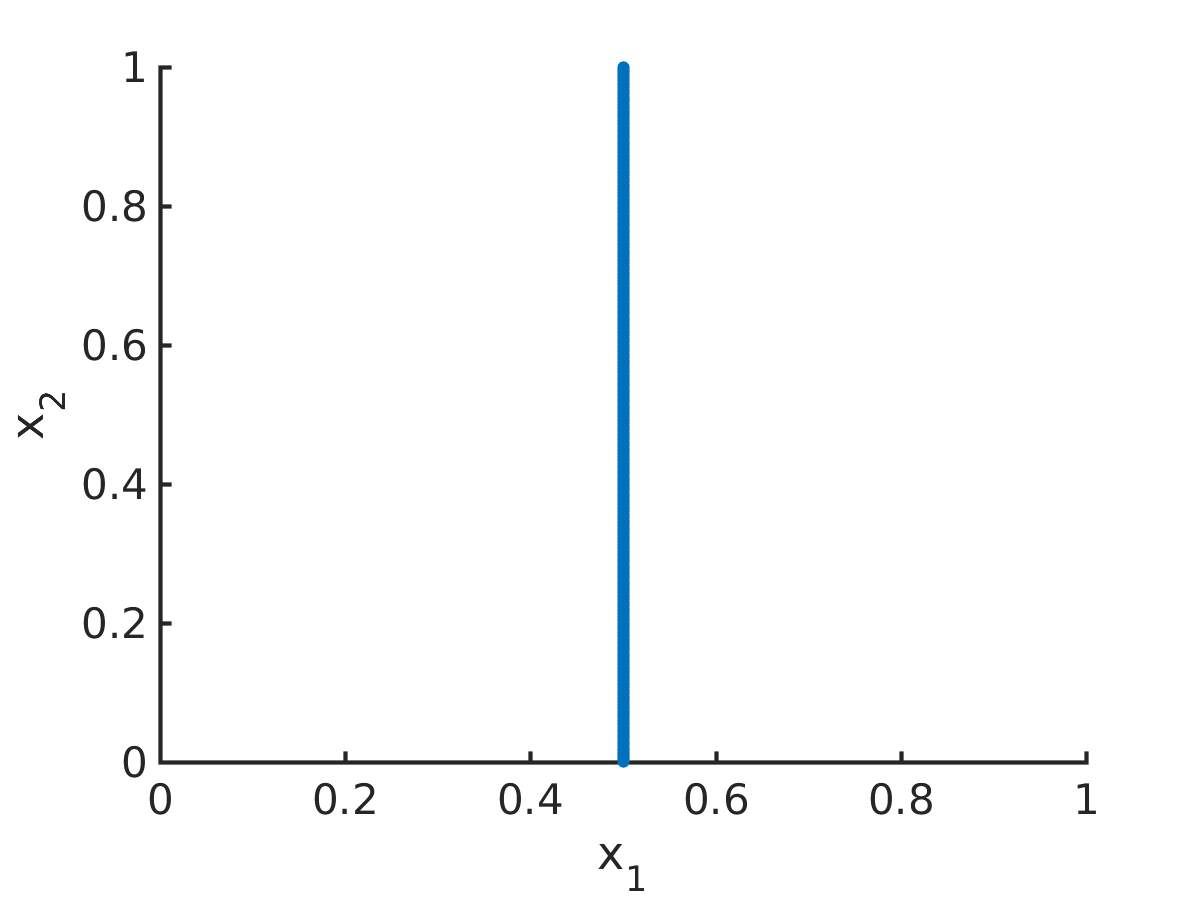}\label{fig:expcoeffnumericssmallrc}}
	\subfloat[${e_s=100}$, ${R_c=0.5}$, $f_l$ exp.]{
		\includegraphics[width=0.32\textwidth]{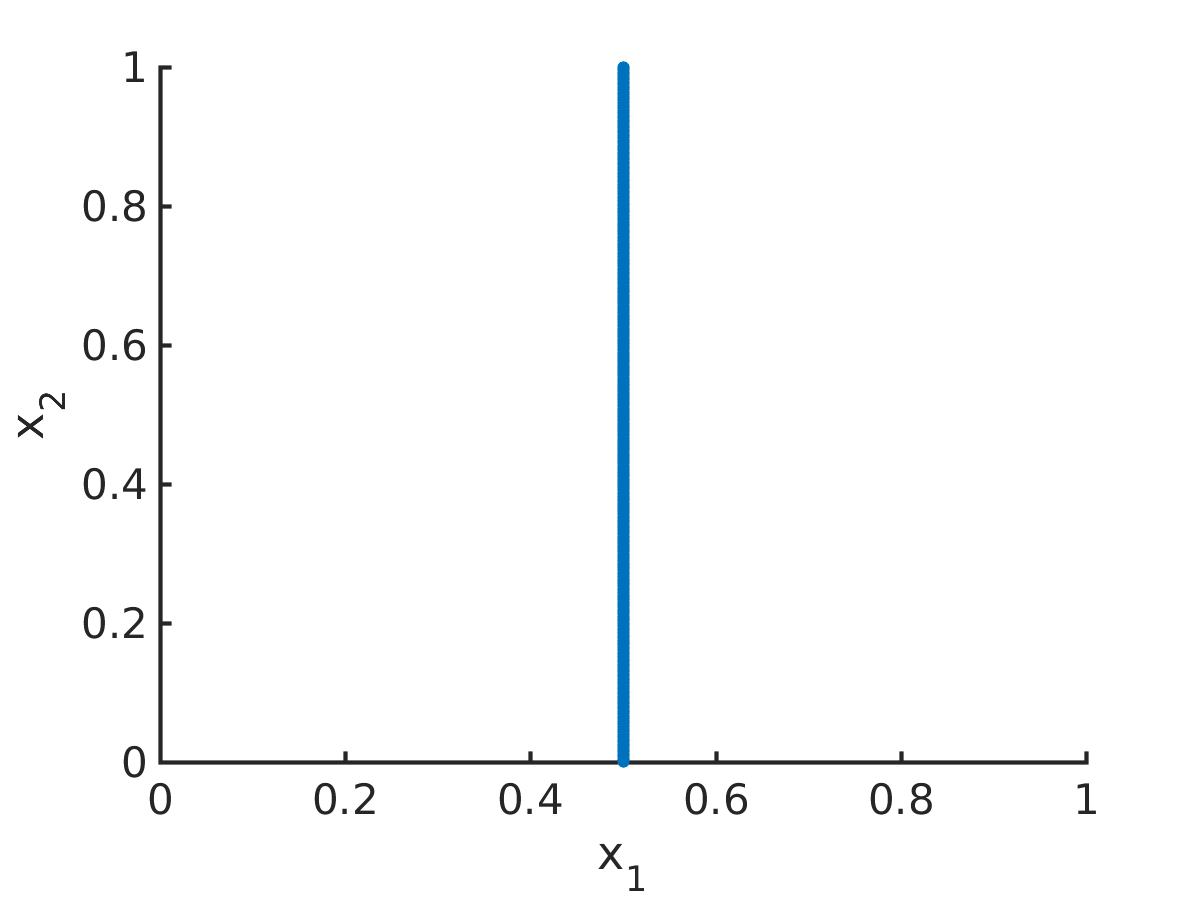}\label{fig:expcoeffnumericslargercexpl}}        
	\caption{Stationary solution to the model \eqref{eq:particlemodelperiodic} for total force \eqref{eq:totalforcehom} with exponential force coefficient $f_s^\epsilon(|d|)=c\exp(-e_s |d|)-c\exp(-e_s (R_c-\epsilon))$ for $|d|\in[0,R_c-\epsilon]$ along $s$, defined in  \eqref{eq:expforcedef},    and $f_l^\epsilon(|d|)=0.1-3|d|$ or $f_l^\epsilon(|d|)=0.13\exp(-100|d|)-0.03\exp(-10|d|)$ for $|d|\in[0,R_c-\epsilon]$ along $l$ with cutoff $R_c$ in the limit $\epsilon\to 0$.}\label{fig:expcoeffnumerics}
\end{figure}

\begin{figure}[htbp]
	\centering 
	\subfloat[Initial data]{\includegraphics[width=0.32\textwidth]{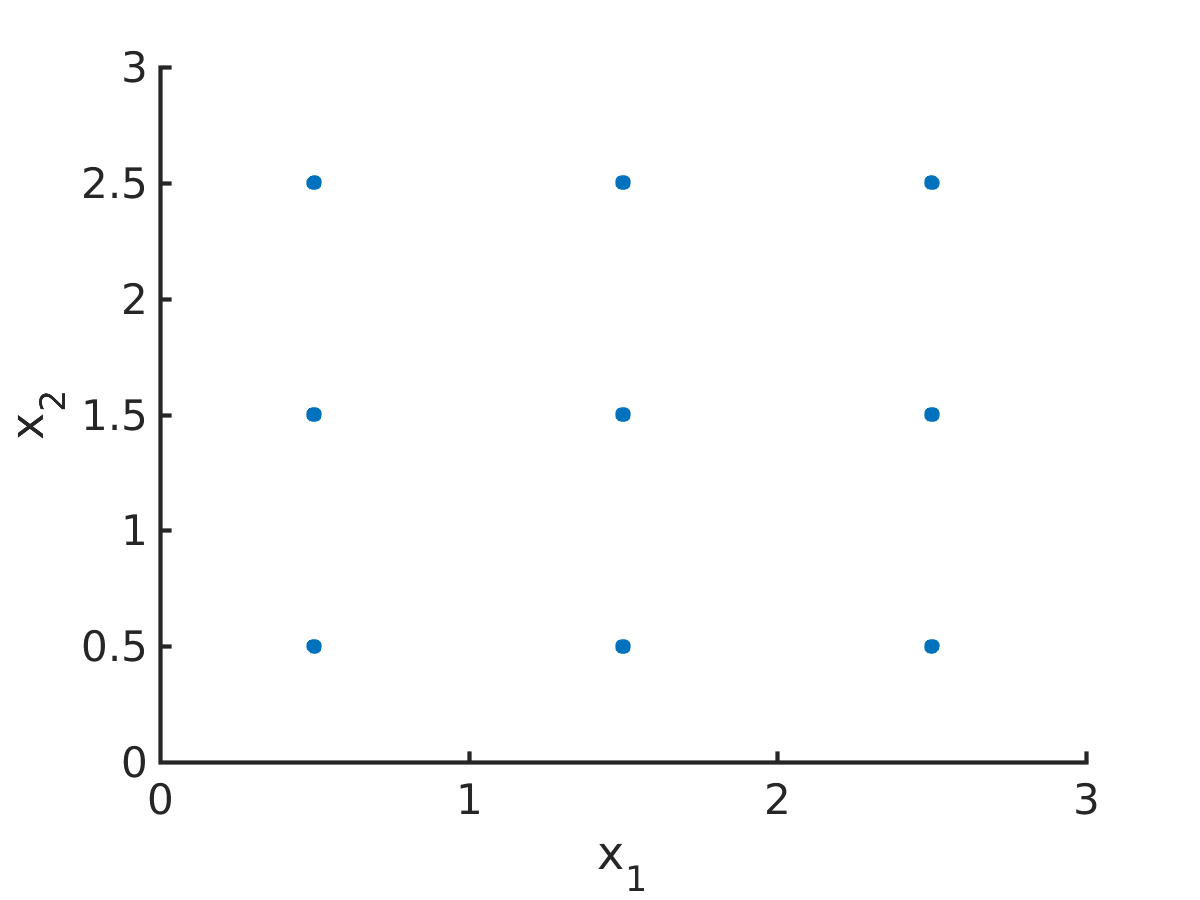}}
	\subfloat[$R_c=0.1$]{
		\includegraphics[width=0.32\textwidth]{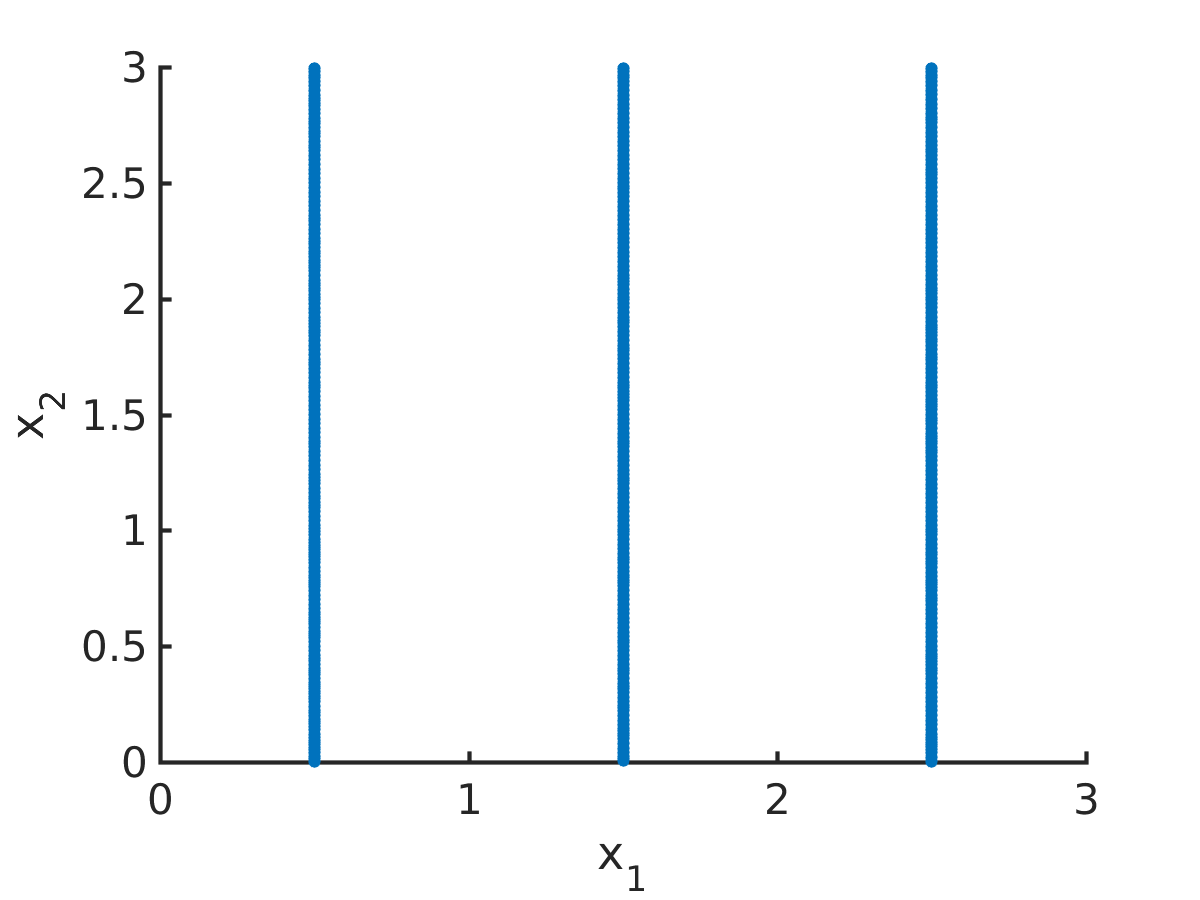}}       
	\caption{Stationary solution to the model \eqref{eq:particlemodelperiodic} for total force \eqref{eq:totalforcehom} with  exponential force coefficients $f_s^\epsilon(|d|)=c\exp(-e_s |d|)-c\exp(-e_s (R_c-\epsilon))$ in \eqref{eq:expforcedef}  and  ${f_l^\epsilon(|d|)=0.13\exp(-100|d|)-0.03\exp(-10|d|)}$ with cutoff $R_c$ on the domain $[0,3]^2$.}\label{fig:expcoeffnumericlargedomain}
\end{figure}

\begin{figure}[htbp]
	\centering 
	\subfloat[$s=(1,1)/\sqrt{2}$,  $ l=(-1,1)/\sqrt{2}$]{\includegraphics[width=0.32\textwidth]{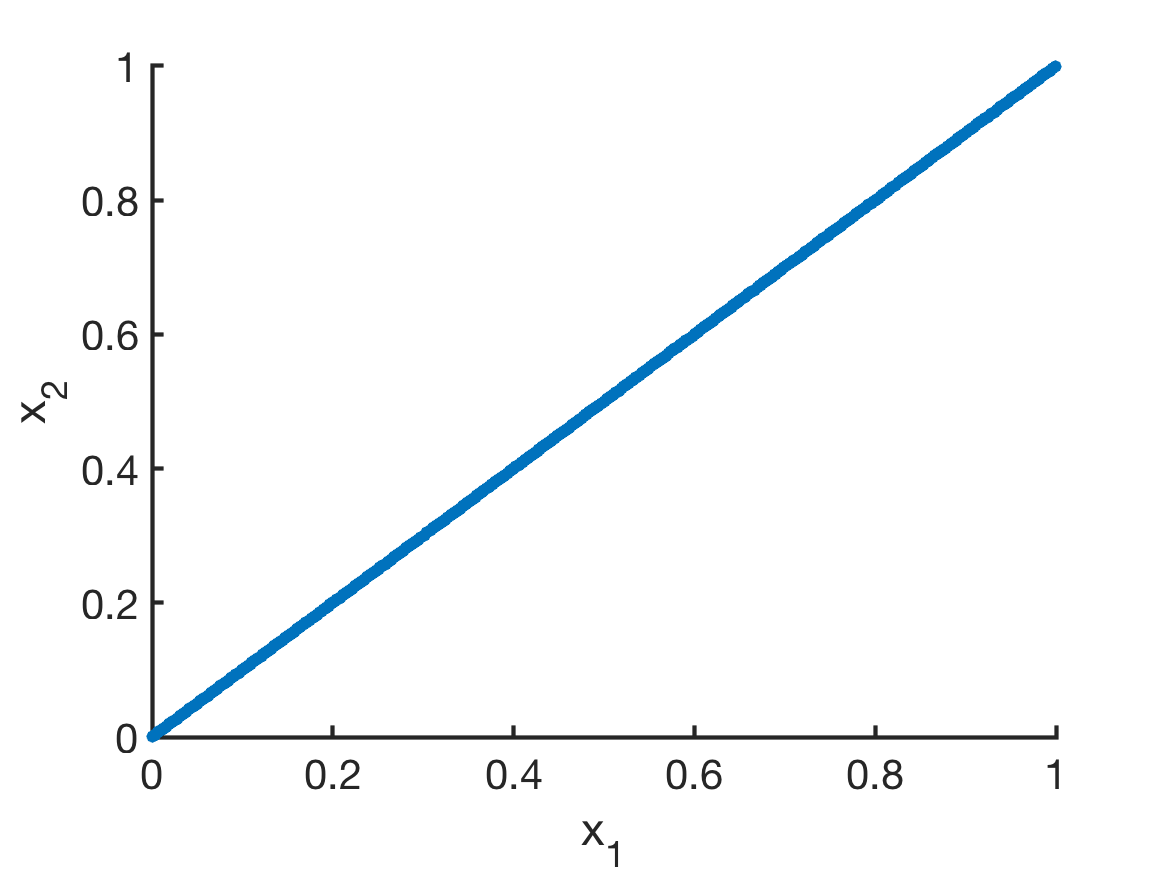}\label{fig:tensorfielddiag}}
	\subfloat[$s=(1,2)/\sqrt{5},l=(-2,1)/\sqrt{5}$]{	\includegraphics[width=0.32\textwidth]{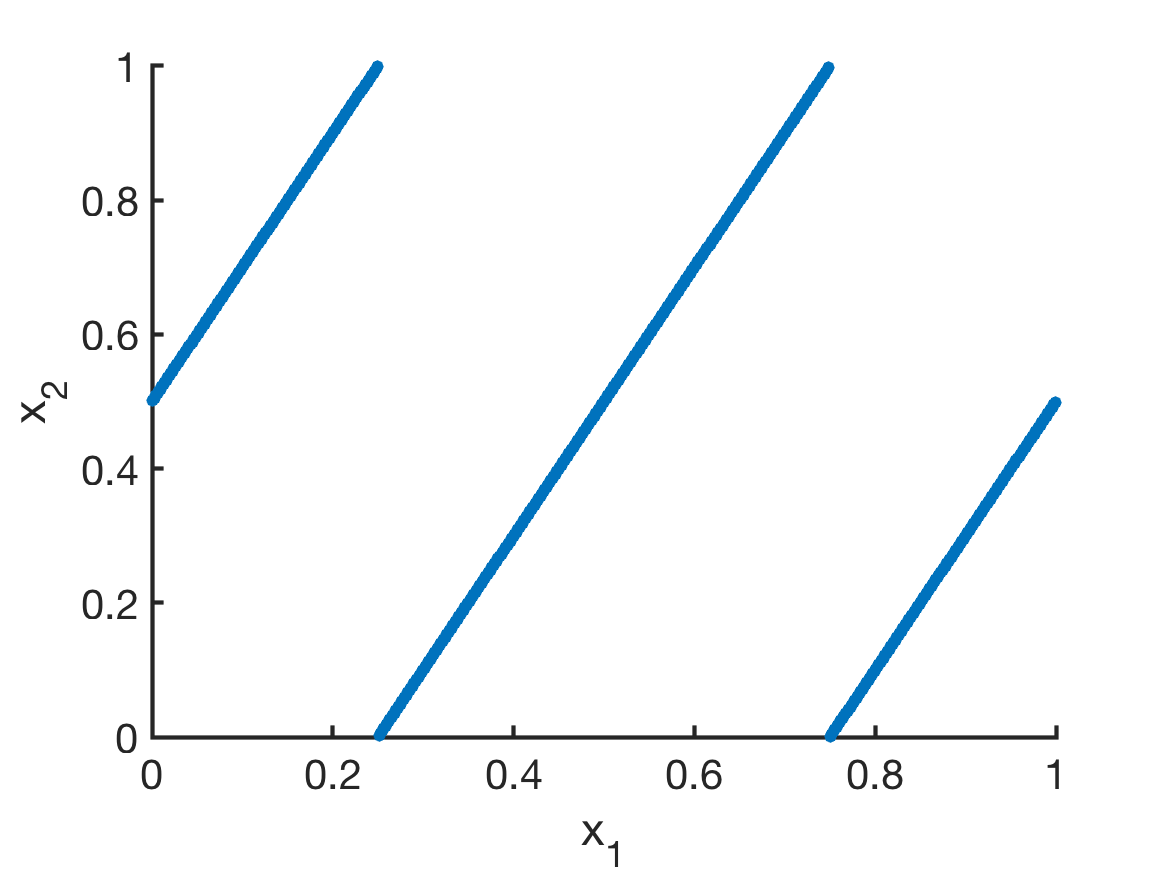}\label{fig:tensorfieldtwo}}    
	\subfloat[$s=(1,5)/\sqrt{26},l=(-5,1)/\sqrt{26}$]{	\includegraphics[width=0.32\textwidth]{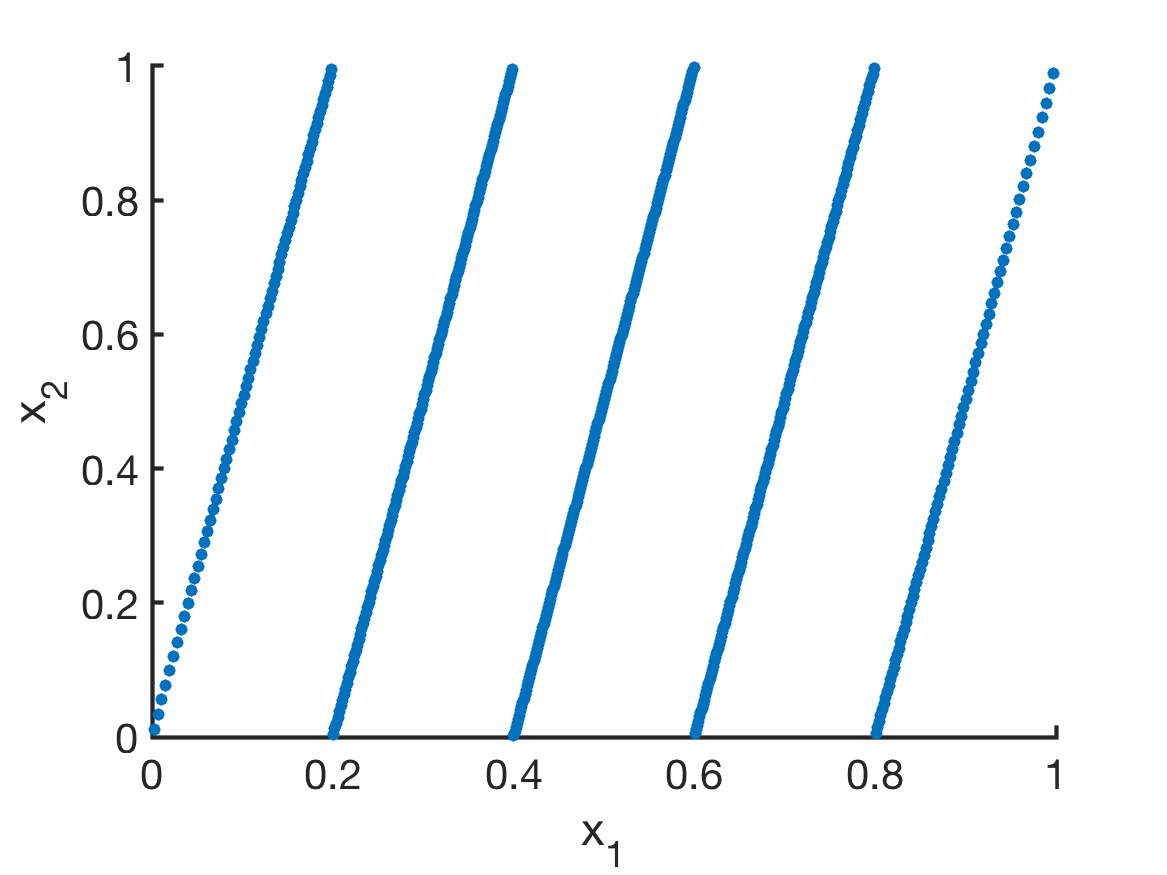}\label{fig:tensorfieldfive}}    
	\caption{Stationary solution to the model \eqref{eq:particlemodelperiodic} for different tensor fields $T$, given by $s,l$, and total force \eqref{eq:totalforcehom} with  exponential force coefficients $f_s^\epsilon(|d|)=c\exp(-e_s |d|)-c\exp(-e_s (R_c-\epsilon))$ for $|d|\in[0,R_c-\epsilon]$ in \eqref{eq:expforcedef}  and  ${f_l^\epsilon(|d|)=0.13\exp(-100|d|)-0.03\exp(-10|d|)}$ for $|d|\in[0,R_c-\epsilon]$ with cutoff $R_c=0.1$ in the limit $\epsilon\to 0$.\label{fig:expcoeffrotatedtensor}}
\end{figure}

\begin{figure}[htbp]
	\centering 
	\subfloat[$e_{l_1}=20$]{\includegraphics[width=0.32\textwidth]{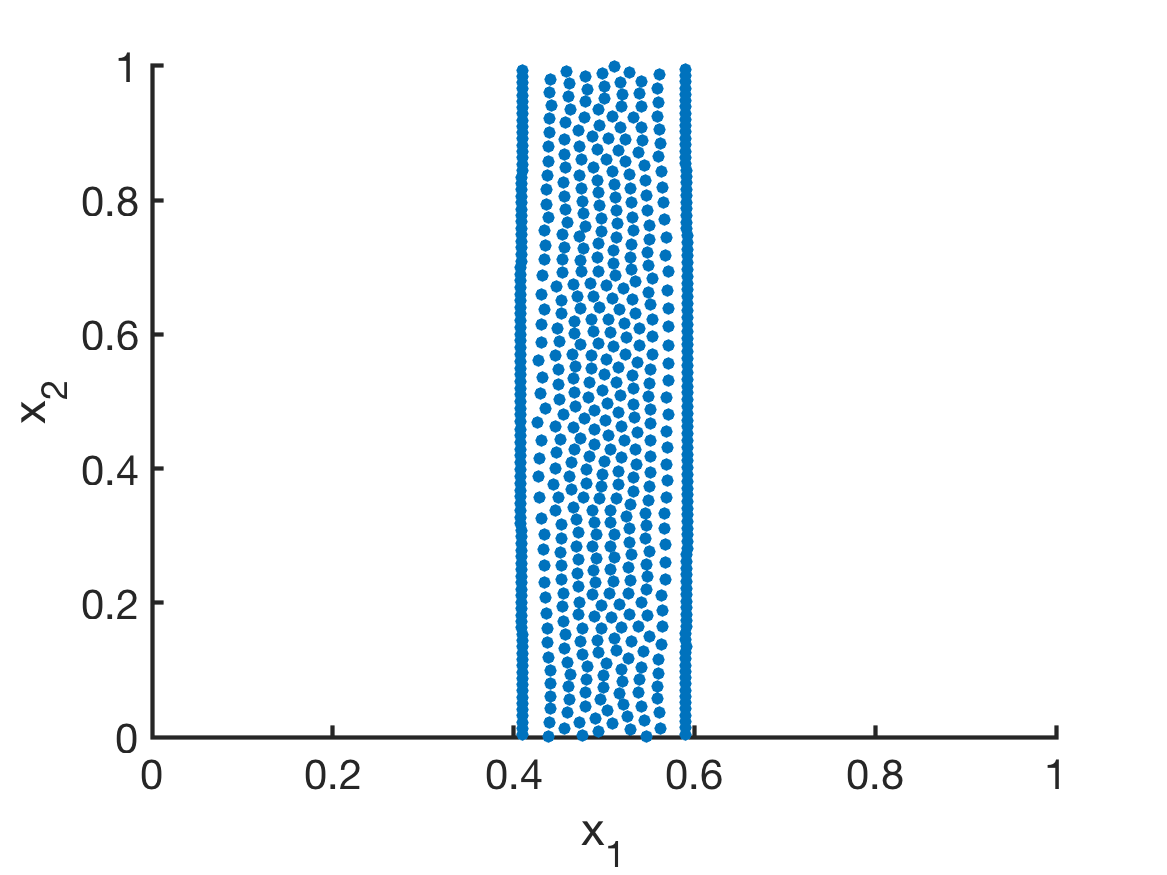}}
	\subfloat[$e_{l_1}=30$]{	\includegraphics[width=0.32\textwidth]{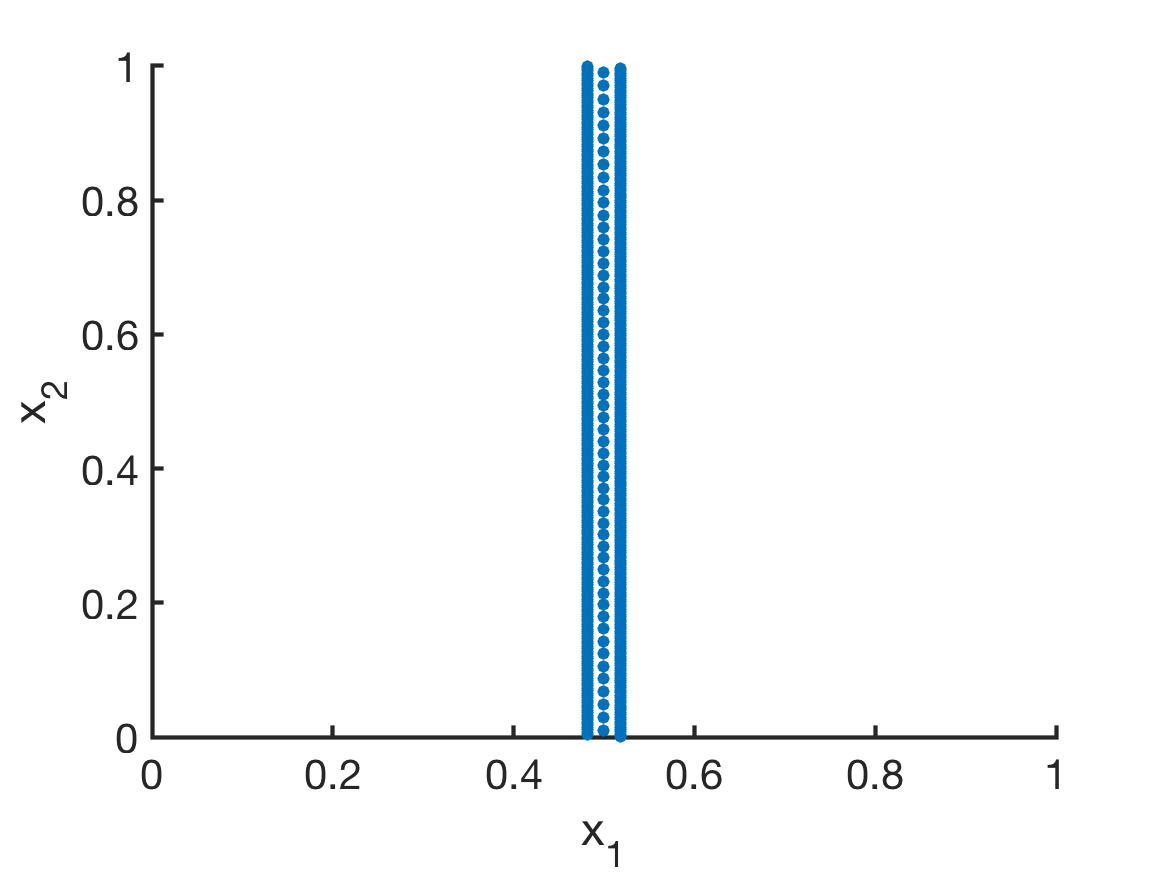}}    
	\subfloat[$e_{l_1}=50$]{	\includegraphics[width=0.32\textwidth]{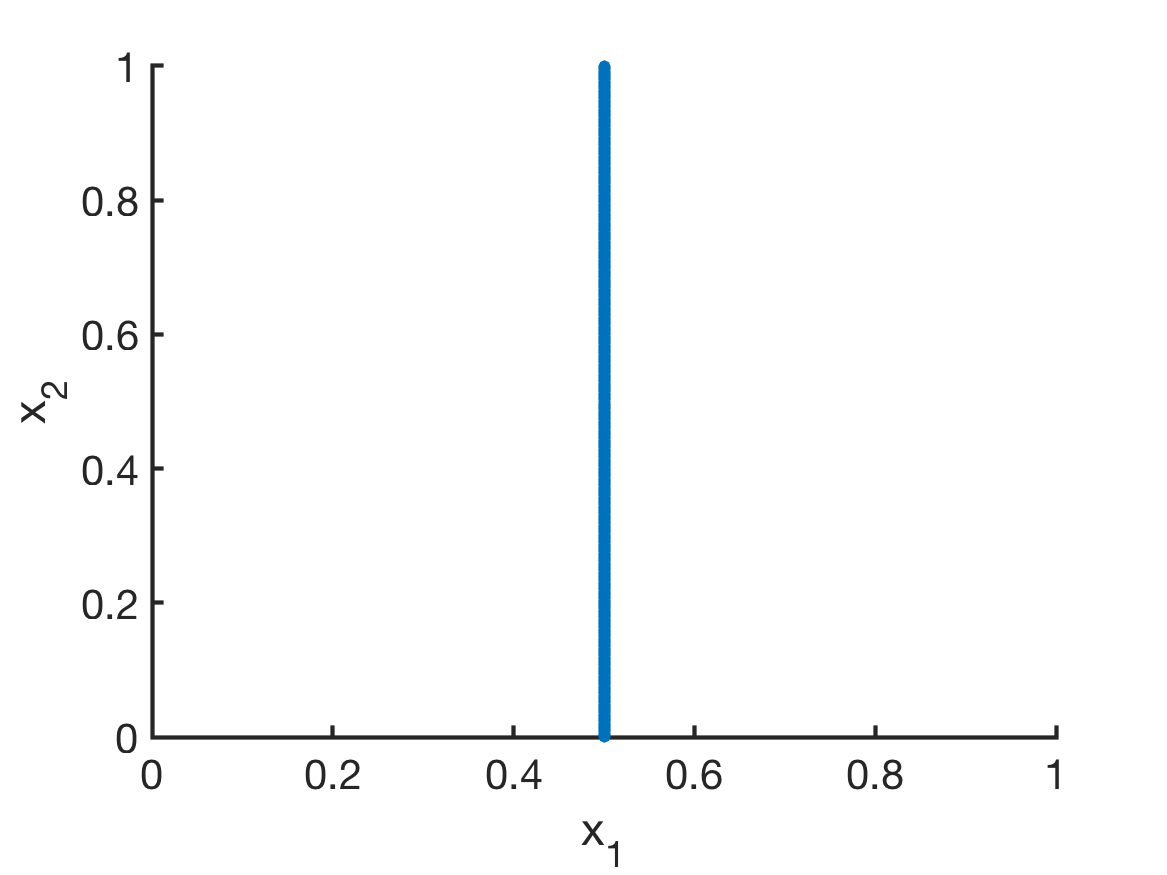}}    
	\caption{Stationary solution to the model \eqref{eq:particlemodelperiodic} for tensor field $T$ with $s=(0,1),l=(1,0)$ and  total force \eqref{eq:totalforcehom} with  exponential force coefficients $f_s^\epsilon(|d|)=c\exp(-e_s |d|)-c\exp(-e_s (R_c-\epsilon))$ for $|d|\in[0,R_c-\epsilon]$, defined in \eqref{eq:expforcedef},  and  ${f_l^\epsilon(|d|)=0.13\exp(-e_{l_1}|d|)-0.03\exp(-10|d|)}$ for $|d|\in[0,R_c-\epsilon]$ with cutoff $R_c=0.5$.}\label{fig:expcoeffrepulsion}
\end{figure}

Numerical results are shown in Figures \ref{fig:linearcoeffnumerics}--\ref{fig:expcoeffrepulsion}. For all numerical simulations we consider $N=600$ particles which are initially  equiangular distributed on a circle with centre $(0.5,0.5)$ and radius $0.005$ as illustrated in Figure \subref*{fig:linearcoeffnumericsinitial}. The stationary solution for the linear force coefficient $f_s^\epsilon$ in \eqref{eq:linearcoefffunctions}, i.e.
\begin{align*}
f_s^\epsilon (|d|)=a_s|d|+b_s, \qquad f_l^\epsilon(|d|)=0.1-3|d|, \qquad |d|\in[0,R_c-\epsilon],
\end{align*}
for different values of $a_s,b_s$, is shown in Figure \ref{fig:linearcoeffnumerics} in the limit $\epsilon\to 0$. As proven in Section~\ref{sec:linearcase}  equidistantly distributed particles along the vertical straight line form an unstable steady state for $N\in\N$ sufficiently large for $R_c\in(0,0.5)$. Hence, the stationary solutions are no lines of uniformly distributed particles and we obtain different clusters or line patterns instead. In Figure \subref*{fig:linearcoeffnumericssmallrc}, we consider $R_c=0.3$, resulting in clusters of particles along the vertical axis. For $R_c=0.5$ and $a_s,b_s$ chosen as $a_s=-\frac{b_s}{R_c}$, the requirement in \eqref{eq:linearcoeffcond} for the necessary stability condition to be satisfied with equality,  the real part of one of the eigenvalues of the stability matrix is equal to zero. The resulting steady states is shown for different scalings of the parameters $a_s,b_s$ in Figures \subref*{fig:linearcoeffnumericseq1} and \subref*{fig:linearcoeffnumericseq2}. One can see  that the particles align along a vertical line along the entire interval $[0,1]$, but are not equidistantly distributed along the vertical axis and thus the vertical straight line is an unstable steady state for any $N\in\N$ sufficiently large. For $a_s>-\frac{b_s}{R_c}$ and $a_s<-\frac{b_s}{R_c}$, respectively, with $R_c=0.5$ the corresponding steady states are shown in Figures  \subref*{fig:linearcoeffnumericssmaller} and \subref*{fig:linearcoeffnumericslarger}, resulting in clusters along the vertical axis.

	In Figure \ref{fig:linearcoeffnumericseps}, we consider the linear force coefficient $f_s^\epsilon$ in \eqref{eq:linearcoefffunctions} for different values of $a_s,b_s$ and $R_c$ where $\epsilon=0.01$ is fixed in contrast to $\epsilon\to 0$ in Figure \ref{fig:linearcoeffnumerics}, i.e.\ we consider the total force \eqref{eq:totalforcehom} with linear force coefficients $f_l^\epsilon(|d|)=a_l|d|+b_l$, $f_s^\epsilon(|d|)=a_s|d|+b_s$ for $|d|\in[0,R_c-\epsilon]$ in \eqref{eq:linearcoefffunctions} with $a_l=-3,b_l=0.1$. In Figure \subref*{fig:linearcoeffnumericssmallrceps}, we consider the same parameter values as in Figure \subref*{fig:linearcoeffnumericssmallrc}, i.e.\ $a_s=-0.2, b_s=0.1$ and $R_c=0.3$, resulting in the same  stationary solution for $\epsilon=0.01$ and $\epsilon\to 0$. In particular, the straight vertical line is unstable both for $\epsilon=0.01$ and $\epsilon\to 0$. For cutoff radius $R_c=0.5$, we obtain different stationary solutions for $\epsilon=0.01$ and $\epsilon\to 0$. In Figure \subref*{fig:linearcoeffnumericseq1eps}, we show the stationary solutions for  $a_s=-0.2, b_s=0.1$ and $R_c=0.5$ as in Figure \subref*{fig:linearcoeffnumericseq1}, i.e.\ $a_s=-\frac{b_s}{R_c}$. Even though stability/instability could not be determined analytically  the numerical results illustrate that straight vertical line is unstable both for $\epsilon=0.01$ and $\epsilon\to 0$.  The stationary solution for $a_s=-0.1, b_s=0.1$ and $R_c=0.5$ is shown in Figure \subref*{fig:linearcoeffnumericssmallereps} for $\epsilon=0.01$ and in Figure \subref*{fig:linearcoeffnumericssmaller} for $\epsilon\to 0$. Our analytical results show that the stationary solution is unstable in this case which is also consistent with the numerical results. In particular, we obtain the same instability results for $\epsilon=0.01$ as in Figure \ref{fig:linearcoeffnumerics} where the limit $\epsilon\to 0$ is considered.

For the exponentially decaying force coefficient $f_s^\epsilon$ along $s$ in  \eqref{eq:expforcedef}, given by
\begin{align*}
f_s^\epsilon(|d|)=c\exp(-e_s |d|)-c\exp(-e_s (R_c-\epsilon)), \qquad |d|\in[0,R_c-\epsilon],
\end{align*} 
for $\epsilon>0$, we consider the parameter values $c=0.1$ and $e_s=100$ if not stated otherwise.  The initial data is given by equiangular distributed particles on a circle with centre $(0.5,0.5)$ and radius $0.005$ in Figure~\subref*{fig:expcoeffnumericsinitial}. In Figures \subref*{fig:expcoeffnumericssmallessmallrcaddconst}--\subref*{fig:expcoeffnumericslargercexpl} the stationary solution for the exponentially decaying force coefficient $f_s^\epsilon$ in the limit $\epsilon\to 0$ is shown. As expected, for small values of $e_s$ and $R_c\in(0,0.5)$, e.g.\ $e_s=10$ as in Figure~\subref*{fig:expcoeffnumericssmallessmallrcaddconst}, the equidistantly distributed particles along the vertical axis are an unstable steady state. In this case, the steady state is given by clusters along the vertical axis and $\Re (\lambda_2(m))\leq 0$ for $m<12$ only. For $R_c=0.5$ the straight vertical line is stable as shown in Figure \subref*{fig:expcoeffnumericssmalleslargercaddconst}. Note that the additional constant in the definition of $f_s^\epsilon$  leads to $f_s^\epsilon(R_c-\epsilon)=f_s^\epsilon(R_c)=0$ and is necessary for the stability of the straight vertical line. In Figure~\subref*{fig:expcoeffnumericssmalleslargerc} we consider $f_s^\epsilon$ without this additional constant, i.e.\ $f_s^\epsilon(|d|)=c\exp(-e_s |d|)$ for $|d|\in[0,R_c-\epsilon]$, where the straight vertical line is clearly unstable and we have $\Re (\lambda_2(m))\leq 0$ for $m<9$ only. If $e_s$ is chosen sufficiently large, e.g. $e_s=100$ as in Figures~\subref*{fig:expcoeffnumericssmallrc} and \subref*{fig:expcoeffnumericslargercexpl}, the straight vertical line appears to be stable even for $R_c<0.5$. An explicit calculation of the eigenvalues for $R_c=0.1$ reveals, however, that $\Re(\lambda_2(m))\leq 0$ for $m< 73723$ only. Note that we obtain stability for a much larger number of modes as in Figures~\subref*{fig:expcoeffnumericssmallessmallrcaddconst} and \subref*{fig:expcoeffnumericssmalleslargerc}. This is also consistent with a straight vertical line as steady state in Figure~\subref*{fig:expcoeffnumericslargercexpl}, while we have clusters as steady states in Figures~\subref*{fig:expcoeffnumericssmallessmallrcaddconst} and \subref*{fig:expcoeffnumericssmalleslargerc}. Further note that $\Re (\lambda_2(73723))=8.3225\cdot 10^{-15}$ and hence it is numerically zero. As discussed in Remark~\ref{rem:expstability} this explains why for $\exp(e_sR_c)\gg 1$, e.g.\ $e_s=100$ and $R_c=0.1$, the straight vertical line appears to be stable.
Finally, we also obtain the straight vertical line as a steady state if we consider exponentially decaying force coefficients $f_l^\epsilon(|d|)=0.13\exp(-100|d|)-0.03\exp(-10|d|)$ instead of $f_l^\epsilon(|d|)=0.1-3|d|$ for $|d|\in[0,R_c-\epsilon]$ in the limit $\epsilon\to 0$ as shown in Figure \subref*{fig:expcoeffnumericslargercexpl}. Note that we also obtain a straight vertical line as stationary solution in  Figures \subref*{fig:expcoeffnumericssmallrc} and \subref*{fig:expcoeffnumericslargercexpl} if $f_s^\epsilon(|d|)=c\exp(-e_s |d|)-c\exp(-e_s (R_c-\epsilon))$ for $|d|\in[0,R_c-\epsilon]$ is replaced by $f_s^\epsilon(|d|)=c\exp(-e_s |d|)$ since $\exp(-e_sR_c)\ll 1$ for $e_s\gg 1$.

In Figure \ref{fig:expcoeffnumericlargedomain} the stationary solution is shown on the domain $[0,3]^2$ instead of the unit square. Here, we consider the same force coefficients  as in Figure \subref*{fig:expcoeffnumericslargercexpl}, i.e.\ exponentially decaying force coefficients along $l$ and $s$. We define the initial data on $[0,3]^2$ by considering the initial data on the unit square, i.e.\ equiangular distributed particles on a circle with centre $(0.5,0.5)$ and radius $0.005$, and extending these initial conditions to $[0,3]^2$ by using the periodic boundary conditions.   As expected we obtain three parallel lines as stationary solution.

For the underlying tensor field $T$ with $s=(0,1)$ and $l=(1,0)$, we have seen that vertical straight patterns are stable. More generally, stripe states along any angle can be obtained by rotating the spatially homogeneous tensor field $T$ appropriately. Examples of rotated stripe patterns are shown in Figure \ref{fig:expcoeffrotatedtensor} where  the vector fields $s=(1,1)/\sqrt{2}, l=(-1,1)/\sqrt{2}$ in Figure \subref*{fig:tensorfielddiag}, $s=(1,2)/\sqrt{5},l=(-2,1)/\sqrt{5}$ in Figure \subref*{fig:tensorfieldtwo} and $s=(1,5)/\sqrt{26},l=(-5,1)/\sqrt{26}$ in Figure \subref*{fig:tensorfieldfive} are considered. Due to the periodicity of the forces, the resulting patterns are also periodic.

Until now, we looked at numerical examples for stable state aligned along a line (or lines). However, the model \eqref{eq:particlemodelperiodic}  is also able to produce  two-dimensional states which can result as an instability of a vertical line. To obtain two-dimensional patterns, we vary the force along $l$. In particular, the force along $l$ has to be less attractive to avoid the concentration along line patterns. In Figure \ref{fig:expcoeffrepulsion}, we vary parameter $e_{l_1}$ in the force coefficient ${f_l^\epsilon(|d|)=0.13\exp(-e_{l_1}|d|)-0.03\exp(-10|d|)}$ for $|d|\in[0,R_c-\epsilon]$. Here, smaller values of $e_{l_1}$ lead to stronger repulsive forces over short distance, resulting in a horizontal spreading of the solution for the tensor field $T$ with $s=(0,1)$ and $l=(1,0)$.

\section*{Acknowledgments}
JAC was partially supported by the EPSRC through grant number EP/P031587/1.
BD has been supported by the Leverhulme Trust research project grant `Novel discretizations for higher-order nonlinear PDE' (RPG-2015-69). LMK was supported by the UK Engineering and Physical Sciences Research Council (EPSRC) grant
EP/L016516/1 and the German Academic Scholarship Foundation (Studienstiftung des Deutschen Volkes). CBS acknowledges support from Leverhulme Trust project on Breaking the non-convexity barrier, EPSRC grant Nr. EP/M00483X/1, the EPSRC Centre Nr. EP/N014588/1, the RISE projects CHiPS and NoMADS, the Cantab Capital Institute for the Mathematics of Information and the Alan Turing Institute.

\bibliographystyle{siamplain}
\bibliography{references}

\end{document}